\pdfoutput=1
\documentclass[11pt]{article}
\usepackage[letterpaper,top=1.0in,bottom=1.05in,left=1.1in,right=1.1in,headheight=22pt,headsep=10pt]{geometry}
\usepackage[T1]{fontenc}
\usepackage[utf8]{inputenc}
\usepackage{lmodern}
\usepackage{amsmath,amssymb,amsthm,mathtools}
\usepackage{microtype}
\usepackage{graphicx,booktabs,array,longtable,xurl,enumitem,needspace,etoolbox,flafter}
\usepackage{tikz}
\usetikzlibrary{arrows.meta,calc,patterns}
\usepackage{xcolor}
\definecolor{Ink}{HTML}{162B44}
\definecolor{LinkBlue}{HTML}{1B4F72}
\definecolor{CiteTeal}{HTML}{116A65}
\definecolor{RuleGray}{HTML}{CBD2D9}
\usepackage{titlesec}
\usepackage{fancyhdr}
\usepackage[unicode,colorlinks=true,linkcolor=LinkBlue,citecolor=CiteTeal,urlcolor=LinkBlue,bookmarksnumbered=true,pdfpagelabels=true]{hyperref}
\usepackage{bookmark}

\newcommand{\R}{\mathbb R}
\newcommand{\norm}[1]{\lVert#1\rVert}
\DeclareMathOperator{\tr}{tr}
\DeclareMathOperator{\aff}{aff}
\DeclareMathOperator{\conv}{conv}
\DeclareMathOperator{\dist}{dist}
\DeclareMathOperator{\argmin}{argmin}

\setlist[itemize]{topsep=4pt,itemsep=2pt,parsep=0pt,leftmargin=1.6em}
\setlist[enumerate]{topsep=4pt,itemsep=3pt,parsep=0pt,leftmargin=1.8em}
\allowdisplaybreaks[1]
\AtBeginEnvironment{thebibliography}{\small\setlength{\itemsep}{3pt plus 1pt}}

\newtheoremstyle{bfgsthm}{8pt plus 2pt minus 1pt}{8pt plus 2pt minus 1pt}{\itshape}{}{\bfseries\color{Ink}}{.}{0.5em}{\thmname{#1}\thmnumber{ #2}\thmnote{ \normalfont(#3)}}
\newtheoremstyle{bfgsdef}{8pt plus 2pt minus 1pt}{8pt plus 2pt minus 1pt}{\normalfont}{}{\bfseries\color{Ink}}{.}{0.5em}{\thmname{#1}\thmnumber{ #2}\thmnote{ \normalfont(#3)}}
\theoremstyle{bfgsthm}
\newtheorem{theorem}{Theorem}
\newtheorem{lemma}{Lemma}
\newtheorem{proposition}{Proposition}
\newtheorem{corollary}{Corollary}
\theoremstyle{bfgsdef}
\newtheorem{definition}{Definition}

\let\amsthmorigqed\qed
\makeatletter
\renewenvironment{proof}[1][\proofname]{\par
  \pushQED{\amsthmorigqed}%
  \normalfont \topsep6\p@\@plus6\p@\relax
  \trivlist
  \item[\hskip\labelsep
        \itshape
    #1\@addpunct{.}]\ignorespaces
}{%
  \popQED\endtrivlist\@endpefalse
}
\makeatother
\renewcommand{\qed}{\qedhere}

\titleclass{\part}{straight}
\titleformat{\part}[display]{\needspace{8\baselineskip}\vspace{10pt}\color{Ink}\hrule height 0.7pt\vspace{8pt}\normalfont\sffamily\bfseries\footnotesize\color{CiteTeal}}{PART \thepart}{3pt}{\normalfont\Large\bfseries\color{Ink}}[\vspace{7pt}{\color{RuleGray}\hrule height 0.5pt}]
\titlespacing*{\part}{0pt}{16pt}{12pt}
\titleformat{\section}{\large\bfseries\color{Ink}}{\thesection.}{0.55em}{}
\titleformat{\subsection}{\normalsize\bfseries\color{Ink}}{\thesubsection}{0.55em}{}
\titleformat{\paragraph}[runin]{\normalsize\bfseries\color{Ink}}{\theparagraph}{0.5em}{}
\titlespacing*{\section}{0pt}{16pt plus 3pt minus 2pt}{7pt plus 1pt}
\titlespacing*{\subsection}{0pt}{12pt plus 2pt minus 1pt}{5pt plus 1pt}
\titlespacing*{\paragraph}{0pt}{7pt plus 2pt minus 1pt}{0.5em}

\renewcommand{\headrulewidth}{0.3pt}
\renewcommand{\headrule}{\hbox to\headwidth{\color{RuleGray}\leaders\hrule height \headrulewidth\hfill}}
\fancypagestyle{plain}{\fancyhf{}\fancyfoot[C]{\footnotesize\thepage}\renewcommand{\headrulewidth}{0pt}}

\newcommand{\kwsep}{\unskip\nobreak\hspace{0.55em}\textperiodcentered\hspace{0.55em}\ignorespaces}
\newcommand{\keywords}[1]{\par\vspace{8pt}\noindent{\footnotesize\raggedright\textbf{\color{Ink}Keywords.}\hspace{0.5em}#1\par}}
\newcommand{\subclass}[1]{\par\vspace{3pt}\noindent{\footnotesize\raggedright\textbf{\color{Ink}MSC 2020.}\hspace{0.5em}#1\par}}

\hypersetup{
  pdftitle={BFGS on Convex Nonsmooth Functions: Convergence and Failure},
  pdfauthor={Yuchen Xie},
  pdfsubject={Unmodified full BFGS; global and local convergence; deterministic failure},
  pdfkeywords={BFGS, convex nonsmooth optimization, global convergence, local convergence, counterexample}
}

\title{\vspace{-1.2em}{\LARGE\bfseries\color{Ink} BFGS on Convex Nonsmooth Functions:\\[3pt] Convergence and Failure}}
\author{\normalsize\textbf{Yuchen Xie}\\[2pt]
{\small Independent Researcher}\\[1pt]
{\small\href{mailto:yxie@u.northwestern.edu}{\texttt{yxie@u.northwestern.edu}}}}
\date{\small November 2026}

\begin{document}
\maketitle
\thispagestyle{plain}

\begin{abstract}
\noindent Unmodified full BFGS with an Armijo--weak-Wolfe line search routinely succeeds
on nonsmooth convex functions even as gradient jumps drive its inverse-Hessian
approximations toward singularity and invalidate the smooth secant bound used by
classical trace--determinant potentials. We organize the analysis into two
independent parts sharing a common algorithmic setup. In Part~I, we estimate the
scalar sequence $(g_k^TH_kg_k)^{-1}$ in two ways: along any bounded-gradient run
with vanishing predicted decreases, an adjacent-degree hierarchy for the
principal minors of $H_k$ forces its geometric and arithmetic means to diverge
without smoothness or conditioning assumptions, whereas on structured nonsmooth
objectives a persistent optimality gap keeps that average bounded via a
comparison energy $w_k^TB_kw_k$ that filters out orthogonal gradient jumps and
bounds field rotation by the decrease in $F$, proving convergence for every
infinite admissible run (with differentiability at accepted iterates) without
assuming iterate convergence or spectral asymptotics. Static radial energies settle
whole-sequence convergence on the open benchmark $F(u,v)=u^2+|v|$ of Guo and
Lewis (2018) within an $n$-dimensional smooth-plus-homogeneous theorem, prove
global convergence for all coercive positively homogeneous convex functions
(including all norms and coercive linear maxima), and yield matrix-independent
local sublevel basins around a unique minimizer for finite polyhedral $C^{1,1}$
partitions such as $\ell_1$-regularized losses; moving comparison fields then
resolve nondegenerate planar smooth maxima with curved switching curves on which
radial secants blow up. Conversely, in Part~II we construct a coercive
four-piece polyhedral function in $\R^2$ on which, at a corner where zero lies
in the affine hull of active gradients outside their convex hull, unmodified
BFGS with $H_0=I$ and the unit-first bracketing search enters a self-similar
three-step cycle and converges to a nonstationary point. Two distinct local
constructions place failing runs from $H_0=I$ arbitrarily near the unique
minimizer: one on a maximum of four convex $C^{2,1}$ branches, and another on a
maximum of four uniformly strongly convex $C^{1,1}$ branches.
\subclass{90C30 \kwsep 90C53 \kwsep 49J52}
\end{abstract}

\newpage
{\small\linespread{0.97}\selectfont
\makeatletter
\renewcommand{\l@part}[2]{%
  \ifnum \c@tocdepth >-2\relax
    \addpenalty{-\@highpenalty}%
    \addvspace{8pt \@plus 2pt}%
    \setlength\@tempdima{2.2em}%
    \begingroup
      \parindent \z@ \rightskip \@pnumwidth
      \parfillskip -\@pnumwidth
      {\leavevmode\normalsize\bfseries\color{Ink}Part #1\hfil \hb@xt@\@pnumwidth{\hss #2}}\par
      \nobreak
      \global\@nobreaktrue
      \everypar{\global\@nobreakfalse\everypar{}}%
    \endgroup
  \fi}
\renewcommand{\l@section}[2]{%
  \ifnum \c@tocdepth >\z@
    \addpenalty\@secpenalty
    \addvspace{3.5pt \@plus 1pt}%
    \setlength\@tempdima{1.5em}%
    \begingroup
      \parindent \z@ \rightskip \@pnumwidth
      \parfillskip -\@pnumwidth
      \leavevmode \bfseries
      \advance\leftskip\@tempdima
      \hskip -\leftskip
      #1\nobreak\hfil \nobreak\hb@xt@\@pnumwidth{\hss #2}\par
    \endgroup
  \fi}
\makeatother
\tableofcontents}
\newpage

\section{Introduction, Algorithmic Setup, and Main Results}\label{art:introduction}

The Broyden--Fletcher--Goldfarb--Shanno (BFGS) method builds curvature models
from gradient differences, yet in computational practice an unmodified BFGS
implementation with an inexact Armijo--weak-Wolfe line search routinely succeeds
on nonsmooth convex problems~\cite{LO}. Classical global convergence proofs for
smooth convex optimization, from Powell~\cite{FND_P76} and Byrd and
Nocedal~\cite{FND_BN89} to Nocedal's survey~\cite{FND_N92}, track the matrix
potential $\psi(B_k)=\tr(B_k)-\ln\det(B_k)$ of the Hessian approximations
$B_k=H_k^{-1}$. That argument relies on a smooth secant bound
$\norm{y_k}^2/(y_k^Ts_k)\le M$ to prevent $\tr(B_k)$ from growing faster than
the determinant term can compensate. On a nonsmooth objective, however, gradients can jump by $O(1)$ across kinks as
the step $s_k=x_{k+1}-x_k$ shrinks: the ratio $\norm{y_k}^2/(y_k^Ts_k)$ becomes
unbounded, $\tr(B_k)$ can diverge, and the inverse-Hessian approximation $H_k$
degenerates toward singularity, whether through unbounded ill-conditioning as
it crushes directions normal to a kink~\cite{LO} or through uniform contraction
$H_k\to0$.

Existing theory leaves a wide gap between empirical success and rigorous
guarantees. In one dimension, Lewis and Overton~\cite[Section~5.1]{LO} analyzed
$|x|$ directly, and on an unbounded two-piece affine model in $\R^n$, Xie and
W\"achter~\cite[Theorem~IV.1 and Remark~IV.1]{XW} showed that BFGS finds an unbounded
descent ray in finitely many steps. In general dimension, Guo and
Lewis~\cite[Theorems~3.2, 4.1, and~4.3]{GL18} extended Powell's method of proof
by replacing the nonsmooth objective with smooth convex surrogates on regions
visited by the iterates. That technique proves convergence on the Euclidean
norm $\norm{x}_2$ in $\R^n$~\cite[Corollary~4.2]{GL18}, where the kink is an
isolated point, but on the canonical mixed smooth--nonsmooth benchmark
\[
 F(u,v)=u^2+|v|
\]
it yields only $\liminf_{k\to\infty}|v_k|=0$, showing that a subsequence
approaches the nonsmooth axis $v=0$ without giving a corresponding bound on
$u_k$~\cite[Corollary~4.4]{GL18}.
Whether every full-BFGS sequence on $u^2+|v|$ converges to the minimizer
$(0,0)$ was left explicitly open by Guo and Lewis~\cite[p.~1309]{GL18} and
highlighted again as unresolved by Gebken~\cite{BG}, whose local criticality
results for piecewise smooth selections instead assume \emph{a priori} iterate
convergence together with prescribed active-set and spectral
asymptotics~\cite[Theorem~3.1 and~(B1)]{BG}.

Even less was settled on the negative side for convex problems. Known divergence
or non-Clarke-stationary examples for BFGS either rely on nonconvex objectives,
as in Mascarenhas~\cite{NEG_M04,NEG_M07}, Diao~\cite{NEG_DIAO26}, and
Lai~\cite{NEG_LAI26} (while Lewis and Zhang~\cite{NEG_LZ15} construct a
nonconvex example converging to a Clarke-stationary point that admits directions
of linear descent), or alter the algorithm itself, whether by imposing exact
line searches and specialized subgradient choices as in Lewis and
Overton~\cite{NEG_LO08} or by replacing full BFGS with a scaled memoryless
update as in Asl and Overton~\cite{NEG_AO21}. While Lai~\cite{NEG_LAI26} recently disproved
the nonconvex semialgebraic stationarity conjecture of Lewis and Overton's
Challenge~7.1~\cite{LO} using a nonconvex homogeneous cycle with negative
transverse curvature, it remained open whether unmodified full BFGS with a
standard inexact Armijo--weak-Wolfe line search can converge to a nonstationary
point on a coercive \emph{convex} function, either globally or from starting
points arbitrarily close to the unique minimizer of a smooth-branch maximum. We
construct convex counterexamples of both kinds.

This paper develops two complementary, self-contained analyses of unmodified
full BFGS on convex nonsmooth functions, organized as Part~\ref{part:convergence}
and Part~\ref{part:failure} following a shared algorithmic setup below. In
Part~\ref{part:convergence}, we establish sufficient geometric conditions under
which every infinite admissible run (for arbitrary $H_0\succ0$ and fixed Wolfe
parameters $0<c_1<c_2<1$) converges to minimality, settling $u^2+|v|$ alongside
broad global and local classes in arbitrary dimension and curved switching
interfaces in the plane. In Part~\ref{part:failure}, we fix one standard
deterministic line-search implementation with $H_0=I$ and construct global and
local convex objectives on which the iterates converge to nonstationary points.
Because the principal-minor growth law of Part~\ref{part:convergence} continues
to hold even on these nonstationary runs, convergence vs.\ failure turns
directly on whether objective geometry supplies a comparison direction that
filters out gradient jumps.

\subsection{Algorithmic setup and basic invariants}\label{art:setup}

To state the convergence theory and the counterexamples on a common footing, we
first record the unmodified full-BFGS protocol and its elementary one-step
invariants; after this subsection, Part~\ref{part:convergence} and
Part~\ref{part:failure} can be read in either order.

\begin{definition}[Admissible runs]\label{art:admissible}
Let $F:\R^n\to\R$ be finite and convex, with $n\ge1$.
An \emph{admissible run} is an infinite sequence of differentiable iterates
$x_k$, nonzero gradients $g_k=\nabla F(x_k)$, and inverse-Hessian approximations
$H_k$ generated from $H_0\succ0$ by step lengths $\alpha_k>0$ and the
unmodified full-BFGS update
\begin{align}
 s_k&=-\alpha_kH_kg_k,\qquad x_{k+1}=x_k+s_k,
       \qquad y_k=g_{k+1}-g_k,\notag\\
 H_{k+1}&=\left(I-\frac{s_ky_k^T}{t_k}\right)
 H_k\left(I-\frac{y_ks_k^T}{t_k}\right)
       +\frac{s_ks_k^T}{t_k},\qquad t_k=y_k^Ts_k,\label{art:update}
\end{align}
satisfying the Armijo and weak-Wolfe conditions for fixed constants
$0<c_1<c_2<1$:
\begin{equation}\label{art:aw}
 F(x_{k+1})\le F(x_k)+c_1g_k^Ts_k,
 \qquad g_{k+1}^Ts_k\ge c_2g_k^Ts_k.
\end{equation}
\end{definition}
This is the standard nonsmooth full-BFGS protocol of Lewis and
Overton~\cite{LO}; see also~\cite{XW}. Encountering $\nabla F(x_k)=0$ is
exact finite termination at a global minimizer.

Throughout, we write
\[
 \begin{gathered}
 B_k=H_k^{-1},\qquad q_k=g_k^TH_kg_k,\qquad
 d_k=-g_k^Ts_k=\alpha_kq_k,\\
 \Delta_k=F(x_k)-F(x_{k+1}),\qquad \omega=1-c_2>0.
 \end{gathered}
\]
Here $q_k$ is the squared $H_k$-norm of the gradient and $d_k$ is the
first-order predicted decrease along $s_k$.

\begin{lemma}[Basic BFGS invariants]\label{art:run-invariants}
Every admissible run satisfies $H_k\succ0$, $q_k,d_k,t_k>0$, and
\begin{equation}\label{art:decrease-bounds}
 c_1d_k\le\Delta_k\le d_k,\qquad t_k\ge\omega d_k.
\end{equation}
Moreover, the Hessian approximations $B_k=H_k^{-1}$ evolve by
\begin{equation}\label{art:hessian-update}
 B_{k+1}=B_k-\frac{g_kg_k^T}{q_k}+\frac{y_ky_k^T}{t_k},
 \qquad \widehat B_k:=B_k-\frac{g_kg_k^T}{q_k}\succeq0,
 \qquad \widehat B_ks_k=0.
\end{equation}
If $F(x_k)$ is bounded below (hence converges to a finite limit), then
$\sum_{k=0}^\infty d_k<\infty$ and in particular $d_k\to0$.
\end{lemma}
\begin{proof}
Given $H_k\succ0$ and $g_k\ne0$, we have $q_k>0$, $d_k=\alpha_kq_k>0$, and
$t_k=(g_{k+1}-g_k)^Ts_k\ge(1-c_2)d_k=\omega d_k>0$ by weak Wolfe. For any
$v\in\R^n$, $v^TH_{k+1}v$ in~\eqref{art:update} is the sum of two nonnegative
terms; if both vanish, then $s_k^Tv=0$ and consequently $(I-y_ks_k^T/t_k)v=v=0$,
so $H_{k+1}\succ0$ by induction. In~\eqref{art:decrease-bounds}, the lower bound
is Armijo and the upper bound follows from convexity,
$F(x_{k+1})\ge F(x_k)+g_k^Ts_k=F(x_k)-d_k$. Summing $c_1d_k\le\Delta_k$ yields
$\sum_kd_k<\infty$ whenever $F(x_k)$ has a finite limit.

Multiplying the right-hand side of~\eqref{art:hessian-update} by~\eqref{art:update}
and using $B_ks_k=-\alpha_kg_k$ and $s_k^Ty_k=t_k$ gives the classical direct
BFGS update~\cite[Section~5]{FND_N92}. Finally, Cauchy--Schwarz in the
$B_k$-inner product gives $(g_k^Tv)^2=(s_k^TB_kv/\alpha_k)^2\le q_kv^TB_kv$, so
$\widehat B_k\succeq0$, while $B_ks_k=-\alpha_kg_k$ and $g_k^Ts_k=-\alpha_kq_k$
give $\widehat B_ks_k=0$.\qed
\end{proof}

Quantifying over all admissible runs in Part~\ref{part:convergence} separates
convergence of valid trajectories from the existence of a differentiable
weak-Wolfe step along a given search ray. Appendix~\ref{art:hypothesis-examples}
gives both a general countability criterion under which infinite admissible runs
exist for every choice of differentiable starting point, positive-definite
initial matrix $H_0\succ0$, and Wolfe parameters $0<c_1<c_2<1$
(Lemma~\ref{art:examples-run-existence}, covering our curved quadratic maxima)
and a coercive three-piece cone where a differentiable starting point has a
search ray containing no differentiable weak-Wolfe step
(Proposition~\ref{art:examples-no-step}). Conversely, the counterexamples in
Part~\ref{part:failure} leave no freedom in step selection: they execute the
concrete deterministic unit-first doubling/bisection algorithm of Lewis and
Overton~\cite[Algorithm~4.6]{LO}, with every accepted and rejected trial point
certified to lie in the open smooth interior of a single active piece.

\subsection{Core analytical ideas: conflicting curvature estimates and self-similar cycles}\label{art:geometry}

Why does full BFGS continue to make progress on $F(u,v)=u^2+|v|$ even as its
matrices degenerate, and what geometric change allows a run to converge to a
nonstationary corner? Both phenomena turn on whether the objective supplies a
comparison direction that reconciles the BFGS matrix update with nonsmooth
gradient jumps.

\paragraph{1. A universal principal-minor hierarchy and reciprocal growth.}
Coupling $\tr(B_k)$ directly to $\det(B_k)$ fails as soon as a single
eigenvalue of $B_k$ explodes across a kink. Instead, we track the entire tower
of elementary symmetric polynomials $e_1(H_k),\ldots,e_n(H_k)$, namely the sums
of $j\times j$ principal minors of the \emph{inverse} matrices $H_k$. We prove
that every full-BFGS step satisfies an adjacent-degree inequality
(Theorem~\ref{art:minor-hierarchy}) which, after multiplying by
$q_k=g_k^TH_kg_k$ under weak Wolfe, bounds $q_ke_j(H_{k+1})$ by a multiple of
$\norm{y_k}^2e_{j+1}(H_k)$ plus a vanishing multiple $O(d_k)e_j(H_k)$, where
$d_k=-g_k^Ts_k$ is the predicted decrease. The gradient jump
$y_k=g_{k+1}-g_k$ enters only through $\norm{y_k}^2$, which stays bounded on
any compact sublevel of a finite convex function regardless of nonsmoothness.
As $d_k\to0$, a geometrically weighted Lyapunov combination of
$e_1(H_k),\ldots,e_n(H_k)$ forces both the geometric and arithmetic means of
$q_k^{-1}=(g_k^TH_kg_k)^{-1}$ to diverge to $+\infty$
(Theorem~\ref{art:reciprocal-growth}):
\[
 \frac1N\sum_{k<N}\frac{1}{g_k^TH_kg_k}\longrightarrow+\infty.
\]
This law holds along every bounded-gradient BFGS run with $d_k\to0$,
without any smoothness, strong convexity, or condition-number assumption. In
fact, it continues to hold along our nonstationary counterexamples, where
$H_k\to0$ drives $q_k\to0$ (exponentially on the homogeneous seed cycle).

\paragraph{2. Directional and transported comparison energies.}
By Cauchy--Schwarz in the $H_k$-inner product, any vector $w_k$ satisfies
$w_k^TB_kw_k\ge(g_k^Tw_k)^2/q_k$. Under bounded gradients and $d_k\to0$, the
divergence $\frac1N\sum_{k<N}q_k^{-1}\to+\infty$ therefore forces the average of
$w_k^TB_kw_k$ to diverge along any field maintaining $|g_k^Tw_k|\ge\sqrt m>0$.
Testing $B_k=H_k^{-1}$ solely along such a vector $w_k$ (or a
positive-semidefinite tensor $W(x_k)$) avoids bounding $\tr(B_k)$ altogether.
In the Hessian update~\eqref{art:hessian-update},
the rank-one downdate $\widehat B_k=B_k-g_kg_k^T/q_k\succeq0$ annihilates the
step direction: $\widehat B_ks_k=0$. Evaluating $E_k=r_k^TB_kr_k$ along the
radial vector $r_k=x_k-x_*$ therefore makes the shift $r_{k+1}=r_k+s_k$
algebraically free and gives the exact one-step identity
(Lemma~\ref{art:reference-energy})
\[
 E_{k+1}=E_k-\frac{(g_k^Tr_k)^2}{q_k}+\frac{(y_k^Tr_{k+1})^2}{t_k}.
\]
At value gap $\eta>0$, the downdate subtracts at least $\eta^2/q_k$ before the
radial secant term is added, and the gradient jump enters only through its
radial projection $y_k^Tr_{k+1}$, not through $\norm{y_k}$. If
$(y_k^Tr_{k+1})^2/t_k\le C$, summing over $k<N$ and using $E_N\ge0$ forces
$\eta^2\sum_{k<N}q_k^{-1}\le E_0+CN$, a linear upper bound that directly
contradicts $\frac1N\sum_{k<N}q_k^{-1}\to+\infty$. Notably, this budget imposes
no upper bound on directions orthogonal to $r_k$; only the nonnegativity
$E_N\ge0$ is used.

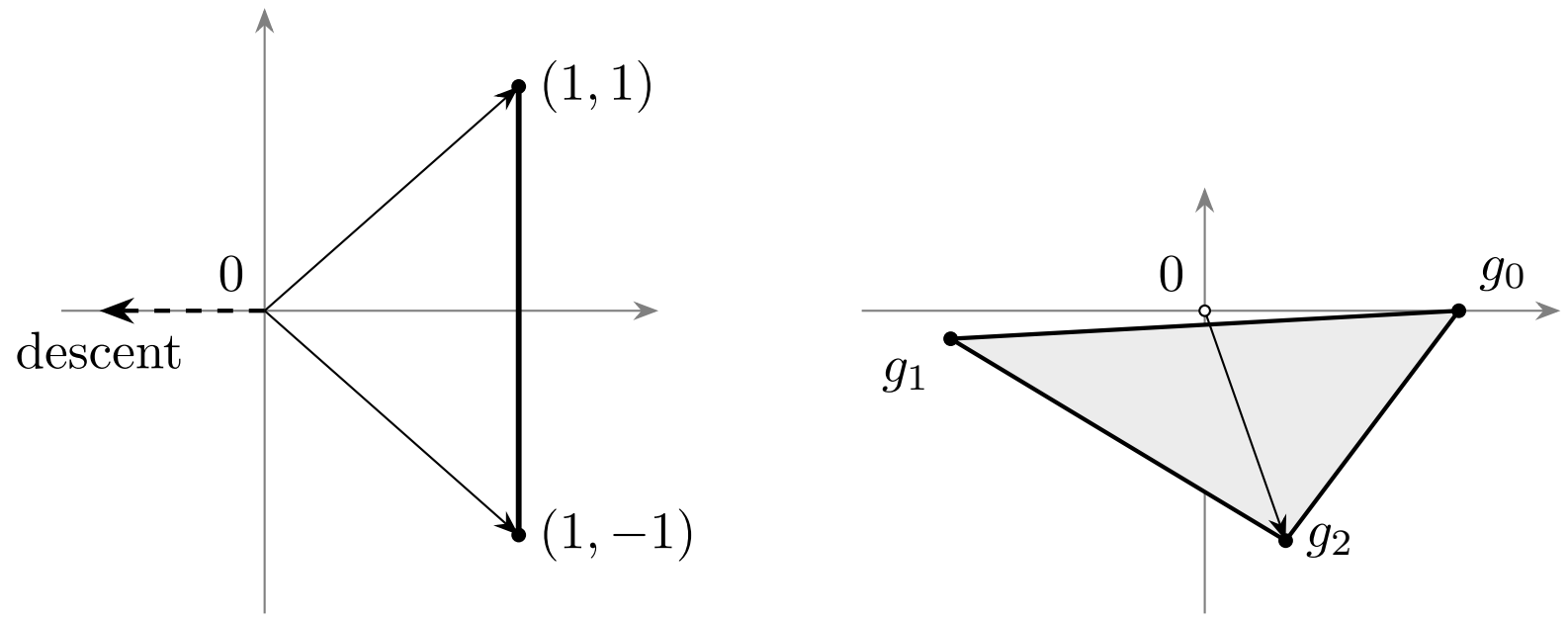
\begin{figure}[tbp]
\centering
\begin{tikzpicture}[x=1.7cm,y=1.5cm,>=Stealth,font=\small]
\begin{scope}[shift={(-1.9,0)}]
\draw[->,gray] (-.8,0)--(1.55,0);\draw[->,gray] (0,-1.35)--(0,1.35);
\draw[line width=1.1pt] (1,-1)--(1,1);
\draw[->] (0,0)--(1,1);\draw[->] (0,0)--(1,-1);
\fill (1,1) circle (1.4pt) node[right] {$(1,1)$};
\fill (1,-1) circle (1.4pt) node[right] {$(1,-1)$};
\draw[->,thick,dashed] (0,0)--(-.65,0) node[below] {descent};
\node at (.4,-1.65) {(a) $u^2+|v|$, $u=1/2$};
\node[above left] at (0,0) {$0$};
\end{scope}
\begin{scope}[shift={(1.8,0)}]
\draw[->,gray] (-1.35,0)--(1.4,0);\draw[->,gray] (0,-1.35)--(0,.55);
\fill[gray!15] (1,0)--(-1,-.125)--(.3184814276,-1.0252321197)--cycle;
\draw[thick] (1,0)--(-1,-.125)--(.3184814276,-1.0252321197)--cycle;
\draw[->] (0,0)--(.3184814276,-1.0252321197);
\fill (1,0) circle (1.4pt) node[above right] {$g_0$};
\fill (-1,-.125) circle (1.4pt) node[below left] {$g_1$};
\fill (.3184814276,-1.0252321197) circle (1.4pt) node[right] {$g_2$};
\fill[white] (0,0) circle (1pt);\draw (0,0) circle (1pt);
\node[above left] at (0,0) {$0$};
\node at (0,-1.65) {(b) Nonstationary three-slope corner};
\end{scope}
\end{tikzpicture}
\caption{Subdifferential geometry in gradient space. (a)~On $F(u,v)=u^2+|v|$ at
$u=1/2$, the limiting gradients $(1,1)$ and $(1,-1)$ across $v=0$ jump
vertically, so the horizontal vector $(-u,0)$ pairs identically ($-2u^2<0$)
with both. (b)~The nonstationary three-slope corner of Section~\ref{art:failure}:
$0$ lies outside $\conv\{g_0,g_1,g_2\}$ (at distance $1/\sqrt{257}$ above the
edge $[g_1,g_0]$, so $(-1,16)$ is a strict common descent direction), yet
$\aff\{g_0,g_1,g_2\}=\R^2$, destroying any jump-canceling comparison direction.}
\label{art:gradient-geometry}
\end{figure}

Figure~\ref{art:gradient-geometry}(a) illustrates why this directional filtering
succeeds on $F(u,v)=u^2+|v|$. Approaching the nonsmooth line $v=0$ need not
make $u^2$ small, because the gradient jumps vertically between $(2u,1)$ and
$(2u,-1)$. For $u\ne0$, both gradients have the same negative pairing $-2u^2$
with the horizontal direction $(-u,0)$, and along any inward ray toward $(0,0)$
the nonsmooth term $|v|$ scales linearly with zero Taylor remainder: for $v\ne0$
and $0\le\tau\le1$,
\begin{equation}\label{art:mixed-inward-remainder}
 F((1-\tau)(u,v))-F(u,v)+\tau\nabla F(u,v)^T(u,v)=\tau^2u^2.
\end{equation}
Only the smooth component $u^2$ contributes a quadratic remainder. A general
Bregman-duality estimate (Theorem~\ref{art:radial-equivalence}) converts any
such one-point inward quadratic bound into the two-point secant bound
$(y_k^Tr_{k+1})^2\le C_\eta t_k$ across arbitrary kinks. When kinks instead lie
along \emph{curved} switching interfaces, radial rays cross them transversely
and radial secant bounds can fail. Even there, the identity
$\widehat B_ks_k=0$ implies that rotating a comparison field $w(x)$ along
level-set contours costs zero energy in the step direction $s_k$, while the
decrease in $F$ bounds the transverse multiplicative cost $1+Ld_k$ through the
Armijo summability $\sum_kd_k<\infty$ (Theorem~\ref{art:certificate}).

\paragraph{3. Affine vs.\ convex subdifferential balance and self-similar cycles.}
Figure~\ref{art:gradient-geometry}(b) reveals the exact geometric obstruction
where comparison energies break down. For a nonempty finite set of active
slopes $\{a_i\}$, ordinary nonstationarity means the existence of a common
descent direction ($a_i^Tw<0$ for all $i$), which by closest-point separation
holds if and only if $0\notin\conv\{a_i\}$. By contrast, screening out all
leading gradient jumps $(a_i-a_j)^Tw=0$ while maintaining strict visibility
requires an \emph{equal-rate} descent direction ($a_i^Tw=\kappa<0$ for all $i$),
which exists if and only if $0\notin\aff\{a_i\}$. In $\R^2$, three affinely
independent slopes $g_0=(1,0)^T$, $g_1=(-1,-1/8)^T$, and $g_2=(c_\circ,d_\circ)^T$
have $\aff\{g_0,g_1,g_2\}=\R^2$ while $0\notin\conv\{g_0,g_1,g_2\}$: the vector
$(-1,16)^T$ decreases all three linear pieces, yet every nonzero vector
$w\in\R^2$ sees an $O(1)$ jump across at least one switching step. We show that
$(c_\circ,d_\circ)$ can be chosen algebraically so that three steps of unmodified
BFGS from $H_0=I$ under Lewis--Overton's standard bracketing line search return
both the matrix and the iterate to an exact scalar multiple $\rho\approx0.0602$
of their initial state: $(z_3,H_3)=(\rho z_0,\rho I)$. Homogeneity then locks
BFGS into an infinite geometric contraction whose summable steps converge to the
nonstationary origin.

\subsection{Overview of main results and independent readability}

\paragraph{Part~\ref{part:convergence}: Curvature budgets, comparison energies, and convergence (Sections~\ref{art:energy}--\ref{art:local}).}
Section~\ref{art:energy} establishes the adjacent-degree principal-minor
hierarchy for a single inverse-BFGS update (Theorem~\ref{art:minor-hierarchy})
and derives the divergence of the geometric and arithmetic means of
$(g_k^TH_kg_k)^{-1}$ whenever gradients are bounded and $d_k\to0$
(Theorem~\ref{art:reciprocal-growth}). Notably, although
Theorem~\ref{art:reciprocal-growth} is stated for admissible runs, its proof
uses no convexity of $F$: both theorems are purely algebraic properties of the
BFGS update under $y_k^Ts_k\ge(1-c_2)d_k$, bounded gradients, and $d_k\to0$.

\begin{table}[tbp]
\centering
\begin{tabular}{@{}>{\raggedright\arraybackslash}p{.41\textwidth}>{\raggedright\arraybackslash}p{.23\textwidth}>{\raggedright\arraybackslash}p{.28\textwidth}@{}}
\toprule
Objective class & Basin of attraction & Guarantee\\
\midrule
Centered smooth-plus-homogeneous $F(x)=\psi(x)+\gamma(x-x_*)$, Theorem~\ref{art:centered}
& Any compact sublevel inside domain of representation
& $F(x_k)\to F_*$, $\dist(x_k,\argmin F)\to0$; $x_k\to x_*$ if unique\\\addlinespace
Coercive degree-one homogeneous (all norms, gauges, coercive $\max_ia_i^Tx$), Theorem~\ref{art:homogeneous}
& Global (every admissible $x_0$)
& $F(x_k)\to0$ and $x_k\to0$\\\addlinespace
Finite polyhedral $C^{1,1}$ partition with unique $x_*$ (including smooth loss $+\;\lambda\norm{x}_1$), Theorem~\ref{art:polyhedral}
& Sublevel $\{F\le F_*+\varepsilon\}$ independent of $H_0$
& $F(x_k)\to F_*$ and $x_k\to x_*$\\\addlinespace
Nondegenerate planar smooth $\max_if_i$ (unique $x_*$, strict contours), Corollary~\ref{art:planar-max}
& Sublevel $\{F<F_*+\varepsilon\}$ independent of $H_0$
& $F(x_k)\to F_*$ and $x_k\to x_*$\\
\bottomrule
\end{tabular}
\caption{Summary of main positive convergence guarantees in $\R^n$ (first three
rows) and $\R^2$ (fourth row). Every guarantee holds for all infinite
admissible full-BFGS runs with arbitrary $H_0\succ0$ and any fixed
$0<c_1<c_2<1$. For $F(u,v)=u^2+|v|$ the first row applies on all of $\R^2$.}
\label{art:positive-outcomes}
\end{table}

Sections~\ref{art:global} and~\ref{art:local} progress from static radial
vectors to moving comparison fields across the four regimes in
Table~\ref{art:positive-outcomes}:
\begin{enumerate}
\item \emph{Smooth-plus-homogeneous decompositions and $u^2+|v|$ (Section~\ref{art:global}).}
Via the static radial identity (Lemma~\ref{art:reference-energy}) and Bregman
inward-to-secant duality (Theorem~\ref{art:radial-equivalence}),
Theorem~\ref{art:centered} covers $F(x)=\psi(x)+\gamma(x-x_*)$ in $\R^n$ for
convex $\psi\in C^{1,1}_{\rm loc}$ and degree-one homogeneous convex $\gamma$,
even when neither summand is minimized at $x_*$ alone, settling whole-sequence
convergence $(u_k,v_k)\to(0,0)$ for every admissible run on $F(u,v)=u^2+|v|$
(Corollary~\ref{art:mixed-benchmark}).
\item \emph{Coercive positively homogeneous functions in $\R^n$ (Section~\ref{art:global}).}
Theorem~\ref{art:homogeneous} and Corollary~\ref{art:homogeneous-examples} yield
global convergence $x_k\to0$ for all norms, bounded-unit-ball gauges, and
coercive linear maxima $\max_{1\le i\le m}a_i^Tx$; for $n\ge2$ this includes
$\norm{x}_1$ and $\norm{x}_\infty$, whose nonsmooth loci contain coordinate
hyperplanes and diagonal facets, respectively, rather than an isolated point.
\item \emph{Polyhedral partitions and composite regularizers in $\R^n$ (Section~\ref{art:global}).}
Near a unique minimizer $x_*$, only polyhedral cells containing $x_*$ meet a
small ball, so inward rays stay within single cells; sublevel localization
(Lemma~\ref{art:localization}) gives an $H_0$-independent basin for finite
polyhedral $C^{1,1}$ partitions, PLQ penalties, Lasso-type $\ell_1$ losses, and
nonpolyhedral Euclidean group norms (Theorem~\ref{art:polyhedral},
Corollaries~\ref{art:plq} and~\ref{art:local-centered}).
\item \emph{Moving comparison fields and curved planar interfaces (Section~\ref{art:local}).}
When switching curves bend, radial secant bounds can fail on every neighborhood
of the minimizer, as shown on an explicit three-quadratic maximum
(Theorem~\ref{art:nonradial-example}). Motivated by this radial obstruction,
Theorem~\ref{art:certificate} and Lemma~\ref{art:tail-comparison} introduce
transported band certificates $W_\eta(x)$ and local vector fields $w(x)$ without
condition-number bounds, smooth surrogates, or prior iterate or spectral
assumptions. In $\R^2$, a normal-angle contour-velocity field whose angular
derivative cancels normal components forces at least one cluster point into the
affine-critical set $\{p:0\in\aff\partial F(p)\}$ (Theorem~\ref{art:planar-cluster}).
When the active gradients at $x_*$ are affinely independent with strictly
positive optimal weights, every affine-critical point near $x_*$ is optimal,
matching the multiplier nondegeneracy condition of nonsmooth Newton and bundle
theory~\cite{LW} and giving local convergence
(Corollaries~\ref{art:planar-convergence} and~\ref{art:planar-max}) with an
explicit basin $\{F<1/100\}$ on the three-quadratic example.
\end{enumerate}

\paragraph{Part~\ref{part:failure}: Self-similar cycles and nonstationary convergence (Sections~\ref{art:failure}--\ref{art:failure-local}).}
The structural hypotheses of Part~\ref{part:convergence} cannot be replaced by
convexity or branch smoothness alone:
\begin{enumerate}
\item \emph{Global failure on a four-piece coercive polyhedral function (Section~\ref{art:failure}).}
Capping the three-slope cone with a fourth affine piece $v-17$ yields a coercive
polyhedral function $\overline G(u,v)=\max\{u,-u-v/8,c_\circ u+d_\circ v,v-17\}$
with unique minimizer $(-1,16)$ and minimum $-1$, on which unmodified BFGS from
$z_0$ with $H_0=I$ and Lewis--Overton's unit-first bracketing search converges
to the nonstationary origin with value $0$ for every $0<c_1\le1/3$, $c_1<c_2<1$
(Theorem~\ref{art:failure_polyhedral}). Homogeneity in
Theorem~\ref{art:homogeneous} is therefore essential: nonzero affine offsets can
destroy the centered representation and separate a nonstationary corner from
the global minimizer.
\item \emph{Two local counterexamples near the minimizer of smooth-branch and
strongly convex maxima (Section~\ref{art:failure-local}).}
Although Theorem~\ref{art:polyhedral} guarantees local convergence on any finite
polyhedral objective with a unique minimizer, finitely many nonlinear smooth
branches can replicate the three-slope corner at shrinking scales
$\bar x_\ell\to0$ around the unique minimizer of a single fixed function. Using
the scale-covariance relation that $F(\bar x+\lambda z)=C+\sigma\mathcal G(z)$
preserves $H_0=I$ when $\sigma=\lambda^2$, Theorem~\ref{art:local-failure}
constructs a coercive maximum of four globally convex $C^{2,1}(\R^2)$ branches
($0\preceq\nabla^2f_i\preceq64I$) with failing starts $x_0^{(\ell)}\to0$ from
$H_0=I$ whose query trajectories shrink to $0$, independently of $(c_1,c_2)$ for
$0<c_1\le1/3$, $c_1<c_2<1$. Imposing uniform strong convexity breaks exact
homogeneity; in a separate construction, Theorem~\ref{art:local-failure-lower}
extracts perturbed failing orbits from the local stable manifold of a
scale-normalized three-step BFGS return map and embeds them into a maximum of
four uniformly strongly convex $C^{1,1}(\R^2)$ branches satisfying
$\underline\mu I\preceq\nabla^2F\preceq\overline L I$ wherever twice
differentiable.
\end{enumerate}

\paragraph{Reader's guide and formal verification.}
Part~\ref{part:convergence} (Sections~\ref{art:energy}--\ref{art:local}, with
supporting proofs in Appendices~\ref{art:matrix-algebra}--\ref{art:hypothesis-examples})
and Part~\ref{part:failure} (Sections~\ref{art:failure}--\ref{art:failure-local},
with rational interval and return-map proofs in
Appendices~\ref{art:seed_certificate}--\ref{art:lower_curvature}) use
independent proof techniques after Section~\ref{art:setup} and may be read in
either order. Section~\ref{art:verification} reunites both parts to map the
exact geometric boundaries between convergence and failure and pose four open
problems. In addition to the self-contained mathematical proofs given here, all
results and rational interval certificates of this paper and the accompanying
companion monograph have been formalized and kernel-checked from first
principles in Lean~4 on top of \textsf{Mathlib} using only the standard axioms
of classical mathematics (Section~\ref{art:verification}; archived at
\mbox{\url{https://doi.org/10.5281/zenodo.24648859}}).

\part{Curvature Budgets, Comparison Energies, and Convergence}\label{part:convergence}
\section{Universal Principal-Minor Hierarchy and Reciprocal Curvature Growth}\label{art:energy}

In smooth optimization, convergence analyses of BFGS rely on Lipschitz
continuity of the gradient to bound the full secant ratio
$\norm{y_k}^2/(y_k^Ts_k)$, keeping the trace of the Hessian approximation
$B_k=H_k^{-1}$ under control~\cite{FND_P76,FND_BN89,FND_N92}. When gradients jump
across nonsmooth kinks, $\norm{y_k}^2/(y_k^Ts_k)$ can diverge to $+\infty$ even
along steps approaching a minimizer. Part~\ref{part:convergence} analyzes
unmodified full BFGS directly in this singular regime by estimating the partial
sums of the same scalar sequence $(g_k^TH_kg_k)^{-1}$ in two opposite ways:
here we prove a universal lower-growth law forcing its arithmetic and geometric
means to diverge, and in Sections~\ref{art:global}--\ref{art:local} we pair that
divergence with static radial and moving comparison energies that bound those
same averages linearly whenever an optimality gap persists.

For a symmetric $n\times n$ matrix $H$, let $e_j(H)$ denote the $j$th elementary
symmetric polynomial of its eigenvalues, which by multilinearity of
$\det(I+zH)$ equals the sum of all $j\times j$ principal minors:
\[
 e_j(H)=\sum_{\substack{I\subset\{1,\ldots,n\}\\|I|=j}}\det H_{I,I}
 \quad(1\le j\le n),\qquad e_{n+1}(H)=0.
\]
These quantities are orthogonally invariant, with $e_1(H)=\tr H$ and
$e_n(H)=\det H$. Classical smooth BFGS theory couples $e_1$ directly to $e_n$,
which fails when nonsmooth jumps make $\norm{y_k}^2/t_k$ unbounded. Tracking all
intermediate degrees $e_1(H),\ldots,e_n(H)$ of the \emph{inverse} matrix $H_k$
resolves this obstacle: each inverse-BFGS update shifts weight upward by at most
one minor degree with factor $\norm{y}^2/(\beta^2q)$, while retaining degree $j$
only with the step-to-curvature ratio $\alpha/\beta$.

\begin{theorem}[Principal-minor hierarchy]\label{art:minor-hierarchy}
Let $H\succ0$, $g,y\in\R^n$, $\alpha>0$, and $s=-\alpha Hg$, and suppose
\[
 q=g^THg>0,\qquad \beta=-\frac{y^THg}{q}>0.
\]
Then the inverse-BFGS update $H_+$ in~\eqref{art:update} with $t=y^Ts=\alpha\beta q$
satisfies
\begin{equation}\label{art:hierarchy}
 e_j(H_+)\le\frac{\norm{y}^2}{\beta^2q}e_{j+1}(H)
                   +\frac\alpha\beta e_j(H),\qquad 1\le j\le n.
\end{equation}
At $j=n$, the exact determinant identity $\det H_+=(\alpha/\beta)\det H$ holds.
\end{theorem}
\begin{proof}
Let $Z\in\R^{n\times j}$ have orthonormal columns. The subspace Gram compression
inequality (Lemma~\ref{art:gram-compression} in Appendix~\ref{art:matrix-algebra})
gives
\[
 \det(Z^TH_+Z)\le
 \frac{\det\operatorname{Gram}_H(y,Z)}{\beta^2q}
       +\frac\alpha\beta\det(Z^THZ),
\]
where $\operatorname{Gram}_H(y,Z)=[y\ Z]^TH[y\ Z]$. Working in an orthonormal
eigenbasis of $H$ with eigenvalues $\lambda_1,\ldots,\lambda_n>0$ and writing
$y_i$ for the coordinates of $y$ in that basis, we sum over all $\binom{n}{j}$
coordinate column selections $Z_I$:
\begin{align*}
 \sum_{|I|=j}\det\operatorname{Gram}_H(y,Z_I)
 &=\sum_{|I|=j}\left(\prod_{i\in I}\lambda_i\right)
                  \sum_{\ell\notin I}\lambda_\ell y_\ell^2\\
 &=\sum_{|J|=j+1}\left(\prod_{i\in J}\lambda_i\right)
                  \sum_{\ell\in J}y_\ell^2
 \le\norm{y}^2e_{j+1}(H),
\end{align*}
where the first equality evaluates the diagonal Schur complement and the second
regroups pairs $(I,\ell)$ by $J=I\cup\{\ell\}$. Orthogonal invariance of $e_j$
yields~\eqref{art:hierarchy}, and the $j=n$ case of
Lemma~\ref{art:gram-compression} gives $\det H_+=(\alpha/\beta)\det H$.\qed
\end{proof}

Theorem~\ref{art:minor-hierarchy} is a purely algebraic property of
the rank-two BFGS update: it makes no reference to convexity or smoothness of an
objective function. Whenever gradients remain bounded along a run and the
predicted decreases $d_k$ tend to zero, cascading~\eqref{art:hierarchy} across
$j=1,\ldots,n$ forces both the geometric and arithmetic means of
$q_k^{-1}=(g_k^TH_kg_k)^{-1}$ to diverge to $+\infty$.

\begin{theorem}[Universal reciprocal growth]\label{art:reciprocal-growth}
If an admissible run has bounded gradients and $d_k\to0$, then
\begin{equation}\label{art:reciprocal-limits}
 \left(\prod_{k<N}q_k^{-1}\right)^{1/N}\longrightarrow+\infty,
 \qquad \frac1N\sum_{k<N}q_k^{-1}\longrightarrow+\infty.
\end{equation}
\end{theorem}
\begin{proof}
Choose $\Lambda>0$ with $\norm{g_k}\le\Lambda$ for all $k$, and set
$K_0=4\Lambda^2/\omega^2$. Along the run,
\[
 \beta_k=-\frac{y_k^TH_kg_k}{q_k}=\frac{t_k}{d_k}\ge\omega,\qquad
 \norm{y_k}\le2\Lambda,\qquad \alpha_k=\frac{d_k}{q_k}.
\]
Multiplying~\eqref{art:hierarchy} by $q_k$ yields the normalized hierarchy
\begin{equation}\label{art:normalized-hierarchy}
 q_ke_j(H_{k+1})\le
 K_0e_{j+1}(H_k)+\frac{d_k}{\omega}e_j(H_k),\qquad 1\le j\le n.
\end{equation}

Given any $\eta>0$, choose $J$ so that $d_k/\omega\le\eta$ for all $k\ge J$, and
define the weighted potential
\[
 V_\eta(H)=\sum_{j=1}^n(K_0/\eta)^{j-1}e_j(H)>0\qquad(H\succ0).
\]
Summing~\eqref{art:normalized-hierarchy} with weights $(K_0/\eta)^{j-1}$ for
$k\ge J$, the $e_{j+1}(H_k)$ terms telescope into at most $\eta V_\eta(H_k)$
(since $e_{n+1}(H_k)=0$), and the $e_j(H_k)$ terms contribute at most
$\eta V_\eta(H_k)$ as well. Since $q_k\le\Lambda^2\tr H_k=\Lambda^2e_1(H_k)\le\Lambda^2V_\eta(H_k)$,
we obtain $V_\eta(H_{k+1})\le(2\eta/q_k)V_\eta(H_k)$ for all $k\ge J$. Iterating
over $m-1$ steps from $J$ and bounding $q_{J+m-1}\le\Lambda^2V_\eta(H_{J+m-1})$
at the final index gives
\[
 \prod_{k=J}^{J+m-1}q_k^{-1}
 \ge\frac{1}{\Lambda^2V_\eta(H_J)(2\eta)^{m-1}}.
\]
Taking $m$th roots and accounting for the prefix $0\le k<J$ gives
$\liminf_{N\to\infty}(\prod_{k<N}q_k^{-1})^{1/N}\ge(2\eta)^{-1}$. Since $\eta>0$
is arbitrary, the geometric mean diverges to $+\infty$, and the arithmetic--geometric
mean inequality gives the second limit in~\eqref{art:reciprocal-limits}.\qed
\end{proof}

Theorem~\ref{art:reciprocal-growth} extends product-type BFGS estimates
developed for smooth and restricted strongly convex
settings~\cite[Theorem~3.7]{S16}, \cite[Proposition~D.1]{FND_JJM24AW} to
arbitrary nonsmooth runs with bounded gradients and $d_k\to0$ (weaker even than
the summability $\sum_kd_k<\infty$ implied by Armijo on bounded-below
objectives; non-asymptotic finite-horizon bounds are recorded in Chapter~2 of
the companion).

\section{Static Radial Energies: Homogeneous and Polyhedral Convergence in \texorpdfstring{$\R^n$}{R\^{}n}}\label{art:global}

The nullspace identity $\widehat B_ks_k=0$ of~\eqref{art:hessian-update} erases
the prior Hessian approximation along the step $s_k$ before the new secant term
$y_ky_k^T/t_k$ is installed. Evaluating $B_k$ along the static radial vector
$r_k=x_k-x_*$ therefore makes the step shift $r_{k+1}=r_k+s_k$ algebraically
free and immediately establishes global and local convergence across every
class in $\R^n$ whose kinks align along rays or local polyhedral cones around
$x_*$.

\subsection{The radial reference energy and the Guo--Lewis benchmark \texorpdfstring{$F(u,v)=u^2+|v|$}{F(u,v) = u\^{}2 + |v|}}

\begin{lemma}[Reference energy identity]\label{art:reference-energy}
For any fixed point $x_*\in\R^n$, set $r_k=x_k-x_*$ and
\[
 E_k=r_k^TB_kr_k,\qquad a_k=g_k^Tr_k,\qquad
 R_k=\frac{(y_k^Tr_{k+1})^2}{t_k}.
\]
Then every step of an admissible run satisfies the exact identity
\begin{equation}\label{art:energy-identity}
 E_{k+1}=E_k-\frac{a_k^2}{q_k}+R_k.
\end{equation}
In particular, if $x_*$ minimizes $F$ with minimum value $F_*=F(x_*)$, then for
every $N\ge1$,
\begin{equation}\label{art:energy-budget}
 (F(x_N)-F_*)^2\sum_{k<N}q_k^{-1}
       \le E_0+\sum_{k<N}R_k.
\end{equation}
\end{lemma}
\begin{proof}
Because $r_{k+1}=r_k+s_k$ and $\widehat B_ks_k=0$, evaluating~\eqref{art:hessian-update}
at $r_{k+1}$ gives
\[
 r_{k+1}^TB_{k+1}r_{k+1}
 =r_{k+1}^T\widehat B_kr_{k+1}+\frac{(y_k^Tr_{k+1})^2}{t_k}
 =r_k^T\widehat B_kr_k+R_k
 =E_k-\frac{a_k^2}{q_k}+R_k,
\]
which is~\eqref{art:energy-identity}. When $x_*$ minimizes $F$, convexity and
monotonicity of $F(x_k)$ give $a_k\ge F(x_k)-F_*\ge F(x_N)-F_*\ge0$ for all
$k<N$; summing~\eqref{art:energy-identity} over $0\le k<N$ and dropping
$E_N\ge0$ yields~\eqref{art:energy-budget}.\qed
\end{proof}

Inequality~\eqref{art:energy-budget} isolates the radial requirement of this
section: if $F(x_k)-F_*\ge\eta>0$ for all $k$ while the radial secant error
satisfies $R_k\le C$, then $\eta^2\sum_{k<N}q_k^{-1}\le E_0+CN$, bounding the
average of $q_k^{-1}$ uniformly in $N$ while imposing no upper bound on the
other $n-1$ directions of $B_k$. Pairing~\eqref{art:energy-budget} with
Theorem~\ref{art:reciprocal-growth} immediately settles the benchmark
$F(u,v)=u^2+|v|$ on which Guo and Lewis~\cite[Corollary~4.4]{GL18} proved
$\liminf_{k\to\infty}|v_k|=0$ while leaving whole-sequence convergence open.

\begin{corollary}[Whole-sequence convergence on $u^2+|v|$]\label{art:mixed-benchmark}
For $F(u,v)=u^2+|v|$, every infinite admissible full-BFGS run, with arbitrary
$H_0\succ0$ and any fixed $0<c_1<c_2<1$, satisfies
\[
 F(u_k,v_k)\longrightarrow0,\qquad (u_k,v_k)\longrightarrow(0,0).
\]
\end{corollary}
\begin{proof}
Set $b=F(u_0,v_0)>0$; Armijo keeps $x_k=(u_k,v_k)$ inside the compact sublevel
$K=\{F\le b\}$, where gradients $g_k=(2u_k,\operatorname{sign}v_k)$ satisfy
$\norm{g_k}^2=4u_k^2+1\le4b+1$ and $\sum_kd_k<\infty$. Consider the reference
energy $E_k=x_k^TB_kx_k$ centered at $x_*=(0,0)$. For any two differentiable
points $x=(u,v)$ and $z=(r,w)$ with $y=\nabla F(z)-\nabla F(x)$ and
$t=y^T(z-x)>0$:
\begin{itemize}
\item if $v$ and $w$ have the same sign, then $y=(2(r-u),0)$, so
$t=2(r-u)^2$ and $(y^Tz)^2=4(r-u)^2r^2=2r^2t\le2F(z)t$;
\item if $v$ and $w$ have opposite signs, then with $S=|v|+|w|>0$ we have
$t=2(r-u)^2+2S$ and $y^Tz=2(r-u)r+2|w|$, so Cauchy--Schwarz gives
\[
 (y^Tz)^2\le\bigl(2(r-u)^2+2S\bigr)\left(2r^2+\frac{2|w|^2}{S}\right)
 \le2(r^2+|w|)t=2F(z)t.
\]
\end{itemize}
In either case, $R_k=(y_k^Tx_{k+1})^2/t_k\le2F(x_{k+1})\le2b$ for every step
$k$. By~\eqref{art:energy-budget} of Lemma~\ref{art:reference-energy}, a
positive limiting value $\lim_kF(x_k)=\eta>0$ would imply
$\eta^2\sum_{k<N}q_k^{-1}\le E_0+2bN$, contradicting
Theorem~\ref{art:reciprocal-growth}. Hence $u_k^2+|v_k|\to0$ and
$(u_k,v_k)\to(0,0)$.\qed
\end{proof}

The benchmark proof above compares gradients at two points $(x,z)$: the
cross-kink jump $y_k=(2(u_{k+1}-u_k),\pm2)$ satisfies $(y_k^Tx_{k+1})^2\le2bt_k$
even when $s_k\to0$ drives $\norm{y_k}^2/t_k\to+\infty$, because along any
single inward ray toward $(0,0)$ the term $|v|$ scales linearly while $u^2$
contributes only the quadratic Taylor error~\eqref{art:mixed-inward-remainder}
(an explicit admissible run with $H_k=\operatorname{diag}(1,2^{-k})\to\operatorname{diag}(1,0)$
is given in Proposition~\ref{art:examples-mixed}). We next give a condition on
the objective along one inward ray that implies the required two-point radial
secant estimate for arbitrary finite convex functions on $\R^n$.

\subsection{From one-point inward models to two-point secant control}

Fix a minimizer $x_*$ of finite convex $F:\R^n\to\R$, write $F_*=F(x_*)$, and
let $K=\{x:F(x)\le b\}$ be compact with $b>F_*$. For differentiable points
$x,z\in K$, write
\[
 \begin{gathered}
 r_x=x-x_*,\qquad g=\nabla F(x),\qquad h=\nabla F(z),\\
 T=(h-g)^T(z-x),\qquad N=(h-g)^Tr_z,
 \end{gathered}
\]
so $T$ is the secant curvature $t$ of Section~\ref{art:energy} and $N$ is the
radial projection of the gradient jump at $z$. The \emph{radial secant
condition} requires control of $N^2/T$ across all descending differentiable
pairs on each positive sublevel band:
\begin{equation}\label{art:radial-condition}
 \begin{gathered}
 \text{for every }\eta>0\text{ there exists }C_\eta\ge0\text{ such that}\\
 N^2\le C_\eta T
 \quad\text{whenever }F_*+\eta\le F(z)\le F(x)\text{ and }T>0.
 \end{gathered}
\end{equation}
Rather than checking pairs of points $(x,z)$ across different smooth regions by
hand, the \emph{inward model condition} tests $F$ at a single differentiable
point $x$ along the inward ray toward $x_*$:
\begin{equation}\label{art:inward-condition}
 \begin{gathered}
 \text{for every }\eta>0\text{ there exist }M_\eta\ge0\text{ and }
       \theta_\eta\in(0,1/2]\text{ such that}\\
 0\le F(x-\tau r_x)-F(x)+\tau\nabla F(x)^Tr_x
       \le\frac{M_\eta}{2}\tau^2\quad(0\le\tau\le\theta_\eta)
 \end{gathered}
\end{equation}
at every differentiable $x\in K$ with $F(x)\ge F_*+\eta$. By convexity, the
inward sample $x-\tau r_x$ remains in $K$ and need not be differentiable; $\tau$
is a dimensionless fraction of the ray to $x_*$, and both $M_\eta$ and
$\theta_\eta$ may depend on the band $\eta>0$.

\begin{theorem}[Inward models control radial secants]\label{art:radial-equivalence}
For any finite convex $F:\R^n\to\R$ and compact sublevel $K=\{F\le b\}$, the
one-point inward condition~\eqref{art:inward-condition} implies the two-point
radial secant condition~\eqref{art:radial-condition}.
\end{theorem}
\begin{proof}
Because $F$ is finite and convex, it is continuous and bounded on a closed
$\epsilon$-neighborhood $K+\epsilon\overline B$; evaluating
$F(x+\epsilon u)\ge F(x)+g^T(\epsilon u)$ over $\norm{u}=1$ supplies $\Lambda>0$
bounding all subgradient norms on $K$, and compactness supplies $R>0$ bounding
$\norm{x-x_*}$ on $K$. Define the Bregman supporting-plane
residual~\cite{BB}
\[
 D_F(z,x)=F(z)-F(x)-\nabla F(x)^T(z-x)\ge0,
\]
which requires differentiability only at its second argument and satisfies
$T=D_F(z,x)+D_F(x,z)$ when both $x$ and $z$ are differentiable.

Fix $\eta>0$ and let $M=M_\eta$, $\theta=\theta_\eta$ be the constants
from~\eqref{art:inward-condition}. Let $x,z\in K$ be differentiable with
$F_*+\eta\le F(z)\le F(x)$ and $T>0$. If $N\ge0$, comparing the supporting
inequality $F(z-\tau r_z)\ge F(x)+\nabla F(x)^T(z-\tau r_z-x)$ at $x$ with the
upper inward bound~\eqref{art:inward-condition} at $z$ yields
\[
 \tau N\le D_F(z,x)+\tfrac12M\tau^2\qquad(0\le\tau\le\theta).
\]
If $N<0$, reversing the roles of $x$ and $z$ (applying the supporting plane at
$z$ and the inward bound at $x$, where $F(x)\ge F(z)\ge F_*+\eta$) gives
$(g-h)^Tr_x=T-N\ge|N|$ and
\[
 \tau(T-N)\le D_F(x,z)+\tfrac12M\tau^2\qquad(0\le\tau\le\theta).
\]
In both cases there exist $A\in[|N|,2\Lambda R]$ and $D\in\{D_F(z,x),D_F(x,z)\}$
satisfying $\tau A\le D+\tfrac12M\tau^2$ for all $\tau\in[0,\theta]$. Optimizing
over $\tau$: if $A=0$, the bound is immediate; if $A>0$, $M>0$, and
$A\le M\theta$, setting $\tau=A/M$ gives $A^2\le2MD$; if $A>M\theta$ (including
$A>0$ with $M=0$), setting $\tau=\theta$ gives
$D\ge\tau A-\tfrac12M\tau^2\ge\tfrac12\theta A$ and hence
$A^2\le(4\Lambda R/\theta)D$. Since $D\le T$, we conclude that
\begin{equation}\label{art:inward-to-radial}
 N^2\le \max\{2M,4\Lambda R/\theta\}\,T.\qed
\end{equation}
\end{proof}

Chapter~4 of the companion proves that the converse implication also holds on
compact sublevels, so~\eqref{art:radial-condition} and~\eqref{art:inward-condition}
are in fact equivalent characterizations of radial comparison geometry.

\begin{corollary}[Radial convergence]\label{art:radial-convergence}
Under either~\eqref{art:radial-condition} or~\eqref{art:inward-condition} on a
compact sublevel $K=\{F\le b\}$, every admissible run starting at $x_0\in K$
satisfies $F(x_k)\to F_*$ and $\dist(x_k,\argmin F)\to0$; if $x_*$ is the unique
minimizer of $F$, then $x_k\to x_*$.
\end{corollary}
\begin{proof}
By Theorem~\ref{art:radial-equivalence},~\eqref{art:radial-condition} holds on
$K$. Armijo keeps $x_k\in K$ with $F(x_k)\downarrow\ell\ge F_*$, so
Lemma~\ref{art:run-invariants} gives $\sum_kd_k<\infty$, and subgradient
boundedness on $K$ allows us to apply Theorem~\ref{art:reciprocal-growth}. If
$\ell>F_*$, set $\eta=(\ell-F_*)/2>0$; then $F_*+\eta<\ell\le F(x_{k+1})<F(x_k)$
and $T=t_k>0$ for every $k$, so~\eqref{art:radial-condition} gives
$R_k=(y_k^Tr_{k+1})^2/t_k\le C_\eta$, and~\eqref{art:energy-budget} of
Lemma~\ref{art:reference-energy} yields
$\eta^2\sum_{k<N}q_k^{-1}\le E_0+C_\eta N$, contradicting
Theorem~\ref{art:reciprocal-growth}. Thus $F(x_k)\to F_*$, and compactness of
$K$ together with continuity of $F$ gives $\dist(x_k,\argmin F)\to0$ (hence
$x_k\to x_*$ when unique).\qed
\end{proof}

Geometrically,~\eqref{art:inward-condition} demands quadratic first-order
accuracy \emph{only along rays toward $x_*$}: in smooth regions its
infinitesimal form is $r_x^T\nabla^2F(x)r_x$ rather than $\lambda_{\max}(\nabla^2F(x))$,
while its finite-difference form also rules out uncompensated radial derivative
jumps (as on $|x-1|$ relative to $x_*=0$). Likewise, any degree-one positively
homogeneous summand centered at $x_*$ has zero inward error regardless of how
its gradients jump in angular directions.

\subsection{Centered smooth-plus-homogeneous objectives in \texorpdfstring{$\R^n$}{R\^{}n}}

Any finite convex positively homogeneous function $\gamma:\R^n\to\R$ of degree
one satisfies $\gamma((1-\tau)r)=(1-\tau)\gamma(r)$ for $0\le\tau\le1$, so its
inward Taylor remainder at $x_*$ vanishes identically regardless of angular
nonsmoothness. Adding an arbitrary locally $C^{1,1}$ convex function $\psi$
contributes only a smooth quadratic remainder, even when neither $\psi$ nor
$\gamma(\cdot-x_*)$ attains a minimum at $x_*$ on its own.

\begin{theorem}[Centered smooth-plus-homogeneous objectives]\label{art:centered}
Let $F:\R^n\to\R$ be finite and convex, and let $x_*$ minimize $F$. Let
$U\subseteq\R^n$ be an open convex set, let $\psi:U\to\R$ be convex and locally
$C^{1,1}$, and let $\gamma:\R^n\to\R$ be finite, convex, and positively
homogeneous of degree one. Suppose
\[
 F(x)=\psi(x)+\gamma(x-x_*)\qquad(x\in U).
\]
Then every admissible run whose initial sublevel $K=\{x:F(x)\le F(x_0)\}$ is
compact and contained in $U$ satisfies $F(x_k)\to F(x_*)$ and
$\dist(x_k,\argmin F)\to0$; if $x_*$ is the unique minimizer, then $x_k\to x_*$.
\end{theorem}
\begin{proof}
Since $g_0\ne0$, convexity gives $F(x_0)>F(x_*)$. Because $K\subset U$ is
compact and convex, $\nabla\psi$ has a uniform Lipschitz constant $L$ on $K$:
covering $K$ by finitely many open sets on which $\nabla\psi$ is Lipschitz and
subdividing any segment in $K$ into pieces shorter than a Lebesgue number of
that cover bounds the variation of $\nabla\psi$ along the segment by the
maximum local Lipschitz constant. Integrating along the inward segment from $x$
to $x-\tau r_x$ therefore gives
\[
 \left|\psi(x-\tau r_x)-\psi(x)+\tau\nabla\psi(x)^Tr_x\right|
 \le\tfrac12L\tau^2\norm{r_x}^2\qquad(0\le\tau\le1).
\]
At any differentiable point $x\in K$ of $F$, smoothness of $\psi$ implies that
$\gamma$ is differentiable at $r_x=x-x_*$, and Euler's identity for degree-one
homogeneous functions gives $\nabla\gamma(r_x)^Tr_x=\gamma(r_x)$, so
$\gamma((1-\tau)r_x)-\gamma(r_x)+\tau\nabla\gamma(r_x)^Tr_x=0$. Hence
\eqref{art:inward-condition} holds on $K$ with uniform constants
$M_\eta=L\max_{x\in K}\norm{r_x}^2$ and $\theta_\eta=1/2$ (the lower bound $0$
follows from convexity of $F$), and Corollary~\ref{art:radial-convergence}
completes the proof.\qed
\end{proof}

Theorem~\ref{art:centered} recovers
Corollary~\ref{art:mixed-benchmark} immediately in $\R^2$ by taking
$\psi(u,v)=u^2$, $\gamma(u,v)=|v|$, $x_*=(0,0)$, and $U=\R^2$, and extends it
to arbitrary smooth-plus-norm or smooth-plus-support-function sums in $\R^n$
such as $\psi(x)+\norm{Ax}_1$ or $\psi(x)+\max_ia_i^Tx$ whenever the composite
minimizer $x_*$ centers the homogeneous term (or locally after active-sign
splitting, as we exploit for the Lasso in Section~\ref{art:polyhedral-subsec}
below).

\subsection{Coercive positively homogeneous functions, norms, and linear maxima}

When $\psi=0$, the inward remainder vanishes everywhere on $\R^n$. Imposing
coercivity prevents rays of zero or negative slope and turns value convergence
into global iterate convergence from every starting point, with no smoothness or
finite-piece hypothesis.

\begin{theorem}[Coercive positively homogeneous objectives]\label{art:homogeneous}
Let $F:\R^n\to\R$ ($n\ge1$) be finite, convex, positively homogeneous of degree
one, and coercive. Then $0$ is the unique minimizer of $F$, and every infinite
admissible full-BFGS run from any differentiable $x_0\ne0$, with any
$H_0\succ0$ and fixed $0<c_1<c_2<1$, satisfies
\[
 F(x_k)\longrightarrow0,\qquad x_k\longrightarrow0.
\]
\end{theorem}
\begin{proof}
Positive homogeneity gives $F(0)=0$. If $F(x)\le0$ for some $x\ne0$, then
$F(tx)=tF(x)\le0$ for all $t>0$, contradicting coercivity as $t\to\infty$. Thus
$F(x)>0$ for all $x\ne0$, and continuity on the compact Euclidean unit sphere
gives
\[
 \mu=\min_{\norm{x}=1}F(x)>0,\qquad F(x)\ge\mu\norm{x}\quad(x\in\R^n).
\]
Every sublevel of $F$ is therefore compact and $0$ is the unique minimizer.
Applying Theorem~\ref{art:centered} with $x_*=0$, $U=\R^n$, $\psi=0$, and
$\gamma=F$ gives $F(x_k)\to0$, and $F(x_k)\ge\mu\norm{x_k}$ gives $x_k\to0$.\qed
\end{proof}

\begin{corollary}[Norms, gauges, and coercive linear maxima]\label{art:homogeneous-examples}
The global convergence conclusion of Theorem~\ref{art:homogeneous} holds on
$\R^n$ for:
\begin{enumerate}
\item every norm on $\R^n$ (in particular $\norm{x}_1$ and $\norm{x}_\infty$);
\item every finite convex degree-one gauge whose unit sublevel $\{x:F(x)\le1\}$
is bounded; and
\item every finite linear maximum $F(x)=\max_{1\le i\le m}a_i^Tx$ whose zero
sublevel $\{x:F(x)\le0\}$ is bounded.
\end{enumerate}
\end{corollary}
\begin{proof}
Every norm is finite, convex, degree-one homogeneous, and positive on the
Euclidean unit sphere, hence satisfies $F(x)\ge\mu\norm{x}$ with $\mu>0$. For a
degree-one gauge with $\{F\le1\}$ bounded or a finite linear maximum with
$\{F\le0\}$ bounded, no nonzero $x$ can satisfy $F(x)\le0$, since the entire ray
$\{tx:t\ge0\}$ would then lie inside that bounded sublevel; compactness of the
unit sphere again yields $F(x)\ge\mu\norm{x}$ with $\mu>0$, so $F$ is coercive
and Theorem~\ref{art:homogeneous} applies.\qed
\end{proof}

Corollary~\ref{art:homogeneous-examples} goes well beyond isolated kinks: in
dimension $n\ge2$, $\norm{x}_1=\sum_i|x_i|$ is nondifferentiable along every
coordinate hyperplane, and $\norm{x}_\infty=\max_i|x_i|$ has gradient jumps
along every diagonal hyperplane $|x_i|=|x_j|$ where two largest coordinates tie,
recovering and generalizing the Euclidean-norm theorem of Guo and
Lewis~\cite[Corollary~4.2]{GL18} and the one-dimensional $|x|$ analysis of Lewis
and Overton~\cite[Section~5.1]{LO} without any symmetry requirement. When
coercivity is dropped on a positively homogeneous function, Chapter~7 of the
companion proves the global value dichotomy $F(x_k)\to0$ or $F(x_k)\to-\infty$.
By contrast, Section~\ref{art:failure} (Theorem~\ref{art:failure_polyhedral})
shows that absence of nonzero affine intercepts in $\max_ia_i^Tx$ is sharp:
adding a single nonzero offset in $\R^2$ produces a coercive four-piece
polyhedral objective on which full BFGS converges to a nonstationary corner.

\subsection{Local convergence: sublevel localization, polyhedral partitions, and composite penalties}\label{art:polyhedral-subsec}

Many composite objectives in nonsmooth optimization, including smooth losses
regularized by $\ell_1$, hinge, or group-norm penalties, admit structured radial
or polyhedral representations only near their minimizer $x_*$. Because an
arbitrary initial matrix $H_0\succ0$ can send the first BFGS step in almost any
descent direction, the iterates need not remain in a prescribed Euclidean ball
around $x_*$. Convexity and uniqueness of $x_*$ resolve this immediately:
sufficiently low sublevels $\{F\le F_*+\varepsilon\}$ lie inside any prescribed
neighborhood of $x_*$ and trap every Armijo run independently of $H_0$.

\begin{lemma}[Matrix-independent sublevel basin]\label{art:localization}
Let $F:\R^n\to\R$ be finite and convex with a unique minimizer $x_*$. Then $F$
is coercive, every sublevel of $F$ is compact, and for every neighborhood $U$
of $x_*$ there exists $\varepsilon>0$ such that
$\{x:F(x)\le F(x_*)+\varepsilon\}\subset U$. Consequently, every admissible run
starting in that sublevel remains in $U$ for all $k\ge0$, regardless of
$H_0\succ0$.
\end{lemma}
\begin{proof}
Choose $r>0$ with $\overline B_r(x_*)\subset U$. Continuity of $F$, compactness
of the sphere $\{\norm{u}=r\}$, and uniqueness of $x_*$ give
\[
 \delta=\min_{\norm{u}=r}\bigl(F(x_*+u)-F(x_*)\bigr)>0.
\]
For any $x$ with $\norm{x-x_*}\ge r$, evaluating convexity at the point at
distance $r$ along the segment from $x_*$ to $x$ yields
\[
 F(x)-F(x_*)\ge\frac{\norm{x-x_*}}r\delta.
\]
Thus $F$ is coercive with compact sublevels, $\varepsilon=\delta/2$ gives
$\{F\le F(x_*)+\varepsilon\}\subset B_r(x_*)\subset U$, and Armijo descent
preserves this sublevel.\qed
\end{proof}

Near $x_*$, a finite polyhedral partition simplifies because every closed
cell that does not contain $x_*$ stays bounded away from $x_*$, so within a
small ball around $x_*$ every active cell is a convex polyhedron containing
$x_*$ itself. Any inward ray $x-\tau(x-x_*)$ therefore remains inside a single
cell, giving quadratic inward accuracy even when gradients jump across cell
facets and without requiring BFGS to identify active manifolds.

\begin{theorem}[Local convergence on finite polyhedral partitions]\label{art:polyhedral}
Let $n\ge1$ and let $F:\R^n\to\R$ be finite and convex with a unique minimizer
$x_*$. Suppose an open neighborhood $U$ of $x_*$ is covered by finitely many
closed full-dimensional polyhedra $P_i$, and for each $i$ there exist an open
set $V_i\supset U\cap P_i$ and a function $f_i\in C^{1,1}_{\rm loc}(V_i)$
(not necessarily convex outside $U\cap P_i$) such that
\[
 F(x)=f_i(x)\qquad(x\in U\cap P_i).
\]
Then there exists $\varepsilon>0$, depending only on $F$, $U$, and the partition
data (and independent of $H_0$, $c_1$, and $c_2$), such that
every infinite admissible full-BFGS run with arbitrary $H_0\succ0$, fixed
$0<c_1<c_2<1$, and initial value $F(x_0)\le F(x_*)+\varepsilon$ satisfies
\[
 x_k\in U\quad(k\ge0),\qquad F(x_k)\longrightarrow F(x_*),\qquad x_k\longrightarrow x_*.
\]
\end{theorem}

\begin{proof}[Proof of Theorem~\ref{art:polyhedral}]
Each closed polyhedron $P_i$ not containing $x_*$ has an open complement
containing $x_*$; since there are finitely many cells, one ball around $x_*$
avoids all cells with $x_*\notin P_i$. Shrink to a ball
$B_\rho(x_*)\subset U\cap\bigcap_{i:x_*\in P_i}V_i$ on which the gradients
$\nabla f_i$ of all remaining extensions share a common Lipschitz constant $L$.

At any differentiability point $x\in B_\rho(x_*)$ of $F$, pick a cell $P_i$
containing $x$. Since $x_*\in P_i$ and $P_i$ is convex, the entire inward
segment $x-\tau(x-x_*)$ ($0\le\tau\le1$) lies in $U\cap P_i$, where $F=f_i$.
Differentiating along the ray at $\tau=0^+$ gives
\[
 \nabla F(x)^T(x-x_*)=\nabla f_i(x)^T(x-x_*),
\]
and the smooth Taylor bound for $f_i$ along the segment yields
\[
 0\le F(x-\tau r_x)-F(x)+\tau\nabla F(x)^Tr_x
       \le\tfrac12L\tau^2\norm{r_x}^2
       \le\tfrac12L\rho^2\tau^2,
\]
where the lower bound uses convexity of $F$. Lemma~\ref{art:localization}
supplies $\varepsilon>0$ with $K=\{F\le F(x_*)+\varepsilon\}\subset B_\rho(x_*)$,
so~\eqref{art:inward-condition} holds on $K$ with $M_\eta=L\rho^2$ and
$\theta_\eta=1/2$. Corollary~\ref{art:radial-convergence} and Armijo descent
complete the proof.\qed
\end{proof}

\begin{corollary}[Piecewise linear-quadratic and composite objectives]\label{art:plq}
Let $F:\R^n\to\R$ have a unique minimizer $x_*$, and suppose $F$ is either a
finite convex piecewise linear-quadratic function or the sum of a convex locally
$C^{1,1}$ loss and a finite convex piecewise linear-quadratic regularizer. Then
$F$ enjoys the local sublevel convergence guarantee of
Theorem~\ref{art:polyhedral}.
\end{corollary}
\begin{proof}
Here piecewise linear-quadratic means represented by quadratic polynomials on a
finite cover by closed full-dimensional polyhedra (possibly with distinct
Hessians on different cells). Each quadratic polynomial extends smoothly to
$\R^n$, and adding a locally $C^{1,1}$ loss preserves local $C^{1,1}$ cell
extensions, so Theorem~\ref{art:polyhedral} applies.\qed
\end{proof}

Localizing Theorem~\ref{art:centered} via Lemma~\ref{art:localization} yields a
parallel local guarantee that also encompasses \emph{nonpolyhedral} group-norm
regularizers whose kinks have curved cross-sections away from the origin.

\begin{corollary}[Local smooth-plus-homogeneous basin]\label{art:local-centered}
If the local representation $F(x)=\psi(x)+\gamma(x-x_*)$ of
Theorem~\ref{art:centered} holds on an open neighborhood $U$ of a unique
minimizer $x_*$, then there exists $\varepsilon>0$, depending only on $F$ and
$U$ (and independent of $H_0$, $c_1$, and $c_2$), such that every admissible run
with $F(x_0)\le F(x_*)+\varepsilon$ stays in $U$ and satisfies
$F(x_k)\to F(x_*)$ and $x_k\to x_*$.
\end{corollary}
\begin{proof}
Shrink $U$ to an open convex neighborhood $U_0\subseteq U$ of $x_*$ on which the
representation holds; Lemma~\ref{art:localization} supplies a compact sublevel
$\{F\le F(x_*)+\varepsilon\}\subset U_0$, and Theorem~\ref{art:centered}
applies.\qed
\end{proof}

To see both corollaries at work on composite regularization, consider first the
$\ell_1$-regularized least-squares (Lasso) objective
$F(x)=\tfrac12\norm{Ax-b}^2+\lambda\norm{x}_1$ ($A\in\R^{m\times n}$, $\lambda>0$)
with a unique minimizer $x_*$ (already covered by Corollary~\ref{art:plq}).
Setting \mbox{$I_0=\{i:x_{*,i}=0\}$} and $\sigma_i=\operatorname{sign}(x_{*,i})$
for $i\notin I_0$, on the open convex neighborhood
$\{x:\sigma_ix_i>0\text{ for all }i\notin I_0\}$ we have the exact centered
splitting
\begin{equation}\label{art:local-l1-decomposition}
 \begin{aligned}
 F(x)&=\psi(x)+\gamma(x-x_*),\\
 \psi(x)&=\tfrac12\norm{Ax-b}^2+\lambda\sum_{i\notin I_0}\sigma_ix_i,\qquad
 \gamma(r)=\lambda\sum_{i\in I_0}|r_i|,
 \end{aligned}
\end{equation}
placing it directly under Corollary~\ref{art:local-centered} as well. Now
replace $\lambda\norm{x}_1$ by a sum of Euclidean group-norm residuals
$\sum_j\norm{A_jx-b_j}_2$ added to any convex locally $C^{1,1}$ loss: near a
unique minimizer $x_*$, every group with $A_jx_*-b_j\ne0$ is smooth, while every
active group with $A_jx_*-b_j=0$ equals the centered homogeneous seminorm
$\norm{A_j(x-x_*)}_2$. Corollary~\ref{art:local-centered} therefore guarantees
local convergence even though $\norm{A_j(x-x_*)}_2$ need not be polyhedral.

\section{Moving Comparison Fields and Curved Switching Interfaces}\label{art:local}\label{art:positive}

Section~\ref{art:global} showed that the static radial energy
$E_k=(x_k-x_*)^TB_k(x_k-x_*)$ suffices whenever nonsmooth kinks align along rays
from $x_*$ or lie on the facets of a local polyhedral partition containing
$x_*$. When an objective is instead a pointwise maximum
$F=\max_{1\le i\le m}f_i$ of smooth strongly convex functions whose active
Hessians differ at $x_*$, the switching interfaces curve parabolically away from
their tangent rays at $x_*$. Inward rays toward $x_*$ then cross the switching
locus transversely, and the radial secant condition~\eqref{art:radial-condition}
can fail on every neighborhood of $x_*$. For instance, on the explicit
three-quadratic maximum of Theorem~\ref{art:nonradial-example} below, the tie
curve $f_0(u,v)=f_1(u,v)$ satisfies $v=u-\tfrac12u^2+O(u^3)$, along which the
gradient jump pairs nontrivially with the radial vector,
$(\nabla f_0-\nabla f_1)^T(u,v)=\tfrac12u^2-v^2\ne0$, driving
$(y^Tz)^2/t\to+\infty$ for short eligible steps straddling the curve.

Overcoming this radial obstruction requires replacing the fixed vector $x_k-x_*$
by a \emph{moving} comparison field that rotates with the switching interface.
We first formulate a general $\R^n$ convergence criterion for transported matrix
and vector comparison fields, and then construct normal-angle contour-velocity
fields in $\R^2$ that resolve curved smooth maxima, including the
three-quadratic benchmark on which radial energies fail.

\subsection{Transported band certificates and tail vector comparison in \texorpdfstring{$\R^n$}{R\^{}n}}

Across curved switching interfaces, the jump-canceling direction rotates with
the iterate. The identity $\widehat B_ks_k=0$ of~\eqref{art:hessian-update}
makes moving comparison geometry viable: writing $\Pi_s=I-ss^T/(s^Ts)$ for the
orthogonal projector onto $s^\perp$, we have
$\widehat B_k=\Pi_{s_k}\widehat B_k\Pi_{s_k}\succeq0$, so any component of field
motion parallel to the step $s_k$ is annihilated for free. Only transverse
distortion perpendicular to $s_k$ increases the energy; when that factor is
$1+O(d_k)$ per step, the decrease in $F$ controls the distortion via
$\sum_kd_k<\infty$ and keeps the cumulative multiplicative factor bounded.

To state this purely in terms of the objective $F$, let $x,z\in\R^n$ be
differentiable points and write $g=\nabla F(x)$, $h=\nabla F(z)$, $s=z-x$,
$y=h-g$, $d=-g^Ts$, $\Delta=F(x)-F(z)$, and $t=y^Ts$. We call the ordered pair
$(x,z)$ \emph{eligible} when
\[
 d>0,\qquad c_1d\le\Delta\le d,\qquad t\ge\omega d.
\]
Every consecutive pair $(x_k,x_{k+1})$ of an admissible run is eligible. On
each positive sublevel band, a comparison field needs to detect the optimality
gap, bound the projected gradient jump by the secant curvature, and bound the
cost of changing the field by the predicted decrease; we refer to these three
requirements as \emph{visibility}, \emph{directional secant control}, and
\emph{projected transport}.

\begin{definition}[Band certificate]\label{art:band-certificate}
Let $K=\{x:F(x)\le b\}$ be compact, with $b>F_*=\min F$.
For each nonempty band $K_\eta=\{x\in K:F(x)\ge F_*+\eta\}$ ($\eta>0$), a
\emph{band certificate} consists of a symmetric positive-semidefinite matrix
field $W_\eta(x)\succeq0$ defined at the differentiable points of $K_\eta$ and
constants $m_\eta>0$, $C_\eta\ge0$, $L_\eta\ge0$, $\delta_\eta>0$ such that
\begin{align}
 g^TW_\eta(x)g&\ge m_\eta,\label{art:visibility}\\
 y^TW_\eta(z)y&\le C_\eta t,\label{art:secant-control}\\
 \Pi_sW_\eta(z)\Pi_s&\preceq
 (1+L_\eta d)\Pi_sW_\eta(x)\Pi_s,\label{art:transport}
\end{align}
where~\eqref{art:visibility} holds at every differentiable $x\in K_\eta$
and~\eqref{art:secant-control}--\eqref{art:transport} hold for all eligible
pairs $(x,z)$ in $K_\eta$ with $d\le\delta_\eta$.
\end{definition}

Here $W_\eta(x)$ need not be continuous, uniformly bounded, or full rank, and
all constants may degenerate as $\eta\downarrow0$, allowing arbitrary
nonsmoothness at the minimizer (setting $W_\eta(x)=(x-x_*)(x-x_*)^T$ with
$L_\eta=0$ recovers the radial case of Corollary~\ref{art:radial-convergence}).

\begin{theorem}[Convergence from band certificates]\label{art:certificate}
If $K=\{F\le b\}$ is compact and admits a band certificate on every nonempty
band $K_\eta$ ($\eta>0$), then every admissible run starting at $x_0\in K$
satisfies
\[
 F(x_k)\longrightarrow F_*,\qquad
 \dist(x_k,\argmin F)\longrightarrow0.
\]
In particular, if $F$ has a unique minimizer $x_*$, then $x_k\to x_*$.
\end{theorem}
\begin{proof}
Armijo keeps $x_k\in K$ with $F(x_k)$ decreasing to some limit $\ell\ge F_*$, so
Lemma~\ref{art:run-invariants} gives $\sum_kd_k<\infty$, and subgradient
boundedness on $K$ (proved in Theorem~\ref{art:radial-equivalence}) ensures that
Theorem~\ref{art:reciprocal-growth} applies.

Suppose toward a contradiction that $\ell>F_*$, and set $\eta=(\ell-F_*)/2>0$.
Since $F(x_k)\ge\ell>F_*+\eta$, every iterate $x_k$ lies in $K_\eta$. Using the
certificate at level $\eta$, choose $J$ so that $d_k\le\delta_\eta$ for all
$k\ge J$; suppressing the subscript $\eta$, set $E_k=\tr(B_kW(x_k))\ge0$. Since
$\widehat B_k=\Pi_{s_k}\widehat B_k\Pi_{s_k}\succeq0$ and $\tr(AQ)=\tr(Q^{1/2}AQ^{1/2})\ge0$
for $A,Q\succeq0$, combining~\eqref{art:hessian-update}
with~\eqref{art:visibility}--\eqref{art:transport} gives, for all $k\ge J$,
\begin{align*}
 E_{k+1}
 &=\tr(\widehat B_kW(x_{k+1}))+\frac{y_k^TW(x_{k+1})y_k}{t_k}\\
 &\le(1+Ld_k)\tr(\widehat B_kW(x_k))+C
 \le(1+Ld_k)\left(E_k-\frac{m}{q_k}\right)+C.
\end{align*}
Setting $P_J=1$ and $P_{k+1}=(1+Ld_k)P_k$, summability of $d_k$ gives
$1\le P_k\le\exp(L\sum_{j\ge J}d_j)=:\overline P<\infty$. Dividing by $P_{k+1}$
and summing from $k=J$ to $N-1$ (as in standard multiplicative quasi-Fej\'er
estimates; cf.~\cite[Lemmas~2.2--2.3]{FND_CV13}) with $E_N\ge0$ yields
\[
 m\overline P^{-1}\sum_{k=J}^{N-1}q_k^{-1}\le E_J+C(N-J),
\]
contradicting~\eqref{art:reciprocal-limits}. Thus $F(x_k)\to F_*$. By
compactness of $K$ and continuity of $F$, every cluster point of $\{x_k\}$ lies
in $\argmin F$, so $\dist(x_k,\argmin F)\to0$ (and $x_k\to x_*$ when
$\argmin F=\{x_*\}$).\qed
\end{proof}

When a rank-one comparison vector field $w(x)$ is constructed only locally near
the cluster set of a run rather than across an entire sublevel band, as in the
planar contour analysis below, the same tail normalization rules out
nonstationary limits directly.

\begin{lemma}[Tail vector comparison]\label{art:tail-comparison}
Let $\{x_k\}$ be an admissible run with bounded gradients and $\sum_kd_k<\infty$.
There cannot exist an index $J$, vectors $w_k\in\R^n$ for $k\ge J$, scalars
$\lambda_k,\beta_k\in\R$, and constants $m>0$, $C,L\ge0$ such that for all
$k\ge J$,
\[
 \begin{gathered}
 (g_k^Tw_k)^2\ge m,\qquad
 \frac{(y_k^Tw_{k+1})^2}{t_k}\le C,\\
 w_{k+1}=\lambda_kw_k+\beta_ks_k,\qquad
 \lambda_k^2\le1+Ld_k.
 \end{gathered}
\]
\end{lemma}
\begin{proof}
Set $E_k=w_k^TB_kw_k\ge0$ for $k\ge J$. Since $\widehat B_ks_k=0$ and
$\widehat B_k\succeq0$, evaluating~\eqref{art:hessian-update} at
$w_{k+1}=\lambda_kw_k+\beta_ks_k$ gives
\begin{align*}
 E_{k+1}
 &=\lambda_k^2(w_k^T\widehat B_kw_k)+\frac{(y_k^Tw_{k+1})^2}{t_k}
  =\lambda_k^2\left(E_k-\frac{(g_k^Tw_k)^2}{q_k}\right)+\frac{(y_k^Tw_{k+1})^2}{t_k}\\
 &\le(1+Ld_k)\left(E_k-\frac{m}{q_k}\right)+C.
\end{align*}
Dividing by $P_{k+1}=\prod_{j=J}^k(1+Ld_j)\le\exp(L\sum_{j\ge J}d_j)<\infty$ and
summing over $J\le k<N$ bounds $\sum_{k=J}^{N-1}q_k^{-1}$ linearly in $N-J$,
contradicting Theorem~\ref{art:reciprocal-growth}.\qed
\end{proof}

\subsection{Contour-velocity fields and affine-critical cluster points in the plane}

To instantiate Lemma~\ref{art:tail-comparison} across curved switching curves, a
vector field $w(x)$ must pair equally with both active gradients across each
curve while changing slowly enough along BFGS steps to satisfy
$\lambda_k^2\le1+Ld_k$. In $\R^2$, under strict contours and local two-branch regularity we can build
such a field on a neighborhood of any compact subset of a positive level curve
disjoint from the \emph{affine-critical set}
\[
 \mathcal C_{\rm aff}=\{p\in\R^2:0\in\aff\partial F(p)\}.
\]
Recall from Section~\ref{art:geometry} why $\mathcal C_{\rm aff}$ is the natural
critical set for equal-rate comparison directions: at any $p\notin\mathcal C_{\rm aff}$,
the minimum-norm vector $a\in\aff\partial F(p)$ is nonzero and satisfies
$g^Ta=\norm{a}^2>0$ for all $g\in\partial F(p)$, so $-a/\norm{a}^2$ decreases
every active supporting model at unit rate.

The geometric construction on a compact non-affine-critical subset of a level
curve $\{F=\ell\}$ ($\ell>F_*$) proceeds in three steps (proved in full in
Lemma~\ref{art:planar-field} of Appendix~\ref{art:planar-appendix}):
\begin{enumerate}
\item \emph{Corner velocities cancel jumps.} Where two convex $C^2$ branches
$f_1,f_2$ tie with linearly independent gradients, the switching curve is
parametrized by the objective level $f_1(p_j(t))=f_2(p_j(t))=t$. As the level
$t$ increases, both active branches change at unit rate along the corner
trajectory $p_j(t)$. Differentiating in $t$ therefore gives a velocity vector
$v_j(t)=p_j'(t)$ satisfying $\nabla f_1(p_j(t))^Tv_j(t)=\nabla f_2(p_j(t))^Tv_j(t)=1$:
both active gradients have unit inner product with $v_j(t)$, so their difference
$\nabla f_2-\nabla f_1$ vanishes in that direction.
\item \emph{Normal-angle extension and tangent preservation.} When positive
sublevels $\{F\le t\}$ contain no boundary line segments, distinct corners on
$\{F=t\}$ have disjoint supporting-normal arcs. We interpolate the positive
corner pairings smoothly across unit-normal angles $\theta$ into a positive
periodic support height $h(t,\theta)$ and define $w(x)=hn+h_\theta T$ along the
intervening smooth arcs (with unit normal $n=\nabla F/\norm{\nabla F}$ and
tangent $T$), which collapses to $w(x)=v_j(F(x))$ across the $j$th switching
curve. Differentiating $hn+h_\theta T$ with respect to $\theta$ cancels the
normal component~\cite{geom_Bogosel2023}, so $Dw$ maps the one-dimensional
level-set tangent line into itself: $\det(s,Dw\,s)=0$ whenever $s$ is tangent to
$\{F=\text{const}\}$.
\item \emph{Bounding field change by functional decrease.} Weak Wolfe keeps
short eligible steps uniformly transverse to $w$. Because $Dw$ preserves
level-set tangents, the transverse change in $w$ along an eligible step $s$
comes only from motion across level curves and is bounded by the variation of
$F$ along the chord, which convexity controls by $2d$. Expressing
$w(z)=\lambda w(x)+\beta s$ in the basis $(w(x),s)$ then yields
$\lambda^2\le1+Ld$ (Appendix~\ref{art:planar-appendix}), so the decrease in $F$
bounds the cost of changing the field without requiring $\sum_k\norm{s_k}<\infty$.
\end{enumerate}

\begin{theorem}[Affine-critical cluster points in the plane]\label{art:planar-cluster}
Let $F:\R^2\to\R$ be finite and convex, attain its minimum $F_*$, and have a
compact sublevel $K=\{x:F(x)\le b\}$ with $b>F_*$. Assume that:
\begin{enumerate}
\item for every $F_*<t\le b$, the boundary of $\{F\le t\}$ contains no
nontrivial line segment; and
\item at every $p\in K\setminus\mathcal C_{\rm aff}$, either $F$ is $C^2$ near
$p$, or $F=\max\{f_1,f_2\}$ near $p$ for convex $C^2$ functions $f_1,f_2$ both
active at $p$ with $\nabla f_1(p)$ and $\nabla f_2(p)$ linearly independent.
\end{enumerate}
Then every infinite admissible full-BFGS run starting in $K$ has at least one
cluster point in $\mathcal C_{\rm aff}$ (with no local structure assumed at
points of $\mathcal C_{\rm aff}$ itself).
\end{theorem}
\begin{proof}
Along the run, $x_k\in K$, $F(x_k)\downarrow\ell\ge F_*$, gradients are
uniformly bounded on $K$, and $\sum_kd_k<\infty$. Let $\Gamma\subset\{F=\ell\}$
be the nonempty compact cluster set of $\{x_k\}$. If $\ell=F_*$, then
$\Gamma\subseteq\argmin F\subseteq\mathcal C_{\rm aff}$. Suppose for
contradiction that $\Gamma\cap\mathcal C_{\rm aff}=\emptyset$; then $\ell>F_*$.

First, $\norm{s_k}\to0$. Otherwise compactness yields a subsequence with
$x_k\to x$, $x_{k+1}\to z\ne x$, and $d_k\to0$. For every $\theta\in[0,1]$,
\[
 F(x_k)-\theta d_k
 \le F((1-\theta)x_k+\theta x_{k+1})
 \le(1-\theta)F(x_k)+\theta F(x_{k+1}),
\]
so passing to the limit gives $F\equiv\ell$ on the segment $[x,z]$. Since
$\ell>F_*$, convexity forces every point of $\{F=\ell\}$ to lie on the boundary
of $\{F\le\ell\}$ (an interior point of $\{F\le\ell\}$ with value $\ell$ would
lie between a minimizer and a point of $\{F\le\ell\}$, forcing $F<\ell$),
contradicting hypothesis~1.

Applying Lemma~\ref{art:planar-field} (Appendix~\ref{art:planar-appendix}) to
the compact set $\Gamma\subset\{F=\ell\}\setminus\mathcal C_{\rm aff}$ supplies
an open neighborhood $U\supset\Gamma$, a vector field $w$, and constants
$\eta,\delta>0$, $C_0,L\ge0$. Since $\dist(x_k,\Gamma)\to0$, $\norm{s_k}\to0$,
and $d_k\to0$, there exists $J$ such that for all $k\ge J$, $x_k,x_{k+1}\in U$,
$\norm{s_k}<\delta$, and $d_k\le1$; Lemma~\ref{art:planar-field} then gives
$w_k=w(x_k)$ satisfying $(g_k^Tw_k)^2\ge\eta^2$, $(y_k^Tw_{k+1})^2/t_k\le C_0$,
and $w_{k+1}=\lambda_kw_k+\beta_ks_k$ with $\lambda_k^2\le1+Ld_k$, contradicting
Lemma~\ref{art:tail-comparison}.\qed
\end{proof}

\begin{corollary}[Planar convergence when affine criticality implies optimality]\label{art:planar-convergence}
Under the hypotheses of Theorem~\ref{art:planar-cluster}, if
$\mathcal C_{\rm aff}\cap K\subseteq\argmin F$, then every admissible run
starting in $K$ satisfies $F(x_k)\to F_*$ and $\dist(x_k,\argmin F)\to0$; if $F$
has a unique minimizer $x_*$, then $x_k\to x_*$.
\end{corollary}
\begin{proof}
Theorem~\ref{art:planar-cluster} provides a cluster point $p\in\mathcal C_{\rm aff}\cap K\subseteq\argmin F$.
Monotonicity of $F(x_k)$ and continuity of $F$ give $\lim_kF(x_k)=F(p)=F_*$, and
compactness of $K$ yields $\dist(x_k,\argmin F)\to0$ (hence $x_k\to x_*$ under
uniqueness).\qed
\end{proof}

\subsection{Nondegenerate planar smooth maxima and an explicit nonradial basin}

For a finite smooth maximum $F=\max_{1\le i\le m}f_i$ in $\R^2$ with a unique
minimizer $x_*$, both the isolation condition $\mathcal C_{\rm aff}=\{x_*\}$
near $x_*$ and two-branch regularity away from $x_*$ follow automatically from
standard multiplier nondegeneracy at $x_*$, namely affine independence of the
active gradients together with strictly positive optimal convex weights.

\begin{corollary}[Nondegenerate planar smooth maxima]\label{art:planar-max}
Let $F:\R^2\to\R$ be finite and convex with a unique minimizer $x_*$, and
suppose $F=\max_{1\le i\le m}f_i$ near $x_*$ for convex $C^2$ branches $f_i$.
Assume that sufficiently low positive sublevels have no boundary line segments
(for instance, if the branches active at $x_*$ share a positive lower Hessian
bound near $x_*$), and that the active index set $I_*=\{i:f_i(x_*)=F(x_*)\}$ has
affinely independent gradients $\{\nabla f_i(x_*)\}_{i\in I_*}$ admitting
strictly positive weights
\[
 \lambda_i^*>0\quad(i\in I_*),\qquad \sum_{i\in I_*}\lambda_i^*=1,
 \qquad \sum_{i\in I_*}\lambda_i^*\nabla f_i(x_*)=0.
\]
Then there exists $\varepsilon>0$, independent of $H_0\succ0$ and of
$0<c_1<c_2<1$, such that every infinite admissible full-BFGS run with
$F(x_0)<F(x_*)+\varepsilon$ satisfies $F(x_k)\to F(x_*)$ and $x_k\to x_*$.
\end{corollary}
\begin{proof}
By continuity, restrict to a neighborhood of $x_*$ on which all $i\notin I_*$
are strictly inactive. Form the $3\times|I_*|$ matrix $A(x)$ with columns
$(\nabla f_i(x),1)^T$ ($i\in I_*$). Affine independence at $x_*$ means $A(x_*)$
has full column rank, which persists nearby, where any solution of
$A(p)\lambda=(0,0,1)^T$ is uniquely given by the continuous formula
$\lambda(p)=(A(p)^TA(p))^{-1}A(p)^T(0,0,1)^T$ with $\lambda(x_*)=\lambda^*>0$.
Shrink the neighborhood so that every coordinate of $\lambda(p)$ stays strictly
positive. Because $\partial F(p)=\conv\{\nabla f_i(p):i\text{ active at }p\}$
(each active gradient supports $F$ globally by convexity, while any vector
outside their convex hull is strictly separated from it, violating the
directional-derivative formula $F'(p;d)=\max_{i\text{ active}}\nabla f_i(p)^Td$),
any $p\in\mathcal C_{\rm aff}$ in this neighborhood has $0$ as an affine
combination of active gradients; extending by $0$ on inactive indices in $I_*$
forces that weight vector to equal $\lambda(p)>0$, so all $i\in I_*$ are active
and $0\in\conv\partial F(p)$, making $p=x_*$.

Moreover, at most three gradients in $\R^2$ can be affinely independent ($|I_*|\le3$).
At any $p\ne x_*$ in this neighborhood, $p\notin\mathcal C_{\rm aff}$, so three
branches cannot be simultaneously active (their affine hull would be all of
$\R^2$, containing $0$); one active branch makes $F$ locally $C^2$, and two
active branches must have linearly independent gradients (otherwise the line
through them would pass through $0$, putting $p\in\mathcal C_{\rm aff}$).
Applying Lemma~\ref{art:localization} and Corollary~\ref{art:planar-convergence}
gives $\varepsilon>0$. Finally, if the branches $i\in I_*$ satisfy $\nabla^2f_i\succeq\mu I\succ0$ on
a convex neighborhood of $x_*$, then $f_i(x)-\tfrac\mu2\norm{x}^2$ is convex for
each $i\in I_*$, so $F$ is $\mu$-strongly convex there and the boundaries of its
positive sublevels inside that neighborhood contain no nontrivial line
segments.\qed
\end{proof}

This multiplier nondegeneracy matches the active-set regularity assumption of
local bundle-Newton and piecewise-smooth analyses (Gebken~\cite{BG}; Lewis and
Wylie~\cite[Propositions~5.5--5.6]{LW}), with a key difference:
Corollary~\ref{art:planar-max} proves convergence of unmodified BFGS from
objective geometry alone, with no assumed iterate or spectral asymptotics. We
conclude the curved-interface theory with an explicit three-quadratic maximum
showing that moving comparison fields strictly reach beyond radial analysis.

\begin{theorem}[Explicit convergence basin with curved, nonradial switching]\label{art:nonradial-example}
For the three-quadratic maximum on $\R^2$ given by
\[
 \begin{aligned}
 f_0(u,v)&=u+u^2+\tfrac12v^2,\qquad
 f_1(u,v)=v+\tfrac12u^2+\tfrac32v^2,\\
 f_2(u,v)&=-u-v+\tfrac32u^2+v^2,\qquad F(u,v)=\max_{0\le i\le2}f_i(u,v),
 \end{aligned}
\]
every infinite admissible full-BFGS run with $F(x_0)<1/100$, arbitrary
$H_0\succ0$, and any fixed $0<c_1<c_2<1$ satisfies $F(x_k)\to0$ and
$x_k\to(0,0)$. Nevertheless, the radial secant condition~\eqref{art:radial-condition}
fails on a positive sublevel band inside every neighborhood of $(0,0)$, even
when restricted to pairs eligible for $(c_1,c_2)=(1/4,1/2)$.
\end{theorem}
\begin{proof}
Each branch has Hessian between $I$ and $3I$, and $(f_0+f_1+f_2)/3=u^2+v^2$, so
$F(u,v)\ge u^2+v^2$, $(0,0)$ is the unique minimizer with $F_*=0$, and all
positive sublevels are strictly convex. At $(0,0)$, all three branches are
active with gradients $(1,0)^T$, $(0,1)^T$, $(-1,-1)^T$ and weights $1/3$.
To certify the explicit basin $\{F<1/100\}$, solve $A(u,v)\lambda=(0,0,1)^T$
explicitly: with $D=\det A(u,v)=3(1+v-uv)$,
\[
 D\lambda_0=1-4u+3v-7uv,\quad
 D\lambda_1=1+2u-3v-uv,\quad
 D\lambda_2=1+2u+3v+5uv.
\]
On the box $|u|,|v|\le1/10$, we have $D\ge267/100$ and the three numerators are
at least $23/100$, $49/100$, and $45/100$, so $\lambda_i(u,v)>0$ everywhere on
that box. Since $\{F\le b\}\subset\{u^2+v^2\le b\}\subset[-1/10,1/10]^2$ for any
$b\in(F(x_0),1/100)$, the argument of Corollary~\ref{art:planar-max} shows that
$\mathcal C_{\rm aff}\cap\{F\le b\}=\{(0,0)\}$ and two-branch regularity holds
throughout $\{F\le b\}\setminus\{(0,0)\}$; Corollary~\ref{art:planar-convergence}
gives $F(x_k)\to0$ and $x_k\to(0,0)$.

To exhibit the radial obstruction, fix small $u>0$ and let
\[
 v=v(u)=\frac{-1+\sqrt{1+4u+2u^2}}2=u-\tfrac12u^2+O(u^3)
\]
solve $f_0(u,v)=f_1(u,v)$ with $f_2(u,v)<f_0(u,v)$. Along this switching curve,
\[
 (\nabla f_0-\nabla f_1)^T(u,v)
 =(1+u)u-(1+2v)v=\tfrac12u^2-v^2=-\tfrac12u^2+O(u^3)\ne0
\]
for small $u>0$. Taking $x=(u,v+e)$ (uniquely on branch~$1$) and $z=(u,v-e)$
(uniquely on branch~$0$) as $e\downarrow0$, direct evaluation gives
\[
 d=2e(1+3v+3e),\qquad t=2e(1+2v+4e),\qquad \Delta=e(1+4v)+e^2.
\]
As $u,e\to0^+$, $\Delta/d\to1/2$ and $t/d\to1$, so for fixed small $u>0$ and all
sufficiently small $e>0$ the pair $(x,z)$ is eligible for $(c_1,c_2)=(1/4,1/2)$
and stays above half the positive tie value $f_0(u,v)>0$. Yet
$y^Tz\to\tfrac12u^2-v^2\ne0$ while $t\to0^+$, so $(y^Tz)^2/t\to+\infty$ on a
positive band inside any given neighborhood of $(0,0)$.\qed
\end{proof}

\subsection{Verifying the hypotheses and existence of admissible runs}

In applying Table~\ref{art:positive-outcomes} to a given convex objective $F$,
all algorithmic side conditions hold automatically along any admissible run on a
compact sublevel: weak Wolfe preserves $H_k\succ0$, Armijo gives
$\sum_kd_k<\infty$, and finite convexity bounds all subgradients on $K$.
Verification therefore reduces to a concrete geometric check on $F$:
\begin{itemize}
\item \emph{Smooth-plus-homogeneous route (Theorem~\ref{art:centered}):} exhibit
a decomposition $F(x)=\psi(x)+\gamma(x-x_*)$ around a minimizer $x_*$ on an open
convex domain $U$ containing $\{F\le F(x_0)\}$ (or locally around a unique
minimizer $x_*$ via Corollary~\ref{art:local-centered}).
\item \emph{Polyhedral route (Theorem~\ref{art:polyhedral}, Corollary~\ref{art:plq}):}
verify uniqueness of $x_*$ and locally $C^{1,1}$ cell formulas on a finite
polyhedral cover near $x_*$.
\item \emph{Planar smooth-maximum route (Corollary~\ref{art:planar-max}):} at
the unique minimizer $x_*$, check that the active gradients are affinely
independent with strictly positive convex multipliers $\lambda_i^*>0$, together
with a positive lower Hessian bound (or strict contours) on the active branches.
\end{itemize}

Appendix~\ref{art:hypothesis-examples} illustrates these certificates with an
explicit degenerating BFGS trajectory on $u^2+|v|$
(Proposition~\ref{art:examples-mixed}) and a full-rank calibration in $\R^n$
(Proposition~\ref{art:examples-full-rank}). A convergence guarantee for
admissible runs does not by itself guarantee that a line search can take a
differentiable step along a given ray. In Appendix~\ref{art:hypothesis-examples},
Lemma~\ref{art:examples-run-existence} proves that on any bounded-below finite
convex objective having no differentiable stationary point and at most countably
many nondifferentiable points along each line through a differentiability point,
infinite admissible runs exist from every differentiable $x_0$, $H_0\succ0$, and
$0<c_1<c_2<1$, covering both a parameterized curved two-quadratic family
(Proposition~\ref{art:examples-two-quadratics}) and the three-quadratic
objective of Theorem~\ref{art:nonradial-example}
(Corollary~\ref{art:examples-nonradial-runs}); Proposition~\ref{art:examples-no-step}
complements this with a three-piece coercive cone whose initial search ray has
no differentiable weak-Wolfe step.

\part{Self-Similar Cycles and Nonstationary Convergence}\label{part:failure}
\section{The Affine-Critical Three-Slope Cycle and Global Polyhedral Failure}
\label{art:failure}

Using only the unmodified BFGS protocol of Section~\ref{art:setup}, Part~\ref{part:failure}
constructs global and local convex objectives on which unmodified full BFGS from
$H_0=I$ under a standard deterministic line search converges to nonstationary
points, showing that convexity, coercivity, and smooth or strongly convex branch
formulas cannot replace the structural conditions of Part~\ref{part:convergence}.

The construction across Sections~\ref{art:failure} and~\ref{art:failure-local}
proceeds in three stages:
\begin{enumerate}
\item First, in this section, at a nonstationary three-slope corner in $\R^2$
where $0\in\aff\{g_0,g_1,g_2\}\setminus\conv\{g_0,g_1,g_2\}$, we construct three
BFGS steps after which the iterate and inverse-Hessian matrix are exact scalar
copies of their initial values; the resulting run from $z_0$ with $H_0=I$
converges to a nonstationary point of a coercive four-piece polyhedral objective
(Theorem~\ref{art:failure_polyhedral}).
\item Second, in Section~\ref{art:failure-local}, we place shrinking copies of
that three-step cycle near the unique minimizer of a single four-branch
$C^{2,1}$ convex objective under the scaling relation $\sigma=\lambda^2$, which
leaves the initial matrix $H_0=I$ unchanged across every scale
(Theorem~\ref{art:local-failure}).
\item Third, also in Section~\ref{art:failure-local}, because uniform strong
convexity breaks exact scaling, we analyze a scale-normalized three-step return
map and extract a perturbed failing orbit on its local stable manifold,
embedding failing runs into a four-branch uniformly strongly convex $C^{1,1}$
objective (Theorem~\ref{art:local-failure-lower}).
\end{enumerate}

\subsection{The line-search implementation}

To establish genuine algorithmic counterexamples rather than adversarially
chosen valid steps, we fix the standard deterministic unit-first bracketing line
search of Lewis and Overton~\cite[Algorithm~4.6]{LO}. At a differentiable
iterate $x$ with $g=\nabla F(x)$ and search direction $p=-Hg$, the search
initializes $\alpha_{\rm lo}=0$, $\alpha_{\rm hi}=+\infty$, and initial trial
step $\alpha=1$, testing Armijo first and weak Wolfe second:
\begin{equation}\label{art:failure_tests}
 F(x+\alpha p)\le F(x)+c_1\alpha g^Tp,
 \qquad \nabla F(x+\alpha p)^Tp\ge c_2g^Tp.
\end{equation}
If Armijo fails, the routine sets $\alpha_{\rm hi}=\alpha$; if Armijo holds and
weak Wolfe fails, it sets $\alpha_{\rm lo}=\alpha$; the next trial is
$(\alpha_{\rm lo}+\alpha_{\rm hi})/2$ when $\alpha_{\rm hi}<\infty$ and
$2\alpha_{\rm lo}$ otherwise, accepting as soon as both inequalities
in~\eqref{art:failure_tests} hold. Along every run constructed below, every
accepted and rejected trial point $x+\alpha p$ lies in the open smooth interior
of a single active piece, and every inequality in~\eqref{art:failure_tests}
holds or fails strictly, so strict and non-strict implementations execute the
identical decision tree.

We write one unmodified inverse-BFGS update as the operator
\begin{equation}\label{art:failure_update}
 \mathcal B(H,g,g^+,\alpha)
 =\left(I-\frac{sy^T}{y^Ts}\right)H
  \left(I-\frac{ys^T}{y^Ts}\right)+\frac{ss^T}{y^Ts},
\end{equation}
where $s=-\alpha Hg$ and $y=g^+-g$, with no resetting or rescaling of $H$.

\subsection{A self-similar three-step cycle and global polyhedral failure}

Consider a three-piece positively homogeneous cone
$G(z)=\max_{0\le i\le2}g_i^Tz$ on $\R^2$. If three consecutive BFGS steps from
$H_0=I$ visit smooth regions of branches $0,1,2,0$ in turn and return both the
inverse Hessian and the iterate to the same scalar multiple $\rho\in(0,1)$ of
their starting values,
\[
 H_3=\rho I,\qquad z_3=\rho z_0,
\]
then scale invariance of $G$ locks the algorithm into an infinite geometric
contraction $(x_{3m+j},H_{3m+j})=(\rho^mz_j,\rho^mM_j)$. Because $\rho<1$, the steps $\norm{s_k}=O(\rho^{k/3})$ are summable and the
iterates converge to the origin, even when $0\notin\conv\{g_0,g_1,g_2\}$ and the
origin is therefore nonstationary.

Counting degrees of freedom shows why three slopes in $\R^2$ just suffice. Fix
$H_0=I$ and set $g_0=(1,0)^T$ and $g_1=(-1,-1/8)^T$, so a unit step $\alpha_0=1$
produces the rational matrix $M_1=\mathcal B(I,g_0,g_1,1)$. The third slope
$g_2=(c_\circ,d_\circ)^T$ supplies two real parameters, which we tune so that
the dyadic bisection step lengths $\alpha_1=1/4$ and $\alpha_2=1/16$ close the
matrix cycle at a scalar multiple $M_3=\rho I$ (two equations for the symmetric
$2\times2$ matrix $M_3$ to be proportional to $I$). Once $g_2$ and $\rho$ are
fixed, the displacement closure $z_0+s_0+s_1+s_2=\rho z_0$ uniquely determines
the starting point $z_0$.

Specifically, let $(c_\circ,d_\circ)$ be the algebraic root isolated to within
$10^{-30}$ by the rational contraction certificate in
Appendix~\ref{art:seed_certificate}, which satisfies $0.318<c_\circ<0.319$ and
$-1.026<d_\circ<-1.025$. Define
\[
 g_0=(1,0)^T,\quad g_1=(-1,-1/8)^T,\quad g_2=(c_\circ,d_\circ)^T,
 \quad G(z)=\max_{0\le i\le2}g_i^Tz,
\]
and set the contraction factor
\begin{equation}\label{art:seed_rho}
 \rho=\frac{129(c_\circ-8d_\circ)}{(16c_\circ-127d_\circ-16)(127c_\circ+16d_\circ+129)},
 \qquad 0.060<\rho<0.061.
\end{equation}
Define the cycle matrices, steps, and base iterates by
\begin{align}
 M_0&=I,\quad
 M_1=\begin{pmatrix}129/256&-1/16\\-1/16&1\end{pmatrix},
 &M_2&=\mathcal B(M_1,g_1,g_2,1/4),\label{art:seed_matrices}\\
 s_0&=(-1,0)^T,\quad s_1=(127/1024,1/64)^T,
 &s_2&=\rho(1-c_\circ,-d_\circ)^T,\label{art:seed_steps}
\end{align}
and
\begin{equation}\label{art:seed_points}
 z_0=\frac1{1-\rho}
 \begin{pmatrix}897/1024-\rho(1-c_\circ)\\ \rho d_\circ-1/64\end{pmatrix},
 \qquad z_{j+1}=z_j+s_j\quad(j=0,1,2).
\end{equation}

\begin{lemma}[Exact cycle closure and strict trial margins]\label{art:seed_cycle}
The matrices $M_0,M_1,M_2$ in~\eqref{art:seed_matrices} are positive definite,
and
\[
 \mathcal B(I,g_0,g_1,1)=M_1,\qquad
 \mathcal B(M_2,g_2,g_0,1/16)=\rho I,\qquad z_3=\rho z_0.
\]
At $z_j$ ($j=0,1,2$), the unique active gradient of $G$ is $g_j$. With
$p_j=-M_jg_j$, for every $0<c_1\le1/3$ and $c_2>0$ the complete line-search
query sequence is given by Table~\ref{art:trial_table}: every rejected trial
strictly exceeds the current objective value $G(z_j)$, and every accepted step
satisfies the strict margins
\begin{equation}\label{art:seed_margins}
 G(z_j)+\tfrac13g_j^Ts_j-G(z_{j+1})>0,\qquad
 g_{j+1}^Ts_j>0>g_j^Ts_j,
\end{equation}
where $g_3=g_0$. Moreover, all nine trial points and $z_0,z_1,z_2$ have a unique
active branch and satisfy, with $[z]_2$ denoting the second coordinate,
\begin{equation}\label{art:seed_query_bounds}
 |g_i^Tz|<3\quad(i=0,1,2),\qquad [z]_2<2,\qquad G(z)>0.
\end{equation}
\end{lemma}

\begin{table}[tbp]
\centering
\setlength{\tabcolsep}{5.8pt}
\begin{tabular}{cccccl}
\toprule
Phase $j$ & \shortstack{$10^4G$\\current} & Trial $\alpha$
 & \shortstack{$10^4G$\\trial} & Active branch & Decision\\
\midrule
0 & $(8884,8885)$ & $1$ & $(1218,1219)$ & 1 & accept\\
\midrule
1 & $(1218,1219)$ & $1$ & $(3845,3846)$ & 0 & reject\\
  &               & $1/2$ & $(1364,1365)$ & 0 & reject\\
  &               & $1/4$ & $(723,724)$ & 2 & accept\\
\midrule
2 & $(723,724)$ & $1$ & $(6687,6691)$ & 0 & reject\\
  &             & $1/2$ & $(3406,3408)$ & 0 & reject\\
  &             & $1/4$ & $(1765,1767)$ & 0 & reject\\
  &             & $1/8$ & $(944,946)$ & 0 & reject\\
  &             & $1/16$ & $(534,536)$ & 0 & accept\\
\bottomrule
\end{tabular}
\caption{Complete line-search query tree across one three-step cycle, with
rigorous rational enclosures for $10^4G(z_j)$ (``current'') and
$10^4G(z_j+\alpha p_j)$ (``trial''). Every rejected trial strictly raises $G$
and therefore fails Armijo for every $c_1>0$; every accepted trial satisfies the
strict Armijo and weak-Wolfe margins~\eqref{art:seed_margins}, making the
decisions identical for all $0<c_1\le1/3$ and $c_2>0$.}
\label{art:trial_table}
\end{table}

\begin{figure}[tbp]
\centering
\includegraphics[width=\textwidth]{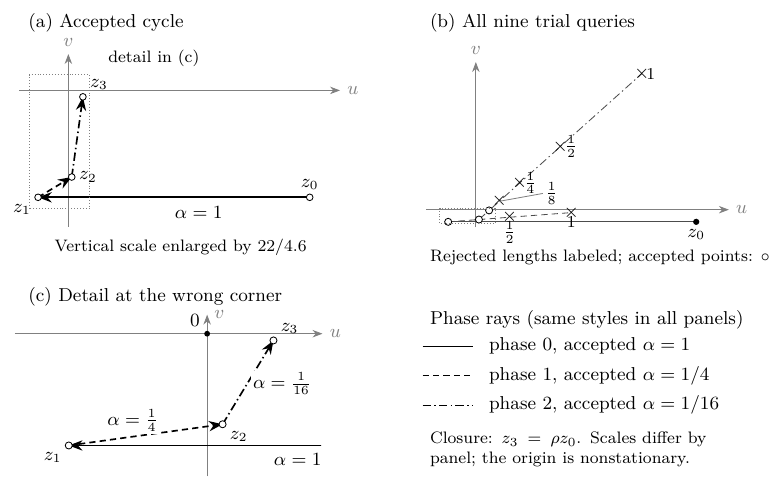}
\caption{Geometry of the self-similar three-step cycle. Solid, dashed, and
dash--dot rays show phases $0,1,2$ with accepted step lengths $1,1/4,1/16$:
(a)~the accepted cycle $z_0\to z_1\to z_2\to z_3=\rho z_0$; (b)~all nine unit
and bisected trial queries; (c)~enlargement near the nonstationary origin.
Vertical-to-horizontal aspect ratios are $22/4.6$, $2.5/4.2$, and $23/21$ in
(a), (b), and (c); exact rational proofs are in
Appendix~\ref{art:seed_certificate}.}
\label{art:cycle-figure}
\end{figure}

The rational root proof, symbolic closure identity, and interval enclosures for
Lemma~\ref{art:seed_cycle} are given in Appendix~\ref{art:seed_certificate},
while Table~\ref{art:trial_table} and Figure~\ref{art:cycle-figure} make the
line-search dynamics immediately readable: phase~$0$ accepts $\alpha=1$ on its
first try, phase~$1$ overshoots onto branch~$0$ at $\alpha=1,1/2$ (strictly
raising $G$) before bisecting to $\alpha=1/4$ on branch~$2$, and phase~$2$
bisects four times from $\alpha=1$ down to $\alpha=1/16$ on branch~$0$. At each
acceptance, $G(z_{j+1})<G(z_j)+\tfrac13g_j^Ts_j$ and $g_{j+1}^Ts_j>0>g_j^Ts_j$,
so Armijo holds strictly for all $c_1\le1/3$ and weak Wolfe holds strictly for
all $c_2>0$ without a single Wolfe rejection.

By positive homogeneity of $G$, Lemma~\ref{art:seed_cycle} propagates across all
cycles $m\ge0$ for every fixed pair $0<c_1\le1/3$, $c_1<c_2<1$:
\begin{equation}\label{art:cone_orbit}
 x_{3m+j}=\rho^mz_j,\qquad H_{3m+j}=\rho^mM_j\qquad(m\ge0,\;j=0,1,2).
\end{equation}
Indeed, scaling $(z_j,M_j)$ by $\rho^m$ scales every trial point $z_j+\alpha p_j$
by $\rho^m$, leaves active branches and gradients unchanged, scales objective
differences and directional derivatives by $\rho^m>0$, and multiplies the
updated matrix~\eqref{art:failure_update} by $\rho^m$. Thus $x_k\to0$ with
$G(x_k)\downarrow0$, even though along $w=(-1,16)^T$ we have
\[
 g_0^Tw=g_1^Tw=-1,\qquad g_2^Tw=-c_\circ+16d_\circ<-1,
\]
so $G(tw)=-t$ for all $t\ge0$ and the limit $0$ is nonstationary. Capping $G$
far below the trajectory by a fourth affine function makes the objective
coercive without disturbing a single query, yielding our global polyhedral
counterexample.

\begin{theorem}[Global failure on a coercive polyhedral function]\label{art:failure_polyhedral}
Define the four-piece polyhedral convex function on $\R^2$ by
\[
 \overline G(u,v)=\max\{u,\;-u-v/8,\;c_\circ u+d_\circ v,\;v-17\}.
\]
Then $\overline G$ is coercive with unique global minimizer $(-1,16)$ and
minimum value $\overline G(-1,16)=-1$. Nevertheless, starting from $x_0=z_0$ and
$H_0=I$, for every fixed pair $0<c_1\le1/3$, $c_1<c_2<1$, the unit-first
bracketing full-BFGS algorithm generates the trajectory~\eqref{art:cone_orbit},
which converges to the nonstationary point $(0,0)$ with limiting value
$\overline G(0,0)=0>-1$.
\end{theorem}
\begin{proof}
At every base point and trial query of the cone trajectory~\eqref{art:cone_orbit},
the bounds~\eqref{art:seed_query_bounds} and $0<\rho<1$ give
$v-17<2\rho^m-17<0<G(u,v)$, so the fourth affine piece $v-17$ is strictly
inactive throughout every line search and $\overline G$ generates the identical
trajectory~\eqref{art:cone_orbit}. Meanwhile, the convex combination of pieces
$0,1,3$ with weights $8/17,8/17,1/17$ vanishes in $(u,v)$ and equals $-1$
identically, with equality if and only if $u=-u-v/8=v-17=-1$, i.e.,
$(u,v)=(-1,16)$, where $c_\circ(-1)+16d_\circ<-16<-1$. Finally, $\overline G(u,v)\le t$
implies $u\le t$, $v\le t+17$, and $-u-v/8\le t$ (hence $u\ge-9t/8-17/8$ and
$v\ge-8t-8u\ge-16t$), so every sublevel of $\overline G$ is compact.\qed
\end{proof}

Theorem~\ref{art:failure_polyhedral} shows that homogeneity in
Theorem~\ref{art:homogeneous} is sharp: shifting one affine piece by $-17$
separates the nonstationary corner $(0,0)$ from the true minimizer $(-1,16)$ and
defeats global convergence from $H_0=I$.

\section{Local Smooth-Branch and Uniformly Strongly Convex Counterexamples}\label{art:failure-local}

\subsection{Scale covariance (\texorpdfstring{$\sigma=\lambda^2$}{sigma = lambda\^{}2}) and smooth \texorpdfstring{$C^{2,1}$}{C\^{}\{2,1\}} local failure}

Near $(-1,16)$, the polyhedral function $\overline G$ of
Theorem~\ref{art:failure_polyhedral} is locally a centered three-piece cone
(whose local convergence is proved in Theorem~\ref{art:polyhedral}). To refute \emph{local} every-start convergence near the unique minimizer of a
single smooth-branch objective, we must embed infinitely many shrinking copies
of the bad corner at centers $\bar x_\ell\to0$ inside one fixed function $F$,
with every failing run initialized from the same matrix $H_0=I$.

A simple scale-covariance calculation governs when rescaled copies share an
initial inverse-Hessian approximation. Suppose that on an open region $U$
containing a current iterate and all trial queries of a line search, $F$ models
a template function $\mathcal G$ via
\begin{equation}\label{art:matched-local-model}
 F(\bar x+\lambda z)=C+\sigma\mathcal G(z)\qquad(\lambda,\sigma>0).
\end{equation}
Then $g_F=(\sigma/\lambda)g_{\mathcal G}$, so for any trial length $\alpha$,
matching physical steps $s_F=\lambda s_{\mathcal G}$ requires scaling the
inverse-Hessian matrix by $\lambda^2/\sigma$:
\begin{equation}\label{art:matched-scaling}
 H_F=\frac{\lambda^2}{\sigma}H_{\mathcal G},\qquad
 p_F=\lambda p_{\mathcal G},\qquad s_F=\lambda s_{\mathcal G}.
\end{equation}
At an accepted step, $y_F=(\sigma/\lambda)y_{\mathcal G}$ and
$y_F^Ts_F=\sigma y_{\mathcal G}^Ts_{\mathcal G}$, so in~\eqref{art:failure_update}
the projection matrices $I-s_Fy_F^T/(y_F^Ts_F)$ are invariant while the rank-one
term $s_Fs_F^T/(y_F^Ts_F)$ scales by $\lambda^2/\sigma$, preserving
$H_F=(\lambda^2/\sigma)H_{\mathcal G}$ across all subsequent updates. Likewise,
both sides of the Armijo and weak-Wolfe tests~\eqref{art:failure_tests} multiply
by $\sigma>0$, preserving every strict trial decision whenever all queried
points remain inside $U$.

\begin{figure}[tbp]
\centering
\includegraphics[width=\textwidth]{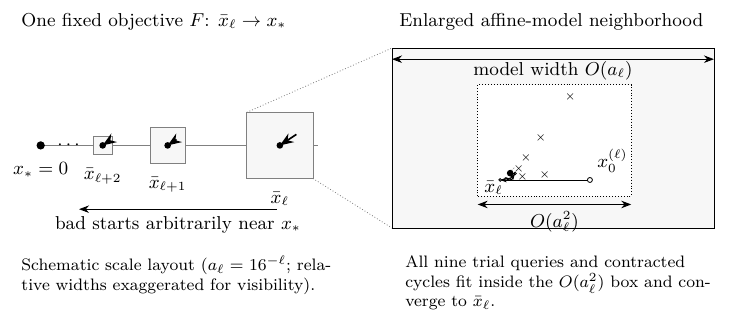}
\caption{Scale separation in the $C^{2,1}$ local counterexample
(Theorem~\ref{art:local-failure}). Inside one fixed objective $F$, nonstationary
corners $\bar x_\ell\to x_*=0$ carry affine-model neighborhoods of width
$O(a_\ell)$, while each copied BFGS trajectory and all of its rejected line-search
trials sit inside a much smaller box of width $O(a_\ell^2)$ and converge to
$\bar x_\ell$.}
\label{art:shrinking-scales}
\end{figure}

Equation~\eqref{art:matched-scaling} isolates the condition
\[
 \sigma=\lambda^2\quad\Longrightarrow\quad H_F=H_{\mathcal G},
\]
which leaves $H_0=I$ unchanged across every scale. To achieve $\sigma=\lambda^2$
together with global $C^{2,1}$ smoothness at the accumulation point $0$, let
$\theta=1/16$ and $a_\ell=\theta^\ell$. Placing scalar affine plateaus of value
$a_\ell^3$ and slope $a_\ell^2$ at nodes $ba_\ell$ gives local slope factor
$a_\ell^2=\sigma/\lambda$; setting spatial scale $\lambda=a_\ell^2$ then yields
$\sigma=a_\ell^4=\lambda^2$. Moreover, each affine plateau has width
$\Theta(a_\ell)$, whereas the entire copied BFGS trial tree has diameter
$\Theta(a_\ell^2)\ll a_\ell$, creating asymptotic room to contain every accepted
and rejected query inside the flat region (Figure~\ref{art:shrinking-scales}).

\begin{lemma}[$C^{2,1}$ scalar plateau spline]\label{art:upper_scalar_lemma}
Let $\theta=1/16$, $a_\ell=\theta^\ell$, and $b_*=2c_\circ-32d_\circ\in(33,34)$.
For each $b\in\{2,b_*\}$ there exists a convex, coercive, nonnegative function
$\varphi_b\in C^{2,1}(\R)$ vanishing only at $0$ such that
\begin{align}
 0&\le\varphi_b''(x)\le\max\{6,102/b\}\qquad(x\in\R),\label{art:scalar_curvature}\\
 \varphi_b(ba_\ell+h)&=a_\ell^3+a_\ell^2h
       \qquad(\ell\ge0,\;|h|\le a_\ell/2),\label{art:scalar_plateau}
\end{align}
and $\varphi_b''$ is globally Lipschitz on $\R$ with constant
$L_b=\max\{6,139264/(5b^2)\}$.
\end{lemma}

The explicit piecewise-quartic formula for $\varphi_b$ and the elementary
calculus verification of Lemma~\ref{art:upper_scalar_lemma} are given in
Appendix~\ref{art:lower_curvature} (equation~\eqref{art:upper_scalar}), side by
side with the companion $C^{1,1}$ spline used for strong convexity. With
$\varphi_b$ defined, we obtain the first local failure theorem.

\begin{theorem}[Local failure with globally smooth convex branches]
\label{art:local-failure}
There exist four globally convex functions $f_i\in C^{2,1}(\R^2)$ satisfying
\[
 0\preceq\nabla^2f_i(x)\preceq64I\qquad(i=0,1,2,3;\;x\in\R^2)
\]
such that $F=\max_{0\le i\le3}f_i$ has compact sublevel sets and a unique
nonsmooth minimizer at $0$, together with starting points $x_0^{(\ell)}\to0$
($\ell\ge1$), all initialized with $H_0=I$, having the following property: for
each fixed $\ell\ge1$, the infinite full-BFGS run from $x_0^{(\ell)}$ generated
by the unit-first bracketing algorithm of Section~\ref{art:failure} is
independent of the choice of $(c_1,c_2)$ in the range
\[
 0<c_1\le\tfrac13,\qquad c_1<c_2<1,
\]
and satisfies
\[
 x_k^{(\ell)}\longrightarrow\bar x_\ell\ne0,\qquad
 F(x_k^{(\ell)})\downarrow F(\bar x_\ell)>F(0)=0,\qquad
 0\notin\partial F(\bar x_\ell).
\]
Every accepted and rejected trial point is differentiable, each line search
terminates in at most five trials with accepted step lengths repeating
$1,1/4,1/16$, and the supremum of the norms of all accepted and rejected query
points in the $\ell$th run tends to zero as $\ell\to\infty$.
\end{theorem}

\begin{proof}[Proof of Theorem~\ref{art:local-failure}]
Define the four-branch objective on $\R^2$ by
\begin{equation}\label{art:upper_objective}
 F(u,v)=\max\{\varphi_2(u),\;\varphi_2(-u-v/8),\;
                    \varphi_{b_*}(c_\circ u+d_\circ v),\;v\}.
\end{equation}
With $(b_0,b_1,b_2)=(2,2,b_*)$, for $i=0,1,2$ we have
$\nabla^2f_i(x)=\varphi_{b_i}''(g_i^Tx)g_ig_i^T\succeq0$; by
\eqref{art:scalar_curvature}, their operator norms are bounded by $51$,
$51(65/64)<52$, and $\max\{6,102/b_*\}\norm{g_2}^2=6\norm{g_2}^2<12$,
respectively, while $f_3(u,v)=v$
has zero Hessian, so $0\preceq\nabla^2f_i\preceq64I$ and each $\nabla^2f_i$ is
globally Lipschitz. Because $u$ and $-u-v/8$ are linearly independent and
$\varphi_2$ is coercive and nonnegative with unique zero at $0$, every sublevel
of $F$ is compact and $F(x)=0$ if and only if $x=0$. At $0$, $\nabla f_i(0)=0$
for $i=0,1,2$ while $\nabla f_3(0)=(0,1)^T$, so $F'(0;h)=\max\{0,h_2\}$ and $0$
is a nonsmooth minimizer.

For each $\ell\ge1$, set $\bar x_\ell=a_\ell(2,-32)^T$, where $g_0^T\bar x_\ell=g_1^T\bar x_\ell=2a_\ell$
and $g_2^T\bar x_\ell=b_*a_\ell$. Whenever $|g_i^Th|\le a_\ell/2$ for $i=0,1,2$,
the plateau identity~\eqref{art:scalar_plateau} puts the first three branches on
their affine segments with values at least $a_\ell^3/2>0$, while $|h_1|\le a_\ell/2$
and $|h_2|\le8(|h_1|+a_\ell/2)\le8a_\ell$ force $f_3(\bar x_\ell+h)=-32a_\ell+h_2<0$.
Hence
\begin{equation}\label{art:upper_local_model}
 F(\bar x_\ell+h)=a_\ell^3+a_\ell^2G(h)
 \qquad\text{whenever }|g_i^Th|\le a_\ell/2\quad(i=0,1,2),
\end{equation}
and in particular $F(\bar x_\ell)=a_\ell^3>0$ with $F'(\bar x_\ell;(-1,16)^T)=-a_\ell^2<0$,
so $0\notin\partial F(\bar x_\ell)$.

Start at $x_0^{(\ell)}=\bar x_\ell+a_\ell^2z_0$ with $H_0=I$. We claim by
induction that for every $m\ge0$ and $j\in\{0,1,2\}$,
\begin{equation}\label{art:upper_orbit}
 x_{3m+j}^{(\ell)}=\bar x_\ell+a_\ell^2\rho^mz_j,\qquad
 H_{3m+j}^{(\ell)}=\rho^mM_j.
\end{equation}
Indeed, at~\eqref{art:upper_orbit} every candidate query on the trial list of
Table~\ref{art:trial_table} has displacement $h=a_\ell^2\rho^mz_{\rm trial}$
from $\bar x_\ell$, which by~\eqref{art:seed_query_bounds} and $a_\ell\le1/16$
satisfies $|g_i^Th|<3a_\ell^2\le(3/16)a_\ell<a_\ell/2$. Thus~\eqref{art:upper_local_model}
holds at the current iterate and at every queried trial point with
$\lambda=a_\ell^2$, $\sigma=a_\ell^4=\lambda^2$, and $C=a_\ell^3$. By the
scale-covariance relations~\eqref{art:matched-scaling} and Lemma~\ref{art:seed_cycle},
the line search visits the exact nine-query tree of Table~\ref{art:trial_table},
makes the identical strict rejections and acceptances for every
$0<c_1\le1/3$, $c_1<c_2<1$, and closes at~\eqref{art:upper_orbit}. Consequently
$x_k^{(\ell)}\to\bar x_\ell$, $F(x_k^{(\ell)})\downarrow a_\ell^3>0$, and every
query norm is bounded by
$a_\ell\norm{(2,-32)^T}+a_\ell^2\max_{\text{seed queries}}\norm{z}\to0$ as
$\ell\to\infty$.\qed
\end{proof}

\subsection{Hyperbolic return-map persistence under strong convexity}

In Theorem~\ref{art:local-failure}, the flat affine plateaus where $\varphi_b''=0$
prevent the branches from being strongly convex. Adding a positive quadratic
$\tfrac\mu2\norm{z}^2$ restores uniform strong convexity, but perturbs the
branch gradients by $\mu z$ and breaks exact scale invariance: the homogeneous
seed orbit is no longer an exact trajectory. Although generic perturbations move
off the cycle, the homogeneous seed fixed point of a scale-normalized three-step
return map admits a two-dimensional local stable manifold extending into
positive scale, on which the nine-query decision tree persists and yields a
second local counterexample under uniform strong convexity.

\begin{theorem}[Local failure under uniform strong convexity]
\label{art:local-failure-lower}
There exist four uniformly strongly convex functions $f_i\in C^{1,1}(\R^2)$ with
globally Lipschitz gradients and constants $0<\underline\mu\le\overline L<\infty$
such that $F=\max_{0\le i\le3}f_i$ is finite and strongly convex on $\R^2$ with
compact sublevel sets, a unique nonsmooth minimizer at $0$, and
\[
 \underline\mu I\preceq\nabla^2F(x)\preceq\overline L I
 \quad\text{at every twice-differentiable point }x\in\R^2
\]
(with the first three branches not twice differentiable at $0$). Moreover, there
exist starting points $x_0^{(\ell)}\to0$ ($\ell\ge1$), all initialized with
$H_0=I$, having the following property: for each fixed $\ell\ge1$, the infinite
full-BFGS run from $x_0^{(\ell)}$ generated by the unit-first bracketing
algorithm of Section~\ref{art:failure} is independent of the choice of
$(c_1,c_2)$ for $0<c_1\le1/3$, $c_1<c_2<1$ and satisfies
\[
 x_k^{(\ell)}\longrightarrow\bar x_\ell\ne0,\qquad
 \bar x_\ell\longrightarrow0,\qquad
 F(x_k^{(\ell)})\downarrow F(\bar x_\ell)>F(0),\qquad
 0\notin\partial F(\bar x_\ell).
\]
Every accepted and rejected trial point is differentiable, the accepted step
lengths repeat $1,1/4,1/16$, and the supremum of all query norms in the
$\ell$th run tends to zero as $\ell\to\infty$.
\end{theorem}

The complete proof in Appendix~\ref{art:lower_curvature} introduces two further
analytical ideas summarized in Table~\ref{art:regularity-map}:
\begin{enumerate}
\item \emph{Scale-normalized return map and stable manifold.} On the strongly
convex prototype $Q_\mu(z)=\tfrac\mu2\norm{z}^2+G(z)$ (whose unique minimizer
 lies away from $0$), we factor out the shrinking scale $\tau>0$ via $x=\tau z$
and $H=\tau M(\zeta,\eta)$ with $M_{11}=1$. Three BFGS steps define a real-analytic
five-dimensional return map $\mathcal T_\mu(z,\zeta,\eta,\tau)$ for which the
homogeneous seed $P_*=(z_0,0,1,0)$ is a fixed point on the boundary $\tau=0$.
Rational interval differentiation shows that $D\mathcal T_\mu(P_*)$ is
hyperbolic, with two unstable eigenvalues equal to $\rho^{-1}$, a saddle pair
$\lambda_s\in(0.038,0.039)$ and $\lambda_u\in(1.564,1.565)$, and scale
eigenvalue $\rho$. Consequently, $P_*$ admits a two-dimensional local stable
manifold that penetrates into the region $\tau>0$. A sufficiently small
positive-$\tau$ point on that stable manifold generates a genuine
perturbed BFGS run converging to the nonstationary origin with identical trial
lists across all $0<c_1\le1/3$, $c_1<c_2<1$ (Proposition~\ref{art:quadratic_failure}).
\item \emph{Orbit-copying spline and global whitening.} Because $Q_\mu$ itself
is not scale-homogeneous, we select a single sufficiently small perturbed
orbit $(x_k^{\rm seed},H_k^{\rm seed})$ and replicate it at centers
$\bar x_\ell=a_\ell(2,-32)^T$ ($a_\ell=128^{-\ell}$) via scalar splines $\psi_b$
with node values $a_\ell^2$ and slopes $a_\ell$. This yields the exact local
model $F(\bar x_\ell+a_\ell y)=a_\ell^2(C_*+Q_\mu(y))$, so $\lambda=a_\ell$ and
$\sigma=a_\ell^2$ once again satisfy $\lambda^2/\sigma=1$, giving the same
initial matrix $H_0^{\rm seed}\succ0$ on every copy. A single global linear
change of coordinates $z\mapsto(H_0^{\rm seed})^{1/2}z$ transforms that common
initial matrix into $H_0=I$ simultaneously for all $\ell\ge1$. Because the ratio
of node slope $a_\ell$ to node position $ba_\ell$ stays $\Theta(1)$ as
$\ell\to\infty$, $\psi_b$ remains globally $C^{1,1}$ while losing second differentiability at
$0$, which is where this second local construction loses $C^2$ regularity
(compare Theorems~\ref{art:local-failure} and~\ref{art:local-failure-lower}).
\end{enumerate}

\begin{table}[tbp]
\centering
\begin{tabular}{@{}>{\raggedright\arraybackslash}p{.29\textwidth}>{\raggedright\arraybackslash}p{.31\textwidth}>{\raggedright\arraybackslash}p{.32\textwidth}@{}}
\toprule
Branch class & Global curvature bounds & Regularity at $x_*=0$\\
\midrule
Four convex $C^{2,1}$ branches (Theorem~\ref{art:local-failure})
& $0\preceq\nabla^2f_i(x)\preceq64I$ on $\R^2$
& Smooth $C^{2,1}$ branches everywhere; curvature vanishes on plateaus\\[3pt]
Four strongly convex $C^{1,1}$ branches (Theorem~\ref{art:local-failure-lower})
& Each $f_i$ uniformly strongly convex with globally Lipschitz $\nabla f_i$; $\underline\mu I\preceq\nabla^2F(x)\preceq\overline L I$ where twice differentiable
& Uniformly strongly convex $C^{1,1}$ branches; $f_0,f_1,f_2$ not $C^2$ at $0$\\
\bottomrule
\end{tabular}
\caption{Comparison of the two local smooth-branch counterexamples. Here
$\nabla^2F(x)$ denotes the classical Fr\'echet derivative of $\nabla F$ defined
on a neighborhood of $x$. Both theorems hold for one fixed objective, $H_0=I$,
and all $0<c_1\le1/3$, $c_1<c_2<1$ simultaneously, with all accepted and
rejected queries shrinking to $x_*=0$.}
\label{art:regularity-map}
\end{table}

\bookmarksetup{startatroot}
\addtocontents{toc}{\protect\vspace{6pt}}
\section{Sharp Geometric Boundaries, Open Problems, and Formal Verification}\label{art:verification}

Bringing Part~\ref{part:convergence} and Part~\ref{part:failure} together shows
how closely the positive and negative regimes meet. The convergence conditions
of Sections~\ref{art:global}--\ref{art:local} are sufficient rather than
necessary: Chapter~5 of the companion gives a one-dimensional convex example
with accumulating kinks on which every admissible full-BFGS run converges
without admitting a uniform sublevel-band certificate. Complementing the
positive theorems, the counterexamples of
Sections~\ref{art:failure}--\ref{art:failure-local} explain why neither
convexity, coercivity, nor smooth or strongly convex branch formulas can replace
those geometric conditions: when $0\in\aff\partial F\setminus\conv\partial F$
eliminates equal-rate descent directions, unmodified full BFGS from $H_0=I$ can
enter a self-similar cycle and converge to a nonstationary point.

\subsection{Why the counterexamples escape the positive theory}

Each counterexample violates a specific hypothesis of the positive theorems:
\begin{enumerate}
\item \emph{Global homogeneity vs.\ affine offsets.} The three-slope cone
$G(z)=\max_{0\le i\le2}g_i^Tz$ is positively homogeneous of degree one, but
satisfies $G(t(-1,16)^T)=-t$ for $t\ge0$ and so fails coercivity; its four-piece
coercive cap $\overline G$ in Theorem~\ref{art:failure_polyhedral} includes the
affine offset $v-17$, breaking homogeneity around the minimizer $(-1,16)$ and
leaving $(0,0)$ as a nonstationary corner with $0\in\aff\partial\overline G(0)\setminus\conv\partial\overline G(0)$.
\item \emph{Finite local polyhedral partitions vs.\ accumulating smooth kinks.}
Near $(-1,16)$, $\overline G$ is locally a three-piece polyhedral cone centered
at $(-1,16)$ and therefore converges locally by Theorem~\ref{art:polyhedral}.
By contrast, in Theorems~\ref{art:local-failure} and~\ref{art:local-failure-lower},
nonlinear scalar splines cause nonstationary three-slope corners $\bar x_\ell\to0$
to accumulate at the minimizer $0$: no neighborhood of $0$ admits a finite
polyhedral partition with smooth extensions (since Theorem~\ref{art:polyhedral}
would otherwise force convergence once $F(x_0^{(\ell)})\le\varepsilon$),
demonstrating that finitely many smooth formulas in $\max_if_i$ do not imply a
finite local polyhedral cell structure.
\item \emph{Active-gradient multiplier nondegeneracy.} Both local smooth-branch
counterexamples are planar four-branch maxima, so we can test the multiplier
hypotheses of Corollary~\ref{art:planar-max} directly at $x_*=0$. In
\eqref{art:upper_objective}, all four branches are active at $0$ with gradients
$0,0,0,e_2$, where $e_2=(0,1)^T$; in~\eqref{art:lower_objective}, the active
gradients at $0$ are again $0,0,0,e_2$ before whitening and $0,0,0,A^Te_2\ne0$
after $x=Az$ with $A=(H_0^{\rm seed})^{1/2}$. In both cases, the repeated zero
gradients violate affine independence, and any convex combination placing
positive weight $\lambda_3>0$ on the fourth branch equals $\lambda_3e_2\ne0$ (or
$\lambda_3A^Te_2\ne0$), violating strict complementarity $\lambda_i^*>0$ as
well; in~\eqref{art:lower_objective}, three branches also lose $C^2$ regularity
at $0$.
\end{enumerate}

A one-line secant identity also clarifies why uniform strong convexity cannot by
itself prevent BFGS from converging to a nonstationary point. On $Q_\mu(x)=\tfrac\mu2\norm{x}^2+G(x)$,
if consecutive differentiable iterates $x$ and $x+s$ activate slopes $g_i$ and
$g_j$, then the secant vector is
\begin{equation}\label{art:jump-secant}
 y=\bigl(\mu(x+s)+g_j\bigr)-\bigl(\mu x+g_i\bigr)
  =\mu s+(g_j-g_i).
\end{equation}
Along the perturbed trajectory of Proposition~\ref{art:quadratic_failure},
$\norm{s_k}\to0$ while $g_0,g_1,g_2$ continue to cycle: the strongly convex
contribution $\mu s_k$ vanishes with the step length, leaving the $O(1)$ jump
$g_j-g_i$ to dominate the rank-two BFGS update.

\subsection{Open problems}

Our results point directly to four open problems in the theory of quasi-Newton
methods for nonsmooth optimization:
\begin{enumerate}
\item \emph{Uniformly strongly convex $C^2$ branch maxima.} Globally,
$Q_\mu(x)=\tfrac\mu2\norm{x}^2+G(x)$ is already a maximum of three $\mu$-strongly
convex $C^\infty$ quadratic functions on which full BFGS can converge to the
nonstationary origin (Proposition~\ref{art:quadratic_failure}). Locally near the
minimizer, however, Table~\ref{art:regularity-map} leaves a sharp regularity
gap: Theorem~\ref{art:local-failure} achieves $C^{2,1}$ branches without a
positive lower Hessian bound, Theorem~\ref{art:local-failure-lower} achieves
uniformly strongly convex $C^{1,1}$ branches that lose $C^2$ regularity at $0$,
and Corollary~\ref{art:planar-max} proves local convergence in $\R^2$ when the
active gradients at $x_*$ are nondegenerate. Does unmodified full BFGS converge
locally from every differentiable start near the minimizer of a maximum of
finitely many uniformly strongly convex $C^2$ functions when active-gradient
nondegeneracy is dropped?
\item \emph{Moving comparison fields in dimension $n\ge3$.} The radial theorems
(Theorems~\ref{art:centered}, \ref{art:homogeneous}, and~\ref{art:polyhedral})
hold in arbitrary dimension $\R^n$, whereas our curved-interface analysis
(Theorem~\ref{art:planar-cluster} and Corollary~\ref{art:planar-max}) uses the
one-dimensionality of level-set tangents in $\R^2$ to cancel normal components
of $\partial_\theta(hn+h_\theta T)$ and charge field rotation to functional
decrease. Extending moving comparison fields, or the higher-rank comparison tensors
$W_\eta(x)$ already permitted by Theorem~\ref{art:certificate}, to nondegenerate
smooth maxima with curved switching manifolds in $\R^n$ ($n\ge3$) is a natural
next step.
\item \emph{Convergence rates via principal minors and comparison energies.}
Theorem~\ref{art:reciprocal-growth} and our comparison-energy inequalities give
asymptotic convergence without conditioning assumptions (finite-horizon
estimates are recorded in the companion). Determining sharp non-asymptotic or
R-linear convergence rates for structured nonsmooth families such as
$F(u,v)=u^2+|v|$ or coercive polyhedral norms remains largely open.
\item \emph{Positive-measure vs.\ zero-measure failure on convex functions
(Lewis--Overton Challenge~7.1).} On \emph{nonconvex} coercive semialgebraic
objectives, Lai~\cite{NEG_LAI26} recently constructed a nonempty open
set of initial states $(x_0,H_0)$ converging to a non-Clarke-stationary point by
using negative transverse curvature to contract perturbations onto a homogeneous
cycle. Under convexity, however, Theorem~\ref{art:homogeneous} rules out global
failure on every coercive positively homogeneous objective, and the failing
starts exhibited in Theorems~\ref{art:failure_polyhedral}, \ref{art:local-failure},
and~\ref{art:local-failure-lower} consist of a single point $z_0$ or discrete
sequences $x_0^{(\ell)}$ (or lie on a two-dimensional stable manifold in the
normalized state space of Appendix~\ref{art:lower_curvature}). Thus our
deterministic counterexamples do not contradict empirical success from typical
initializations: under absolutely continuous random initialization as in
Challenge~7.1~\cite{LO}, can a \emph{convex} nonsmooth function admit a
positive-measure set of initial states $(x_0,H_0)$ on which full BFGS converges
to nonstationary points, or are all convex nonstationary basins of measure zero?
\end{enumerate}

\subsection{Formal verification and reproducibility}

All results and rational interval certificates in this article, together with
the extended noncompact and quantitative theory in the accompanying companion
monograph (\emph{BFGS on Convex Nonsmooth Functions: Complete Companion},
\texttt{companion/}), have been formalized and verified in Lean~4
(\texttt{v4.19.0}) on top of pinned \textsf{Mathlib}
(\texttt{c44e0c8ee63c}). Every verified declaration has been audited to depend
only on the three standard foundational axioms of Lean's classical logic
(\texttt{propext}, \texttt{Classical.choice}, and \texttt{Quot.sound};
\path{bfgs_lean_article/AXIOM_AUDIT.json}), with every compiled project
module additionally re-verified through an independent \texttt{lean4checker}
kernel replay (\path{bfgs_lean_article/KERNEL_REPLAY.json}). A complete
statement-by-statement map from the 38 numbered statements of this article to their
formal Lean~4 declarations (\path{ARTICLE_LEAN_MAP.md}), the item-level
coverage registers (\path{editorial/lean_coverage/register/}), and
reproducible build and audit scripts are publicly archived on Zenodo at
\url{https://doi.org/10.5281/zenodo.24648859}.

\paragraph{Acknowledgements.}
Large language models assisted theorem development, proof reconstruction, Lean~4
formalization, literature searches, software tooling, and manuscript preparation;
the companion report documents this workflow and the formal verification audits
in detail.

\bookmarksetup{startatroot}
\appendix
\addtocontents{toc}{\protect\vspace{6pt}\noindent{\protect\footnotesize\protect\sffamily\protect\bfseries\protect\color{CiteTeal} APPENDICES FOR PART I}\protect\par\protect\vspace{2pt}}
\section{Gram compression for the inverse-BFGS update}\label{art:matrix-algebra}

This appendix establishes the subspace determinant inequality underlying the
principal-minor hierarchy (Theorem~\ref{art:minor-hierarchy}).

\begin{lemma}[Subspace Gram compression]\label{art:gram-compression}
Under the assumptions of Theorem~\ref{art:minor-hierarchy}, let $Z\in\R^{n\times j}$
have orthonormal columns ($1\le j\le n$). Then
\[
 \det(Z^TH_+Z)\le
 \frac{\det\operatorname{Gram}_H(y,Z)}{\beta^2q}
       +\frac\alpha\beta\det(Z^THZ),
\]
and at $j=n$ the exact identity $\det H_+=(\alpha/\beta)\det H$ holds.
\end{lemma}
\begin{proof}
On the chosen $j$-dimensional subspace, the projected inverse-BFGS term reduces
to a rank-two perturbation, so only a two-dimensional determinant has to be
expanded. Set $v=Hg$, $\kappa=g^THy=-\beta q$, and $P=I-vy^T/\kappa$. Because
$y^Ts=\alpha\beta q$, the inverse-BFGS update~\eqref{art:update} decomposes as
\[
 H_+=\mathcal F+\frac{\alpha}{\beta q}vv^T,
 \qquad \mathcal F=PHP^T.
\]
Writing
\[
 \begin{gathered}
 A=Z^THZ\succ0,\quad a=Z^THg,\quad b=Z^THy,\quad r=y^THy,\\
 u=a^TA^{-1}a,\quad w=a^TA^{-1}b,\quad z=b^TA^{-1}b,
 \end{gathered}
\]
the projected term compresses to a rank-two perturbation of $A$:
\[
 Z^T\mathcal FZ
 =A-[a\ b]
 \begin{pmatrix}-r/\kappa^2&1/\kappa\\1/\kappa&0\end{pmatrix}
 [a\ b]^T.
\]
Applying Sylvester's identity $\det(I+UV)=\det(I+VU)$ reduces the determinant
ratio $\det(Z^T\mathcal FZ)/\det A$ to a $2\times2$ determinant whose expansion
gives
\begin{equation}\label{art:compression-determinant}
 \det(Z^T\mathcal FZ)
 =\frac{\det A}{\kappa^2}
       \bigl((\kappa-w)^2+u(r-z)\bigr).
\end{equation}
For $\operatorname{Gram}_H(y,Z)=[y\ Z]^TH[y\ Z]$ and
$\operatorname{Gram}_H(g,y,Z)=[g\ y\ Z]^TH[g\ y\ Z]$, taking Schur complements
with respect to $A\succ0$ yields
\begin{align*}
 \det\operatorname{Gram}_H(y,Z)&=(r-z)\det A\ge0,\\
 \det\operatorname{Gram}_H(g,y,Z)
 &=\bigl((q-u)(r-z)-(\kappa-w)^2\bigr)\det A\ge0,
\end{align*}
both nonnegative as Gram determinants (even if their columns are linearly
dependent). Substituting these two identities into~\eqref{art:compression-determinant}
and using $\kappa^2=\beta^2q^2$ gives
\[
 \det(Z^T\mathcal FZ)
 =\frac{q\det\operatorname{Gram}_H(y,Z)
             -\det\operatorname{Gram}_H(g,y,Z)}{\beta^2q^2}
 \le\frac{\det\operatorname{Gram}_H(y,Z)}{\beta^2q}.
\]

Next, adding $\rho aa^T$ to $Z^T\mathcal FZ$ replaces the top-left entry
$-r/\kappa^2$ of the $2\times2$ coupling matrix by $-r/\kappa^2-\rho$, which by
the same $2\times2$ determinant expansion increases the determinant by
$\rho u\det A$ without requiring $Z^T\mathcal FZ$ to be invertible. Since
$Z^Tv=a$, setting $\rho=\alpha/(\beta q)$ recovers $Z^TH_+Z$; moreover, the
Schur complement of $A$ in $\operatorname{Gram}_H(g,Z)\succeq0$ gives $0\le u\le q$,
so $\rho u\det A\le(\alpha/\beta)\det A$, proving the inequality.

Finally, when $j=n$, $Z$ is orthogonal, so $A^{-1}=Z^TH^{-1}Z$ gives $u=q$,
$w=\kappa$, and $z=r$: the right-hand side of~\eqref{art:compression-determinant}
vanishes, the rank-one increment equals $(\alpha/\beta)\det A=(\alpha/\beta)\det H$,
and orthogonal invariance gives $\det H_+=(\alpha/\beta)\det H$.\qed
\end{proof}

\section{Constructing the planar contour comparison field}\label{art:planar-appendix}

This appendix constructs the planar comparison vector field $w(x)$ for
Theorem~\ref{art:planar-cluster} by prescribing the level-curve velocity at each
switching corner, extending it smoothly along the intervening smooth arcs, and
verifying the required visibility, transversality, transport, and secant bounds.

\begin{lemma}[Planar contour comparison field]\label{art:planar-field}
Let $F:\R^2\to\R$ be finite and convex with minimum $F_*$, and let
$C\subset\{x:F(x)=\ell\}$ be a nonempty compact set at level $\ell>F_*$.
Suppose the boundary of $\{F\le\ell\}$ contains no nontrivial line segment and
every $p\in C$ lies outside $\mathcal C_{\rm aff}$ and satisfies the local $C^2$
or two-branch regularity hypothesis of Theorem~\ref{art:planar-cluster}. Fix
$0<c_1<c_2<1$. Then there exist an open neighborhood $U\supset C$, a compact
buffer $Q\supset\overline U$, a vector field $w$ defined on a neighborhood of
$Q$, and constants $\eta,\kappa,\delta>0$ and $C_0,C_1,L\ge0$ such that:
\begin{enumerate}
\item $w$ is locally Lipschitz, bounded on $Q$, takes a common inner product
with all subgradients at each point of $Q$, and satisfies $g^Tw(p)\ge\eta$ for
every $p\in Q$ and $g\in\partial F(p)$;
\item for every eligible pair $(x,z)$ in $U$ with $\norm{s}<\delta$ ($s=z-x$),
the segment $[x,z]$ lies in $Q$ and
\[
 |\det(s,w(x+\theta s))|\ge\kappa\norm{s}
 \quad(0\le\theta\le1),\qquad \frac{(y^Tw(z))^2}{t}\le C_0;
\]
\item for every such pair there exist $\lambda>0$ and $\beta\in\R$ satisfying
\[
 w(z)=\lambda w(x)+\beta s,\qquad
 |\log\lambda|\le C_1d,
 \qquad d\le1\ \Longrightarrow\ \lambda^2\le1+Ld.
\]
\end{enumerate}
\end{lemma}
\begin{proof}
\emph{1. Corner trajectories and angular separation.}
At any two-branch point of $C$, linear independence of $\nabla f_1,\nabla f_2$
makes $x\mapsto(f_1(x),f_2(x))$ a local $C^2$ diffeomorphism (injective on a
shrunk chart), so the inverse
function theorem parametrizes the switching curve near level $\ell$ by a $C^2$
path $p_j(t)$ solving $f_1(p_j(t))=f_2(p_j(t))=t$. Differentiating with respect
to $t$ shows that its $C^1$ velocity $v_j(t)=p_j'(t)$ satisfies
\begin{equation}\label{art:corner-velocity}
 \nabla f_1(p_j(t))^Tv_j(t)=\nabla f_2(p_j(t))^Tv_j(t)=1.
\end{equation}
Covering the compact set $C$ by finitely many buffered smooth and two-branch
charts shows that $C$ contains only finitely many corners, and on a shrunk
neighborhood of $C$ and level interval around $\ell$, every nonsmooth point lies
on one of these curves $p_j(t)$.

At each level-$\ell$ corner $p_j(\ell)$, the positive cone spanned by its two
active gradients cuts out a compact arc of unit supporting normals of length
strictly less than $\pi$, on which $n^Tv_j(\ell)>0$ by~\eqref{art:corner-velocity}.
For distinct corners on $\{F=\ell\}$, these normal arcs are disjoint: a shared
supporting normal would expose both corners and hence the entire line segment
between them on the boundary of $\{F\le\ell\}$, contradicting strict convexity
of the contour.

\emph{2. Normal-angle blending and tangent preservation.}
Because the level-$\ell$ normal arcs have positive angular separation, we can
choose smooth $2\pi$-periodic cutoffs $0\le\chi_j(\theta)\le1$ with pairwise
disjoint supports such that $\chi_j\equiv1$ on an open enlargement of the $j$th
arc and $n(\theta)^Tv_j(t)>0$ throughout $\operatorname{supp}\chi_j$ on a
shrunk level interval around $\ell$. Writing $n(\theta)=(\cos\theta,\sin\theta)^T$
and $T(\theta)=(-\sin\theta,\cos\theta)^T$, set
\begin{equation}\label{art:planar-height}
 h(t,\theta)=1+\sum_j\chi_j(\theta)\bigl(n(\theta)^Tv_j(t)-1\bigr).
\end{equation}
Since at most one $\chi_j(\theta)$ is nonzero at any $\theta$, $h(t,\theta)$ is
a convex combination of $1$ and $n(\theta)^Tv_j(t)$, hence smooth in $\theta$,
$C^1$ in $t$, and uniformly bounded above and away from zero. At every smooth
point $x$, let $n(\theta(x))=\nabla F(x)/\norm{\nabla F(x)}$ and define
\begin{equation}\label{art:planar-field-formula}
 w(x)=h(F(x),\theta(x))n(\theta(x))
       +h_\theta(F(x),\theta(x))T(\theta(x)),
\end{equation}
which depends periodically on $\theta(x)$. On an open neighborhood $N_j$ of the
$j$th switching curve inside its chart (where levels lie in the shrunk interval
and unit normals satisfy $\chi_j\equiv1$ with all other cutoffs zero), $h=n^Tv_j$,
$h_\theta=T^Tv_j$, and~\eqref{art:planar-field-formula} collapses on both smooth
sides to
\begin{equation}\label{art:planar-corner-field}
 w(x)=v_j(F(x)).
\end{equation}
Defining $w(x)=v_j(F(x))$ on the switching curve itself makes $w$ locally
Lipschitz across every kink and $C^1$ elsewhere. On a compact buffer $Q$ around
$C$ inside that level preimage where $F-F_*>0$, the supporting inequality
$\norm{g}\norm{x-x_*}\ge F(x)-F_*$
bounds $\norm{\nabla F}$ uniformly away from zero at smooth points, where
$\nabla F(x)^Tw(x)=\norm{\nabla F(x)}h(F(x),\theta(x))$, while at switching
points every $g\in\partial F(p)$ satisfies $g^Tw(p)=1$ by~\eqref{art:corner-velocity}.
Thus $w$ is locally Lipschitz and bounded on $Q$, pairs equally with all
subgradients at each point, and satisfies $g^Tw(p)\ge\eta>0$ throughout $Q$.

Differentiating $hn+h_\theta T$ with respect to $\theta$ cancels the normal
terms (Bogosel~\cite[Section~2]{geom_Bogosel2023}):
\begin{equation}\label{art:planar-tangent-cancellation}
 \partial_\theta(hn+h_\theta T)=(h+h_{\theta\theta})T.
\end{equation}
At any smooth point, since $\nabla F$ is orthogonal to $T$, directional
differentiation along $T$ differentiates only $\theta(x)$ and $F(x)$ has zero
derivative, so $Dw$ maps $T$ into $\R T$; in the orthonormal basis $(T,n)$,
\[
 Dw=\begin{pmatrix}a&b\\0&e\end{pmatrix}
\]
with uniformly bounded entries $a,b,e$ on $Q$.

\emph{3. Uniform transversality of eligible steps.}
Fix an open neighborhood $U\supset C$ with $\overline U\subset\operatorname{int}Q$.
Suppose the determinant lower bound in item~2 failed for arbitrarily short
eligible pairs in $U$: then there would exist eligible pairs $(x_m,z_m)$ in $U$
coalescing to some $p\in Q$, unit steps $s_m/\norm{s_m}\to\hat s$, and
$\theta_m\in[0,1]$ with $\det(\hat s,w(p))=0$. Upper semicontinuity of
$\partial F$ gives subsequential limits $\bar g,\bar h\in\partial F(p)$ of the
endpoint gradients. By item~1, $\bar g^Tw(p)=\bar h^Tw(p)\ge\eta>0$; since
$\hat s$ is parallel to $w(p)$ and $d_m/\norm{s_m}\to-\bar g^T\hat s\ge0$, we
must have $\hat s=-w(p)/\norm{w(p)}$. But then
\[
 \frac{d_m}{\norm{s_m}}\longrightarrow-\bar g^T\hat s>0
 \qquad\text{whereas}\qquad
 \frac{t_m}{\norm{s_m}}\longrightarrow(\bar h-\bar g)^T\hat s=0,
\]
contradicting weak Wolfe $t_m\ge(1-c_2)d_m$. Thus $|\det(s,w(x+\theta s))|\ge\kappa\norm{s}$
holds uniformly once $\delta>0$ is small.

\emph{4. Transport and secant bounds.}
Along an eligible segment $q(\tau)=x+\tau s$ ($\tau\in[0,1]$), set
$\varphi(\tau)=F(q(\tau))$ and $D(\tau)=\det(s,w(q(\tau)))$, which is nonvanishing
of constant sign on $[0,1]$. At any smooth point $q(\tau)$, $(T,n)$ has
determinant $-1$, so expanding $s=s_TT+s_nn$ gives
\begin{equation}\label{art:planar-determinant-derivative}
 \det(s,Dw\,s)=s_n\bigl[(a-e)s_T+bs_n\bigr].
\end{equation}
Because $|D(\tau)|\ge\kappa\norm{s}$ and $\varphi'(\tau)=\norm{\nabla F(q(\tau))}s_n$
with $\norm{\nabla F}$ bounded away from zero on $Q$, we have
$|D'(\tau)/D(\tau)|\le A|\varphi'(\tau)|$; at a switching point,
\eqref{art:planar-corner-field} gives $w\circ q=v_j\circ\varphi$, so the
one-sided chain rule for the convex function $\varphi$ gives the identical bound
for right derivatives on $[0,1)$. Since $\varphi'_+(\tau)\ge\varphi'_+(0)=-d$,
we have $|\varphi'_+(\tau)|\le\varphi'_+(\tau)+2d$, so $H(\tau)=\log|D(\tau)|$
satisfies $|H'_+(\tau)|\le A(\varphi'_+(\tau)+2d)$ for all $\tau\in[0,1)$.
Because any continuous function $u$ on $[0,1]$ with $u'_+\le0$ on $[0,1)$
satisfies $u(1)\le u(0)$ (otherwise $u(\tau)-\epsilon\tau$ would attain an
interior minimum before $1$ with nonnegative right derivative), applying this to
$\pm(H(\tau)-H(0))-A(\varphi(\tau)-\varphi(0)+2d\tau)$ and using
$\varphi(1)\le\varphi(0)$ yields $|H(1)-H(0)|\le2Ad$. Writing
$w(z)=\lambda w(x)+\beta s$ in the independent basis $(w(x),s)$ gives
$\lambda=D(1)/D(0)>0$, so $|\log\lambda|\le2Ad$, and convexity of $\exp$ on
$[0,4A]$ gives $\lambda^2\le\exp(4Ad)\le1+(e^{4A}-1)d$ whenever $d\le1$,
proving item~3.

Finally, to bound $(y^Tw(z))^2/t$, refine $\delta$ so that every $\delta$-close
pair in $U$ lies inside a single buffered smooth or two-branch chart on which
branch gradients are $M$-Lipschitz ($M>0$) and $z-(\nabla f(z)-\nabla f(x))/M$
remains in the convex chart. If $x$ and $z$ activate the same branch $f$, setting
$b(q)=f(q)-\nabla f(x)^Tq$ and $y=\nabla f(z)-\nabla f(x)$ gives
$b(x)\le b(z-y/M)\le b(z)-\norm{y}^2/(2M)$ by convexity of $b$ at $x$ and the
$M$-Lipschitz descent bound along the buffered segment from $z$ to $z-y/M$
(cf.~\cite[Section~2]{WW}), yielding the local one-sided Bregman estimate
\[
 \norm{\nabla f(z)-\nabla f(x)}^2
 \le2M\bigl(f(z)-f(x)-\nabla f(x)^T(z-x)\bigr)\le2Mt,
\]
so boundedness of $w$ bounds $(y^Tw(z))^2/t$. If $x$ and $z$ activate distinct
branches $f_i$ and $f_j$ in a two-branch chart, the segment $[x,z]$ contains a
tie point $p$ where $(\nabla f_j(p)-\nabla f_i(p))^Tw(p)=0$; Lipschitz
continuity of $\nabla f_i,\nabla f_j,$ and $w$ gives $|y^Tw(z)|\le A_1\norm{s}$.
It remains to bound $t$ from below in terms of $\norm{s}$. Writing
$r=\nabla f_j(z)-\nabla f_j(x)$ ($\norm{r}\le M\norm{s}$), we
have $\nabla f_i(x)^Ts=-d$ and $\nabla f_j(x)^Ts=t-d-r^Ts$; uniform linear
independence of $\nabla f_i,\nabla f_j$ on the chart and $d\le t/(1-c_2)$ yield
\[
 \norm{s}\le A_2(2d+t+M\norm{s}^2)
 \le A_2\left(1+\frac{2}{1-c_2}\right)t+A_2M\norm{s}^2.
\]
Absorbing $A_2M\norm{s}^2\le\tfrac12\norm{s}$ for small $\delta$ gives
$t\ge a_0\norm{s}$ ($a_0>0$), whence $(y^Tw(z))^2/t\le(A_1^2/a_0)\norm{s}$ is
uniformly bounded.\qed
\end{proof}

\section{Worked certificates and existence of admissible runs}\label{art:hypothesis-examples}

This appendix exhibits explicit band certificates and full-BFGS trajectories for
Definition~\ref{art:band-certificate}, establishes general sufficient conditions
for the existence of infinite admissible runs from prescribed initial data, and
shows that even from a differentiable point on a coercive cone, a BFGS search
ray need not contain any differentiable weak-Wolfe step. Throughout, norms are
Euclidean.

\subsection{Explicit band certificates and BFGS trajectories}

\begin{proposition}[Degenerating trajectory and band certificate on $u^2+|v|$]\label{art:examples-mixed}
Let $F(u,v)=u^2+|v|$ and $K_b=\{x\in\R^2:F(x)\le b\}$ ($b>0$). Then $F$ is
finite and convex with unique minimizer $0$ and compact sublevel $K_b$, and for
every $0<\eta\le b$, Definition~\ref{art:band-certificate} holds on $K_b$ with
\[
 W_\eta(x)=xx^T,\qquad m_\eta=\eta^2,\quad C_\eta=2b,
 \quad L_\eta=0,\quad\delta_\eta=1.
\]
Moreover, for $b=2/3$ and $(c_1,c_2)=(1/4,1/2)$, the sequence
\[
 x_k=\bigl(0,\tfrac23(-\tfrac12)^k\bigr)^T,\qquad
 g_k=\bigl(0,(-1)^k\bigr)^T,\qquad
 H_k=\begin{pmatrix}1&0\\0&2^{-k}\end{pmatrix},\qquad \alpha_k=1
\]
is an infinite admissible run in $K_{2/3}$ starting from $H_0=I$ and converging
to $(0,0)$, along which $\lambda_{\min}(H_k)=2^{-k}\to0$.
\end{proposition}
\begin{proof}
At any differentiable point $x=(u,v)^T\in K_b$ with $F(x)\ge\eta$,
$g=(2u,\operatorname{sign}v)^T$ gives $g^TW_\eta(x)g=(2u^2+|v|)^2\ge F(x)^2\ge\eta^2$.
For any differentiable $x=(u,v)^T$ and $z=(r,w)^T$ in $K_b$ with $t=y^T(z-x)>0$,
the two-case Cauchy--Schwarz calculation in the proof of
Corollary~\ref{art:mixed-benchmark} gives $(y^Tz)^2\le2F(z)t\le2bt$, and
$z=x+s$ gives $\Pi_sz=\Pi_sx$, verifying~\eqref{art:visibility}--\eqref{art:transport}.
Along the displayed sequence, $x_k$ is differentiable with $g_k\ne0$,
$s_k=-H_kg_k=(0,-(-1/2)^k)^T$, $y_k=-2g_k$, $d_k=2^{-k}$,
$F(x_k)-F(x_{k+1})=\tfrac13d_k\ge\tfrac14d_k$, $t_k=2d_k$, and~\eqref{art:update}
maps $\operatorname{diag}(1,2^{-k})$ to $\operatorname{diag}(1,2^{-k-1})$.\qed
\end{proof}

\begin{proposition}[Full-rank band certificate in $\R^n$]\label{art:examples-full-rank}
For $F(x)=\tfrac12\norm{x}^2$ on $\R^n$ ($n\ge1$) and $b>0$, $K_b=\{F\le b\}$ is
compact with unique minimizer $0$, and for every $0<\eta\le b$,
Definition~\ref{art:band-certificate} holds with the rank-$n$ constant field
\[
 W_\eta(x)=I,\qquad m_\eta=2\eta,\quad C_\eta=1,
 \quad L_\eta=0,\quad\delta_\eta=1.
\]
For every $x_0\ne0$, $x_k=2^{-k}x_0$, $g_k=x_k$, $H_k=I$, $\alpha_k=1/2$ is an
infinite admissible run in $K_{F(x_0)}$ for $(c_1,c_2)=(1/4,1/2)$ converging to
$0$.
\end{proposition}
\begin{proof}
Here $g^TIg=\norm{x}^2=2F(x)\ge2\eta$, $y=s$ gives $y^TIy=t$, and constant
$W_\eta=I$ satisfies~\eqref{art:transport} with $L_\eta=0$. Along the stated
history, $y_k=s_k=-\tfrac12x_k$ leaves $H_k=I$ invariant under~\eqref{art:update},
$F(x_k)-F(x_{k+1})=\tfrac34d_k\ge\tfrac14d_k$, and $g_{k+1}^Ts_k=\tfrac12g_k^Ts_k$.\qed
\end{proof}

\subsection{A general existence criterion and curved-interface examples}

\begin{lemma}[Existence of differentiable Armijo--weak-Wolfe steps and runs]\label{art:examples-run-existence}
Let $F:\R^n\to\R$ be finite, convex, and bounded below, and suppose that for
every differentiability point $x$ of $F$ and every direction $p\ne0$, the set
\[
 \{t\in\R:F\text{ is not differentiable at }x+tp\}
\]
is countable (where differentiability refers to $F$ on $\R^n$). Then at any
differentiable $x$ with $g=\nabla F(x)\ne0$, any $H\succ0$, and any
$0<c_1<c_2<1$, the search ray $x+tp$ ($p=-Hg$, $t>0$) contains a differentiable
point $x+\alpha p$ satisfying both~\eqref{art:aw} with strict inequalities. If,
in addition, $\nabla F(x)\ne0$ at every differentiability point of $F$, then an
infinite admissible run exists from every differentiable $x_0$, every
$H_0\succ0$, and every $0<c_1<c_2<1$.
\end{lemma}
\begin{proof}
Set $\phi(t)=F(x+tp)$, $D=g^Tp<0$, and $\sigma(t)=(\phi(t)-\phi(0))/t$ for
$t>0$. Differentiability at $x$ gives $\lim_{t\downarrow0}\sigma(t)=D$, while
boundedness below gives $\liminf_{t\to\infty}\sigma(t)\ge0>c_1D$. Since
$D<c_2D<c_1D<0$, continuity of $\sigma$ yields a nonempty open interval of
$t>0$ on which $c_2D<\sigma(t)<c_1D$. Picking $\alpha$ in that interval outside
the countable nondifferentiable set gives strict Armijo $\phi(\alpha)<\phi(0)+c_1\alpha D$,
and convexity of $\phi$ gives $\nabla F(x+\alpha p)^Tp=\phi'(\alpha)\ge\sigma(\alpha)>c_2D$,
hence strict weak Wolfe. By Lemma~\ref{art:run-invariants}, $H_+\succ0$; if no
differentiable point is stationary, repeating this choice inductively yields an
infinite admissible run.\qed
\end{proof}

\begin{proposition}[Curved two-quadratic family with global convergence and run existence]\label{art:examples-two-quadratics}
For any $a,b,c>0$, define on $\R^2$
\[
 F_{a,b,c}(u,v)=\max\{au^2+bv^2+v,\;cu^2+bv^2-v\},
\]
whose switching locus is the parabola $v=\tfrac12(c-a)u^2$ (a line only when
$a=c$). Then $F_{a,b,c}$ has unique minimizer $(0,0)$, compact sublevels,
satisfies Corollary~\ref{art:planar-convergence} on every positive sublevel and
Corollary~\ref{art:planar-max} at $(0,0)$, admits an infinite admissible run
from every differentiable $x_0$, every $H_0\succ0$, and every $0<c_1<c_2<1$, and
every such run converges in value and iterates to $(0,0)$.
\end{proposition}
\begin{proof}
Both branches have Hessian at least $2\min\{a,b,c\}I\succ0$ and average to
$\tfrac12(a+c)u^2+bv^2$, giving unique minimizer $(0,0)$, compact sublevels, and
strictly convex positive contours. Along the tie curve $v=\tfrac12(c-a)u^2$, the
two branch gradients $(2au,2bv+1)^T$ and $(2cu,2bv-1)^T$ have determinant
$-2u(a+c+b(c-a)^2u^2)$, which vanishes only at $(0,0)$; for $u\ne0$ they are
linearly independent and their affine hull misses $0$. Thus
$\mathcal C_{\rm aff}=\{(0,0)\}$ and Corollary~\ref{art:planar-convergence}
applies on every positive sublevel, while at $(0,0)$ the active gradients
$(0,1)^T$ and $(0,-1)^T$ are affinely independent with weights $1/2$. Finally,
$(0,0)$ is nondifferentiable (so every differentiable point is nonstationary),
and along any line through a differentiable point the branch difference
$(a-c)u^2+2v$ is a nonzero polynomial of degree at most two with at most two
roots, so Lemma~\ref{art:examples-run-existence} guarantees run existence from
every differentiable start.\qed
\end{proof}

\begin{corollary}[Run existence on the three-quadratic nonradial example]\label{art:examples-nonradial-runs}
For the objective $F$ of Theorem~\ref{art:nonradial-example}, every
differentiable $x_0$ with $F(x_0)<1/100$ (for instance $x_0=(1/1000,0)^T$, where
$F(x_0)=1001/10^6<1/100$), every $H_0\succ0$, and every $0<c_1<c_2<1$ admit an
infinite admissible run, and every such run converges to $(0,0)$.
\end{corollary}
\begin{proof}
Here $F(u,v)\ge u^2+v^2$ has unique minimizer $(0,0)$, where $F$ is
nondifferentiable (since $F-f_0\ge0$ vanishes at $0$ with
$\nabla f_0(0)=(1,0)^T\ne0$, whereas $F\ge u^2+v^2$ forces any gradient at $0$
to vanish), so $\nabla F(x)\ne0$ at every differentiable point. Completing
squares gives the pairwise tie loci
\[
 \begin{aligned}
 f_0=f_1&\iff \tfrac12(u+1)^2-\bigl(v+\tfrac12\bigr)^2=\tfrac14,\\
 f_0=f_2&\iff (u-2)^2+(v-1)^2=5,\\
 f_1=f_2&\iff \tfrac12(v+2)^2-\bigl(u-\tfrac12\bigr)^2=\tfrac74.
 \end{aligned}
\]
Two nondegenerate hyperbolas and a circle contain no affine line, so along any
line $x+tp$ ($p\ne0$) each $f_i(x+tp)-f_j(x+tp)$ ($i<j$) is a nonzero polynomial
in $t$ of degree at most two, with at most two roots. Lemma~\ref{art:examples-run-existence}
and Theorem~\ref{art:nonradial-example} complete the proof.\qed
\end{proof}

\subsection{Why theorems quantify over admissible runs: an obstructed initial ray}

\begin{proposition}[A coercive homogeneous objective with no differentiable first step]\label{art:examples-no-step}
The function $G(u,v)=\max\{u,-u+v,-u-v\}$ is finite, convex, coercive, and
positively homogeneous of degree one with unique minimizer $(0,0)$ (hence covered
by Theorem~\ref{art:homogeneous}). At any $x_0=(\varepsilon,0)^T$
($\varepsilon>0$), $G$ is differentiable; with $H_0=I$, however, the search ray
$x_0-\alpha H_0\nabla G(x_0)$ ($\alpha>0$) contains no differentiable
point satisfying weak Wolfe for any $c_2<1$.
\end{proposition}
\begin{proof}
Since $G(u,v)\ge\max\{u,-u\}=|u|$ and $2G(u,v)\ge u+(-u\pm v)=\pm v$ gives
$2G(u,v)\ge|v|$, $G$ is coercive with unique minimizer $(0,0)$. At
$x_0=(\varepsilon,0)^T$, $\nabla G(x_0)=(1,0)^T$, so the search ray is
$(\varepsilon-\alpha,0)^T$ ($\alpha>0$). For $0<\alpha<\varepsilon$,
$\nabla G(\varepsilon-\alpha,0)=(1,0)^T$, so weak Wolfe would require
$-1\ge-c_2$, impossible for $c_2<1$; for $\alpha\ge\varepsilon$, the two branches $-u+v$ and $-u-v$ both evaluate to
$\alpha-\varepsilon\ge\varepsilon-\alpha$ with distinct gradients $(-1,1)^T$ and
$(-1,-1)^T$, so $(\varepsilon-\alpha,0)^T$ is nondifferentiable.\qed
\end{proof}

\addtocontents{toc}{\protect\vspace{6pt}\noindent{\protect\footnotesize\protect\sffamily\protect\bfseries\protect\color{CiteTeal} APPENDICES FOR PART II}\protect\par\protect\vspace{2pt}}
\section{The algebraic seed and its finite rational certificate}
\label{art:seed_certificate}

This appendix proves Lemma~\ref{art:seed_cycle} in three purely rational stages:
isolating the algebraic root $(c_\circ,d_\circ)$, verifying symbolic closure of
the three-step matrix and iterate cycle, and enclosing all nine trial queries in
exact rational interval arithmetic. All decimals below denote exact terminating
rationals.

\subsection{Existence and uniqueness of the defining root}
We choose the third slope $(c_\circ,d_\circ)$ so that the last two accepted step
lengths in the symbolic cycle below satisfy $\beta(c_\circ,d_\circ)=1/4$ and
$\gamma(c_\circ,d_\circ)=1/16$. Clearing positive denominators yields the polynomial system
\begin{equation}\label{art:seed_polynomials}
 P(c_\circ,d_\circ)=
 \begin{pmatrix}264192RK-UV(U+129)\\16VTK-RUY\end{pmatrix}=0,
\end{equation}
where
\begin{align*}
 U&=127c_\circ+16d_\circ,&V&=16c_\circ-127d_\circ-16,&R&=c_\circ-8d_\circ,\\
 K&=8c_\circ^2+8d_\circ^2+d_\circ-8,&T&=c_\circ^2-c_\circ+d_\circ^2,&Y&=(c_\circ-1)^2+d_\circ^2.
\end{align*}
We isolate the root inside the closed $\ell_\infty$ box
$\|(c_\circ,d_\circ)-\bar r\|_\infty\le10^{-30}$ centered at
\begin{align*}
 \bar r_1&=0.31848142762711863064645131250290703717912205023605,\\
 \bar r_2&=-1.02523211967574872548135087001300414284511909318444.
\end{align*}

\begin{lemma}[Rational root certificate]\label{art:root_certificate}
The closed box $\|(c_\circ,d_\circ)-\bar r\|_\infty\le10^{-30}$ contains a unique
root of $P$. Throughout this box all of $U,V,R,K,T,Y,U+129$ are strictly
positive, and~\eqref{art:seed_rho} holds.
\end{lemma}
\begin{proof}
Set $D=DP(\bar r)$ and let $A=\operatorname{adj}(D)/\det D$ be its $2\times2$
rational inverse. Direct rational evaluation of the eighteen Taylor monomials of
$P$ of total degrees $2\le|\nu|\le5$ at $\bar r$ certifies the three exact
inequalities
\begin{equation}\label{art:root_arithmetic}
 \begin{gathered}
 \det D\ne0,\qquad
 \max_i |(AP(\bar r))_i|<10^{-50},\\
 \max_i\sum_{2\le|\nu|\le5}|a_{i,\nu}|\,|\nu|\le7000,
 \qquad
 a_{i,\nu}=-\frac{1}{\nu_1!\nu_2!}\sum_{j=1}^2A_{ij}\,\partial^\nu P_j(\bar r).
 \end{gathered}
\end{equation}
Because $AD=I$, the preconditioned Newton map $N(h)=h-AP(\bar r+h)$ has zero
linear part and exact polynomial expansion
$N_i(h)=-(AP(\bar r))_i+\sum_{2\le|\nu|\le5}a_{i,\nu}h^\nu$. On the closed ball
$\|h\|_\infty,\|k\|_\infty\le r=10^{-30}$, telescoping each monomial gives
$|h^\nu-k^\nu|\le|\nu|r\|h-k\|_\infty$ ($|\nu|\ge2$), so $N$ is a contraction
with Lipschitz constant at most $7000r\ll1$ and maps the ball into itself since
$\|N(h)\|_\infty\le10^{-50}+7000r^2<r$ (cf.\ Rump~\cite[Theorem~13.3]{NEG_RUMP10}).
Invertibility of $A$ identifies its unique fixed point with the unique root of
$P$ in the box.

Since the box lies inside $318/1000<c_\circ<319/1000$ and
$-1026/1000<d_\circ<-1025/1000$, interval evaluation gives strict positivity of
$U,V,R,K,T,Y,U+129$; for instance,
\[
 \frac{2397}{100}\le U\le\frac{24113}{1000},\quad
 \frac{119263}{1000}\le V\le\frac{59703}{500},\quad
 \frac{4259}{500}\le R\le\frac{8527}{1000}.
\]
Substituting these positive bounds into $\rho=129R/[V(U+129)]$ yields
$60/1000<\rho<61/1000$. In particular, all divisions used below are valid.\qed
\end{proof}

\subsection{Symbolic closure of the three-step matrix cycle}
\begin{proof}[Proof of the algebraic part of Lemma~\ref{art:seed_cycle}]
From $s_0=(-1,0)^T$ and $g_1-g_0=(-2,-1/8)^T$, the first update is immediate:
$\mathcal B(I,g_0,g_1,1)=M_1$. For any $(c_\circ,d_\circ)$ where the factors
above are positive, define the rational step lengths and contraction ratio
\[
 \beta=\frac{66048RK}{UV(U+129)},\qquad
 \gamma=\frac{VTK}{RUY},\qquad
 \widetilde\rho=\frac{2064TK}{U(U+129)Y}.
\]
We verify the exact symbolic identities
\begin{equation}\label{art:closure_identity}
 \begin{gathered}
 \mathcal B(\mathcal B(M_1,g_1,g_2,\beta),g_2,g_0,\gamma)
     =\widetilde\rho I,\\
 -\gamma\mathcal B(M_1,g_1,g_2,\beta)g_2
     =\widetilde\rho(g_0-g_2).
 \end{gathered}
\end{equation}
Structurally, in~\eqref{art:failure_update} the projection factors of
$M_2(\beta)=\mathcal B(M_1,g_1,g_2,\beta)$ are independent of $\beta$ while its
rank-one term is linear in $\beta$, so $M_2(\beta)=C+\beta E$. Requiring
$(C+\beta E)g_2$ to be parallel to $g_0-g_2$ is a single linear equation in
$\beta$ solved by the displayed formula; for a unit vector $w\perp(g_0-g_2)$,
rational evaluation gives $w^T(C+\beta E)w=\widetilde\rho$, and the displayed
$\gamma$ gives last step $s_2=-\gamma M_2(\beta)g_2=\widetilde\rho(g_0-g_2)$. At
the third update, $y_2=g_0-g_2$ is parallel to $s_2$, so the projection matrix
$I-s_2y_2^T/(y_2^Ts_2)$ is the orthogonal projector $\Pi=ww^T$ onto
$(g_0-g_2)^\perp$: its projected term equals $\Pi(C+\beta E)\Pi=\widetilde\rho\,\Pi$,
its rank-one term equals $s_2s_2^T/(y_2^Ts_2)=\widetilde\rho(I-\Pi)$, and their
sum is $\widetilde\rho I$, proving~\eqref{art:closure_identity}.

At the root of~\eqref{art:seed_polynomials}, $\beta=1/4$, $\gamma=1/16$, and
\[
 \widetilde\rho-\rho=
 \frac{2064R}{V(U+129)}\left(\gamma-\frac1{16}\right)=0.
\]
Direct substitution for $j=0,1$ and the second identity in~\eqref{art:closure_identity}
for $j=2$ confirm that $s_j=-\alpha_jM_jg_j$ in~\eqref{art:seed_steps} with
$(\alpha_0,\alpha_1,\alpha_2)=(1,1/4,1/16)$, while~\eqref{art:seed_points} solves
$s_0+s_1+s_2=-(1-\rho)z_0$, giving $z_3=\rho z_0$. Finally, the strict margins
certified below imply $(g_{j+1}-g_j)^Ts_j>0$ for $j=0,1,2$, so $M_0,M_1,M_2\succ0$
by Lemma~\ref{art:run-invariants}.\qed
\end{proof}

\subsection{Rational interval certification of all nine queries}
Across one cycle, the unit-first search evaluates exactly nine trial points. We
bound them by exact rational interval arithmetic, using the standard rules
$[a,b]+[c,d]=[a+c,b+d]$, $-[a,b]=[-b,-a]$,
$[a,b][c,d]=[\min(ac,ad,bc,bd),\max(ac,ad,bc,bd)]$, and
$[a,b]^{-1}=[1/b,1/a]$ when $0\notin[a,b]$. Starting from the coarse rational
enclosure
\begin{equation}\label{art:query_box}
 \begin{gathered}
 c_\circ\in\frac{[318481,318482]}{10^6},\qquad
 d_\circ\in\frac{[-1025233,-1025232]}{10^6},\\
 \rho\in\frac{[6019,6021]}{10^5}
 \end{gathered}
\end{equation}
implied by Lemma~\ref{art:root_certificate} and~\eqref{art:seed_rho}, we
evaluate $z_j$ from~\eqref{art:seed_points} and the search directions
\[
 p_0=(-1,0)^T,\qquad p_1=(127/256,1/16)^T,\qquad
 p_2=16\rho(1-c_\circ,-d_\circ)^T.
\]
For each row of Table~\ref{art:trial_table} with designated active branch $i$
and trial point $q=z_j+\alpha p_j$, we bound from below the active-branch,
slope, and acceptance/rejection margins
\[
 (g_i-g_k)^Tq\quad(k\ne i),\qquad -\alpha g_j^Tp_j,\qquad
 \begin{cases}
 g_j^Tz_j+\tfrac13\alpha g_j^Tp_j-g_i^Tq,&\text{acceptance},\\
 g_i^Tq-g_j^Tz_j,&\text{rejection},
 \end{cases}
\]
together with $\alpha g_i^Tp_j$ on accepted rows. Interval evaluation shows that
every margin in each row exceeds the corresponding entry of the positive
rational vector
\begin{equation}\label{art:query_lower_bounds}
 \frac1{10^4}(729,\;2417,\;145,\;75,\;5964,\;2682,\;1041,\;221,\;20).
\end{equation}
The same interval operations give the value enclosures in Table~\ref{art:trial_table},
positive active-branch gaps at $z_0,z_1,z_2$ (hence local differentiability at
every base point and query), and positivity of $3\pm g_k^Tq$ ($k=0,1,2$),
$2-[q]_2$, and $G(q)=g_i^Tq$ (for the certified active branch $i$) across all
base points and queries, establishing~\eqref{art:seed_query_bounds}
and completing the proof of Lemma~\ref{art:seed_cycle}.

\section{Scalar spline interpolations and the strongly convex local embedding}
\label{art:lower_curvature}

This appendix contains the two scalar spline constructions and the
perturbed-orbit persistence proof used in Section~\ref{art:failure-local}.
Subsection~\ref{art:upper_spline_subsec} proves
Lemma~\ref{art:upper_scalar_lemma}, giving the $C^{2,1}$ plateau spline for the
exact-copy construction of Theorem~\ref{art:local-failure}. Under uniform strong
convexity (Theorem~\ref{art:local-failure-lower}), exact scale invariance is
lost: the second subsection therefore constructs a single perturbed failing
orbit via a scale-normalized BFGS return map and its local stable manifold, and
the third subsection copies that orbit across scales using a $C^{1,1}$ spline.
Throughout, an \emph{admissible pair} means any
fixed $(c_1,c_2)$ with $0<c_1\le1/3$ and $c_1<c_2<1$; the objectives, orbits,
and starting points constructed below do not depend on $(c_1,c_2)$.

\subsection{Proof of Lemma~\ref{art:upper_scalar_lemma}: the \texorpdfstring{$C^{2,1}$}{C\^{}\{2,1\}} plateau spline}\label{art:upper_spline_subsec}

To connect infinitely many affine plateaus at nodes $ba_\ell$ ($a_\ell=\theta^\ell$,
$\theta=1/16$, $b\in\{2,b_*\}$ with $b_*=2c_\circ-32d_\circ\in(33,34)$) while
keeping $\varphi_b''$ globally Lipschitz down to $0$, define for transition
half-width $\varepsilon>0$ and center $t_b\in(\varepsilon,1-\varepsilon)$ the
monotone $C^1$ cubic derivative ramp and its primitive
\begin{gather*}
 P_b(\xi)=
 \begin{cases}
 0,&\xi\le t_b-\varepsilon,\\
 3v^2-2v^3,&|\xi-t_b|<\varepsilon,\\
 1,&\xi\ge t_b+\varepsilon,
 \end{cases}\\
 v=\frac{\xi-t_b+\varepsilon}{2\varepsilon},
 \qquad J_b(\xi)=\int_0^\xi P_b(t)\,dt,
\end{gather*}
so $J_b(\xi)=0$ for $\xi\le t_b-\varepsilon$, $J_b(\xi)=2\varepsilon(v^3-v^4/2)$
for $|\xi-t_b|<\varepsilon$, and $J_b(\xi)=\xi-t_b$ for $\xi\ge t_b+\varepsilon$.

\begin{proof}[Proof of Lemma~\ref{art:upper_scalar_lemma}]
Take $\theta=1/16$, $a_\ell=\theta^\ell$, $\varepsilon=1/128$, and
$t_b=[1-(1+\theta+\theta^2)/b]/(1-\theta^2)$. On each dyadic interval
$b\theta a_\ell\le x\le ba_\ell$ ($\ell\ge0$), define
\begin{equation}\label{art:upper_scalar}
 \varphi_b(x)=(\theta a_\ell)^3+
 ba_\ell^3(1-\theta)\bigl\{\theta^2\xi+(1-\theta^2)J_b(\xi)\bigr\},\qquad
 \xi=\frac{x/(ba_\ell)-\theta}{1-\theta},
\end{equation}
and extend $\varphi_b$ outside $(0,b]$ by
\[
 \varphi_b(x)=
 \begin{cases}
 1-3(x+1)+3(x+1)^2,&x<-1,\\
 -x^3,&-1\le x\le0,\\
 1+x-b,&x\ge b.
 \end{cases}
\]
Since $P_b'=6v(1-v)/(2\varepsilon)$ on $|\xi-t_b|<\varepsilon$ and $P_b'=0$
elsewhere, $0\le P_b'\le3/(4\varepsilon)=96$; with $\int_0^1P_b=1-t_b$,
differentiation on $[b\theta a_\ell,ba_\ell]$ gives
\[
 \varphi_b'(x)=a_\ell^2\bigl\{\theta^2+(1-\theta^2)P_b(\xi)\bigr\},\qquad
 \varphi_b''(x)=\frac{a_\ell(1+\theta)}bP_b'(\xi).
\]
Because $\theta^2+(1-\theta^2)(1-t_b)=(1+\theta+\theta^2)/b$, we have
$\varphi_b(ba_\ell)=a_\ell^3$, while $P_b$ is constant ($0$ and $1$) near
$\xi=0$ and $\xi=1$, so $\varphi_b'$ and $\varphi_b''$ match across every node
$ba_\ell$ with values $a_\ell^2$ and $0$. As $x\downarrow0$, $a_\ell\le x/(b\theta)$
gives $\varphi_b(x)=O(x^3)$, $\varphi_b'(x)=O(x^2)$, and $\varphi_b''(x)=O(x)$,
matching $-x^3$ through second order at $0$; values through second order also
match at $x=-1$ and $x=b$. Since $\varphi_b''\ge0$ everywhere ($6$ on
$(-\infty,-1)$, $-6x$ on $[-1,0]$, $(a_\ell(1+\theta)/b)P_b'\ge0$ on $(0,b]$,
and $0$ on $[b,\infty)$) with $\varphi_b'<0$ on $(-\infty,0)$ and $\varphi_b'>0$
on $(0,\infty)$, $\varphi_b$ is convex, coercive, nonnegative, and zero only at
$0$, with $\varphi_b''\le\max\{6,(17/16)\cdot96/b\}=\max\{6,102/b\}$.

Around $ba_\ell$, the flat-slope intervals below and above $ba_\ell$ have
half-widths (divided by $a_\ell$) equal to $b(1-\theta)(1-t_b-\varepsilon)$ and
$b(1-\theta)(t_b-\varepsilon)/\theta$, both strictly exceeding $1/2$ for all
$b\in[2,34]$ (with $[b,\infty)$ covering the upper side at $\ell=0$),
proving~\eqref{art:scalar_plateau}. Finally, $|P_b''|\le3/(2\varepsilon^2)$ on
each transition gives $|\varphi_b'''|\le(17/15)\cdot3/(2\varepsilon^2)b^{-2}=139264/(5b^2)$
on $(0,b)$ and $|\varphi_b'''|\le6$ on $(-\infty,0)$; continuity of $\varphi_b''$
across joins and at $0$ yields the global Lipschitz constant
$L_b=\max\{6,139264/(5b^2)\}$.\qed
\end{proof}

\subsection{A strongly convex prototype and hyperbolic return map}

On the three-slope cone $G(x)=\max_{0\le i\le2}g_i^Tx$ of
Section~\ref{art:failure}, consider the $\mu$-strongly convex objective
\[
 Q_\mu(x)=\tfrac\mu2\norm{x}^2+G(x)\qquad(\mu>0).
\]
The minimum-norm vector in $\conv\{g_0,g_1,g_2\}$ is
\[
 \bar g=\frac{(1,-16)^T}{257}=\frac{129}{257}g_0+\frac{128}{257}g_1,
\]
since $g_0^T\bar g=g_1^T\bar g=\norm{\bar g}^2=1/257<g_2^T\bar g$. Consequently
$Q_\mu$ has unique nonsmooth minimizer and minimum value
\begin{equation}\label{art:quadratic_minimizer}
 x_*=-\frac{\bar g}{\mu},\qquad \min Q_\mu=-\frac{1}{514\mu}<0,
\end{equation}
whereas at the origin $\partial Q_\mu(0)=\conv\{g_0,g_1,g_2\}$ remains
nonstationary with $Q_\mu(0)=0$.

\begin{proposition}[Perturbed failing trajectory under strong convexity]\label{art:quadratic_failure}
Let $\mu>0$. For every neighborhood $V$ of $0\in\R^2$, there exists an initial
state $(x_0,H_0)$ with $H_0\succ0$ such that the infinite deterministic
full-BFGS run on $Q_\mu$ from $(x_0,H_0)$ is independent of $(c_1,c_2)$ across
all admissible pairs $0<c_1\le1/3$, $c_1<c_2<1$, converging to the
nonstationary origin with $Q_\mu(x_k)\downarrow0$,
accepted step lengths repeating $1,1/4,1/16$, and $x_0$ together with all
accepted and rejected trial points differentiable and contained in $V$.
\end{proposition}

To capture perturbed trajectories as scale shrinks, write states in normalized
coordinates
\[
 x=\tau z,\qquad H=\tau M(\zeta,\eta),\qquad
 M(\zeta,\eta)=\begin{pmatrix}1&\zeta\\\zeta&\eta\end{pmatrix}\qquad(\tau>0),
\]
so the homogeneous seed corresponds to the boundary point $P_*=(z_0,0,1,0)$ in
$(z,\zeta,\eta,\tau)\in\R^2\times\R\times\R\times\R$. Holding the input scale
$\tau$ fixed across three substeps $j=0,1,2$ with $\alpha=(1,1/4,1/16)$ and
cyclic indices $g_3=g_0$, initialize $(Z,M)=(z,M(\zeta,\eta))$ and iterate
\begin{align}
 \gamma_j&=g_j+\mu\tau Z,\qquad S=-\alpha_jM\gamma_j,
       \qquad Z^+=Z+S,\label{art:lower_return_substep}\\
 Y&=g_{j+1}-g_j+\mu\tau S,\qquad
 M^+=\left(I-\frac{SY^T}{Y^TS}\right)M
          \left(I-\frac{YS^T}{Y^TS}\right)+\frac{SS^T}{Y^TS}.\notag
\end{align}
Writing $(Z_3,M_3)$ for the result after $j=0,1,2$ and $\hat\rho=(M_3)_{11}$,
define the five-dimensional normalized return map
\begin{equation}\label{art:lower_return}
 \mathcal T_\mu(z,\zeta,\eta,\tau)
 =\bigl(Z_3/\hat\rho,\;(M_3)_{12}/\hat\rho,\;(M_3)_{22}/\hat\rho,\;\tau\hat\rho\bigr).
\end{equation}
By Lemma~\ref{art:seed_cycle}, $\mathcal T_\mu(P_*)=P_*$ with $\hat\rho(P_*)=\rho>0$
and nonvanishing denominators, so $\mathcal T_\mu$ extends real-analytically to
an open neighborhood of $P_*$ (allowing $\tau\le0$ analytically).

\begin{lemma}[Exact search realization and hyperbolicity]\label{art:lower_hyperbolic}
On a sufficiently small neighborhood of $P_*$, every state with $\tau>0$
executes the exact nine-query decision lists of Table~\ref{art:trial_table}
under the unit-first bracketing algorithm for every admissible pair $(c_1,c_2)$,
and~\eqref{art:lower_return} is its three-step return map. Moreover,
$D\mathcal T_\mu(P_*)$ has five real eigenvalues
\[
 \rho^{-1},\quad\rho^{-1},\quad\lambda_s,\quad\lambda_u,\quad\rho,
 \qquad 0.038<\lambda_s<0.039,\quad 1.564<\lambda_u<1.565,
\]
and its two-dimensional stable eigenspace contains an eigenvector with
$\tau$-component equal to $1$.
\end{lemma}
\begin{proof}
At any physical query $\tau Z$ ($\tau>0$), $Q_\mu(\tau Z)/\tau=G(Z)+\tfrac\mu2\tau\norm{Z}^2$
(adding the same $\tfrac\mu2\tau\norm{Z}^2$ to every branch, leaving active-branch
gaps unchanged) and $\nabla Q_\mu(\tau Z)=g_i+\mu\tau Z$. Because all active
branch gaps, rejection value increases, and acceptance margins~\eqref{art:seed_margins}
in Lemma~\ref{art:seed_cycle} are strict at $P_*$ and continuous in
$(z,\zeta,\eta,\tau)$, they persist on one neighborhood of $P_*$ simultaneously
for all $0<c_1\le1/3$, $c_1<c_2<1$; multiplying~\eqref{art:lower_return_substep}
back by $\tau$ recovers the exact BFGS steps and matrix updates.

At $\tau=0$, the matrix substeps are independent of $z$. Let $J$ denote the
$2\times2$ Jacobian at $(\zeta,\eta)=(0,1)$ of the normalized shape map
$(\zeta,\eta)\mapsto((M_3)_{12},(M_3)_{22})/\hat\rho$ at $\tau=0$. Differentiating
the rational map~\eqref{art:lower_return_substep} symbolically at $\tau=0$ and
applying the rational interval arithmetic of Appendix~\ref{art:seed_certificate}
on the tighter enclosure
\begin{align*}
 c_\circ&\in[318481427627118,318481427627119]/10^{15},\\
 d_\circ&\in[-1025232119675749,-1025232119675748]/10^{15}
\end{align*}
(which contains $\|(c_\circ,d_\circ)-\bar r\|_\infty\le10^{-30}$) yields
\begin{equation}\label{art:shape_certificate}
 10^{10}J\in
 \begin{pmatrix}
 [11123986249,11123986253]&[6459701025,6459701029]\\
 [7515015701,7515015708]&[4905136882,4905136890]
 \end{pmatrix}.
\end{equation}
Evaluating $p(\lambda)=\lambda^2-(\tr J)\lambda+\det J$ via~\eqref{art:shape_certificate}
certifies $p(0.038)>0$, $p(0.039)<0$, $p(1.564)<0$, and $p(1.565)>0$, isolating
the real eigenvalues $\lambda_s\in(0.038,0.039)$ and $\lambda_u\in(1.564,1.565)$.

Finally, at $\tau=0$ we have $Z_3=z+v(\zeta,\eta)$ with $\hat\rho$ independent of
$z$, and the fifth component of $\mathcal T_\mu$ is $\tau\hat\rho$, giving the
block-upper-triangular Jacobian
\[
 D\mathcal T_\mu(P_*)=
 \begin{pmatrix}
 \rho^{-1}I_2&*&*\\0&J&*\\0&0&\rho
 \end{pmatrix}.
\]
Its spectrum is $\{\rho^{-1},\rho^{-1},\lambda_s,\lambda_u,\rho\}$, disjoint from
$0$ and the unit circle. Because $\rho\in(0.060,0.061)$ does not lie in the
spectrum $\{\rho^{-1},\lambda_s,\lambda_u\}$ of the upper-left $4\times4$ block,
solving $(D\mathcal T_\mu(P_*)-\rho I_5)v=0$ with fifth coordinate $v_5=1$
produces a stable $\rho$-eigenvector with $\tau$-component $1$.\qed
\end{proof}

\begin{proof}[Proof of Proposition~\ref{art:quadratic_failure}]
Since $\det D\mathcal T_\mu(P_*)\ne0$ and $P_*$ is hyperbolic, $\mathcal T_\mu$
is a local real-analytic diffeomorphism admitting a smooth two-dimensional local
stable manifold through $P_*$ tangent to its stable eigenspace (see, e.g.,
Dyatlov~\cite[Section~3.5, Theorem~3 and Remark~4]{DY}). Because each
of the finitely many physical query maps $\Phi_q(w)=\tau q(w)$ satisfies
$\Phi_q(P_*)=0$, shrink the implementation neighborhood of
Lemma~\ref{art:lower_hyperbolic} inside $\bigcap_q\Phi_q^{-1}(V)$ so that
$\hat\rho>0$ and $M(\zeta,\eta)\succ0$ throughout. Because the stable
$\rho$-eigenvector has $\tau$-coordinate $1$, this local stable manifold is not
confined to the boundary $\tau=0$: it penetrates into positive physical scale
and contains an initial state $w_0$ with $\tau_0>0$ whose forward iterates
$w_m=(\hat z_m,\zeta_m,\eta_m,\tau_m)=\mathcal T_\mu^m(w_0)$ stay in that
neighborhood and converge to $P_*$ with $\tau_m>0$ and $\tau_m\to0$. Unfolding
the substeps yields a single deterministic full-BFGS run valid for every
admissible pair $(c_1,c_2)$, with all queries inside $V$, $x_k\to0$, and
$Q_\mu(x_k)\downarrow0>\min Q_\mu$.\qed
\end{proof}

\subsection{The \texorpdfstring{$C^{1,1}$}{C\^{}\{1,1\}} orbit-copying spline and proof of Theorem~\ref{art:local-failure-lower}}

To copy a perturbed trajectory of $Q_\mu$ at centers $\bar x_\ell=a_\ell c_*$
approaching $0$ while preserving a common initial matrix ($\sigma=\lambda^2$),
spatial scaling $\lambda=a_\ell$ requires objective scaling $\sigma=a_\ell^2$,
hence scalar spline nodes of value $a_\ell^2$ and slope $a_\ell$. We shift the
linear branch directions by $-\mu c_*$ so that expanding $\tfrac\mu2\norm{\bar x_\ell+a_\ell y}^2$
restores the exact slopes $g_i$ inside every copy.

Specifically, fix
\begin{gather*}
 \mu=\frac1{4096},\qquad c_*=(2,-32)^T,\qquad q_i=g_i-\mu c_*,\\
 b_i=q_i^Tc_*,\qquad r=\frac1{128},\qquad a_\ell=r^\ell,
\end{gather*}
so $b_0=b_1=1791/1024$ and $1.7<b_i<34$. For $b\in\{b_0,b_2\}$, reuse the cubic
ramp $P_b$ and primitive $J_b$ from Subsection~\ref{art:upper_spline_subsec} with
$t_b=[1-(1+r)/b]/(1-r)$ and $\varepsilon=1/1024$, and define on each interval
$ba_{\ell+1}\le x\le ba_\ell$ ($\ell\ge0$)
\begin{equation}\label{art:lower_scalar}
 \psi_b(x)=(ra_\ell)^2+ba_\ell^2(1-r)
       \bigl\{r\xi+(1-r)J_b(\xi)\bigr\},\qquad
 \xi=\frac{x/(ba_\ell)-r}{1-r},
\end{equation}
with $\psi_b(0)=0$, $\psi_b(x)=x^2/2$ for $x<0$, and $\psi_b(x)=1+x-b$ for
$x\ge b$. On $[ba_{\ell+1},ba_\ell]$, we have
$\psi_b'(x)=a_\ell\{r+(1-r)P_b(\xi)\}$ and $\psi_b''(x)=P_b'(\xi)/b\in[0,768/b]$.
The identity $r+(1-r)(1-t_b)=(1+r)/b$ matches node values $\psi_b(ba_\ell)=a_\ell^2$
and slopes $\psi_b'(ba_\ell)=a_\ell$ across joins, and summing the uniform second-derivative
bound across intervals down to $0$ shows that $\psi_b\in C^{1,1}(\R)$ is convex,
coercive, nonnegative, and zero only at $0$. The plateau half-widths divided by
$a_\ell$ are $b(1-r)(1-t_b-\varepsilon)$ and $b(1-r)(t_b-\varepsilon)/r$, both
exceeding $1/2$ for $b\in(1.7,34)$, which yields
\begin{equation}\label{art:lower_plateau}
 \psi_b(ba_\ell+h)=a_\ell^2+a_\ell h\qquad(\ell\ge0,\;|h|\le a_\ell/2).
\end{equation}
Along $x_\ell=ba_\ell\to0^+$ and $\tilde x_\ell=(b+1/4)a_\ell\to0^+$, however,
$\psi_b'(x_\ell)/x_\ell=1/b$ whereas $\psi_b'(\tilde x_\ell)/\tilde x_\ell=1/(b+1/4)$,
so $\psi_b$ is not twice differentiable at $0$.

\begin{proof}[Proof of Theorem~\ref{art:local-failure-lower}]
First define the pre-whitened four-branch objective on $\R^2$ by
\begin{equation}\label{art:lower_objective}
 \begin{aligned}
 f_i(x)&=\tfrac\mu2\norm{x}^2+\psi_{b_i}(q_i^Tx)\qquad(i=0,1,2),\\
 f_3(x)&=\tfrac\mu2\norm{x}^2+x_2,\qquad F(x)=\max_{0\le i\le3}f_i(x).
 \end{aligned}
\end{equation}
Each $f_i$ is $\mu$-strongly convex and globally $C^{1,1}$, with gradient
Lipschitz constants bounded by $L_i=\mu+\max\{1,768/b_i\}\norm{q_i}^2$ ($L_0<439$,
$L_1<446$, $L_2<27$, and $L_3=\mu$), so $512$ bounds all four. Restricting $f_i$
($i<3$) to $tv_i$ with $q_i^Tv_i=1$ and subtracting $\tfrac\mu2t^2\norm{v_i}^2$
recovers $\psi_{b_i}(t)$, so $f_0,f_1,f_2$ are not twice differentiable at $0$.
Since $F(x)\ge\tfrac\mu2\norm{x}^2$ and $F(0)=0$, $0$ is the unique minimizer
(nonsmooth since $\nabla f_i(0)=0$ for $i<3$ while $\nabla f_3(0)=(0,1)^T$) and
all sublevels are compact. At any point $x$ where $F$ is twice differentiable,
all active gradients agree with $\nabla F(x)$ and inactive branches stay
strictly inactive nearby, so taking the maximum of the $512$-smooth upper Taylor
bounds on active branches gives $\mu I\preceq\nabla^2F(x)\preceq512I$.

Setting $\bar x_\ell=a_\ell c_*$ and $C_*=1+\tfrac\mu2\norm{c_*}^2$, for every
$y$ in the open neighborhood $V=\{y:\norm{y}_\infty<1,\;|q_i^Ty|<1/2\;(i=0,1,2)\}$,
\eqref{art:lower_plateau} and $q_i+\mu c_*=g_i$ give
$f_i(\bar x_\ell+a_\ell y)=a_\ell^2(C_*+\tfrac\mu2\norm{y}^2+g_i^Ty)\ge\tfrac12a_\ell^2$
for $i=0,1,2$, while $\norm{y}_\infty<1$ gives $\tfrac\mu2\norm{c_*+y}^2<1098/8192<1$
and hence $f_3(\bar x_\ell+a_\ell y)=\tfrac\mu2a_\ell^2\norm{c_*+y}^2+a_\ell(-32+y_2)<a_\ell^2-31a_\ell<0$
for all $\ell\ge1$ since $a_\ell\le1/128$. Hence
\begin{equation}\label{art:lower_local_model}
 F(\bar x_\ell+a_\ell y)=a_\ell^2\bigl(C_*+Q_\mu(y)\bigr)
 \qquad(\ell\ge1,\;y\in V).
\end{equation}
By Proposition~\ref{art:quadratic_failure}, choose a single perturbed failing
trajectory $(x_k^{\rm seed},H_k^{\rm seed})$ on $Q_\mu$ whose initial point and
all accepted and rejected queries lie in $V$ for every admissible pair
$(c_1,c_2)$. Starting at $\bar x_\ell+a_\ell x_0^{\rm seed}$ with initial matrix
$H_0^{\rm seed}\succ0$, scale covariance~\eqref{art:matched-scaling} with
$\lambda=a_\ell$ and $\sigma=a_\ell^2=\lambda^2$ gives by induction
$x_k^{(\ell)}=\bar x_\ell+a_\ell x_k^{\rm seed}$ and $H_k^{(\ell)}=H_k^{\rm seed}$.
For each fixed $\ell\ge1$, the iterates $x_k^{(\ell)}$ therefore converge to the
nonstationary point $\bar x_\ell\ne0$ ($\partial F(\bar x_\ell)=a_\ell\conv\{g_0,g_1,g_2\}\not\ni0$)
with $F(x_k^{(\ell)})\downarrow a_\ell^2C_*>0$ and query norms bounded by
$a_\ell(\norm{c_*}+\sup_{\text{seed queries}}\norm{x^{\rm seed}})\to0$.

Finally, apply the single invertible linear change of coordinates
$A=(H_0^{\rm seed})^{1/2}\succ0$, replacing $F$ by $\widetilde F(z)=F(Az)$ and
each run by $\widetilde x_k^{(\ell)}=A^{-1}x_k^{(\ell)}$ and
$\widetilde H_k^{(\ell)}=A^{-1}H_k^{(\ell)}A^{-T}$. Because $\widetilde g=A^Tg$
and $\widetilde s=A^{-1}s$ leave all scalar products $g^Ts$, function values,
trial decisions, and inverse-BFGS updates~\eqref{art:failure_update} invariant
(cf.~\cite[Section~2]{LO}), every run $\ell\ge1$ now starts from
$\widetilde H_0^{(\ell)}=I$, while $\widetilde F$ inherits four strongly convex
$C^{1,1}$ branches (with the first three not twice differentiable at $0$) and
satisfies $\underline\mu I\preceq\nabla^2\widetilde F\preceq\overline L I$
wherever twice differentiable with $\underline\mu=\mu\lambda_{\min}(H_0^{\rm seed})$
and $\overline L=512\lambda_{\max}(H_0^{\rm seed})$.\qed
\end{proof}

\bookmarksetup{startatroot}

\end{document}